\documentclass[12pt,reqno]{amsart}
\usepackage{amsmath}
\usepackage{amssymb}
\usepackage{verbatim}
\usepackage{mathrsfs}
\usepackage[capitalise]{cleveref}
\usepackage[noadjust]{cite}

\usepackage[top=0.5in,bottom=0.45in,left=0.7in,right=0.7in]{geometry}

\newcommand{\C}{\mathbb{C}}    %complex number field
\newcommand{\N}{\mathbb{N}}    %natural numbers
\newcommand{\NN}{\mathbb{N}_0} %Natural Nonnegative numbers
\newcommand{\R}{\mathbb{R}}    %real number field
\newcommand{\Z}{\mathbb{Z}}    %integers

\newcommand{\dNN}{\mathbb{N}_0^d}
\newcommand{\dR}{\mathbb{R}^d}

\newcommand{\dZ}{\mathbb{Z}^d}

\newcommand{\sD}{\mathsf{D}}

\newcommand{\CH}[1]{\mathscr{C}^{#1}(\mathbb{R})} 
\newcommand{\dCH}[1]{\mathscr{C}^{#1}(\mathbb{R}^d)} \newcommand{\dLp}[1]{L_{#1}(\mathbb{R}^d)}
\newcommand{\cd}{\mathcal{R}}  %cascade/refinement operator C: use C or R.
\newcommand{\sm}{\operatorname{sm}}  %smoothness exponent
\newcommand{\sr}{\operatorname{sr}}  %sum rules
\newcommand{\lpm}{\operatorname{lpm}}  %linear-phase moments

\newcommand{\bp}{ \begin{proof} }
\newcommand{\ep}{ \end{proof} } %\hfill
\newcommand{\be}{ \begin{equation} }
\newcommand{\ee}{ \end{equation} }
\newcommand{\imply}{ \Longrightarrow }
\newcommand{\tp}{\mathsf{T}}  %transpose
\newcommand{\fs}{\operatorname{fsupp}}
\newcommand{\mspan}{\operatorname{span}}

\newcommand{\pp}{\mathsf{p}}
\newcommand{\pq}{\mathsf{q}}

\newcommand{\lrs}[3]{(l_{#1}(\mathbb{Z}))^{#2\times #3}}
\newcommand{\dlp}[1]{l_{#1}(\mathbb{Z}^d)}
\newcommand{\dlrs}[3]{(l_{#1}(\mathbb{Z}^d))^{#2\times #3}}

\newcommand{\PL}{\Pi}   %ring of all polynomials
\newcommand{\PR}{\mathscr{P}}   %use for polynomial rings
\newcommand{\PV}{\mathscr{V}}
\newcommand{\PB}{\mathscr{B}}

\newcommand{\sd}{\mathcal{S}}  %subdivision operator S

\newcommand{\wh}{\widehat}
\renewcommand{\le}{\leqslant}
\renewcommand{\ge}{\geqslant}
\newcommand{\bs}{\backslash}
\newcommand{\ol}{\overline}

\newcommand{\bo}{\mathscr{O}} %standard big O notation
\newcommand{\setsp}{\;:\;}     %set separator

\newcommand{\vgu}{\upsilon} %matching filter for matrix filters
\newcommand{\td}{\boldsymbol{\delta}}
\newcommand{\ind}{\Lambda} %\triangle

\newtheorem{lemma}{Lemma}

\newtheorem{theorem}[lemma]{Theorem}

\newtheorem{example}{Example}

\newtheorem{definition}{Definition}

\numberwithin{equation}{section}
\numberwithin{lemma}{section}
\begin{document}

\title{Multivariate Generalized Hermite Subdivision Schemes}

\author{Bin Han}

\address{Department of Mathematical and Statistical Sciences,
University of Alberta, Edmonton,\quad Alberta, Canada T6G 2G1.
\quad {\tt bhan@ualberta.ca}\quad {\tt http://www.ualberta.ca/$\sim$bhan}
}

\thanks{Contact information: Department of Mathematical and Statistical Sciences,
University of Alberta, Edmonton,\quad Alberta, Canada T6G 2G1.
\quad {\tt bhan@ualberta.ca}\quad {\tt http://www.ualberta.ca/$\sim$bhan}\quad Fax: 1-780-4926826}

\thanks{Research was supported in part by the Natural Sciences and Engineering Research Council of Canada (NSERC) under grant RGPIN-2019-04276.}

%\thanks{Contact information of corresponding author Bin Han: E-mail: bhan@ualberta.ca, Phone: 1-587-8828375, Fax: 1-780-4926826,  Web: http://www.ualberta.ca/$\sim$bhan}

\makeatletter \@addtoreset{equation}{section} \makeatother

\begin{abstract}
Due to properties such as interpolation, smoothness, and spline connections, Hermite subdivision schemes employ fast iterative algorithms for geometrically modeling curves/surfaces in CAGD and for building Hermite wavelets in numerical PDEs.
In this paper we introduce a notion of generalized Hermite (dyadic) subdivision schemes and then we characterize their convergence, smoothness and  underlying matrix masks with or without interpolation properties. We also introduce the notion of linear-phase moments for achieving the polynomial-interpolation property. For any given integer $m\in \N$, we constructively prove that there always exist convergent smooth generalized Hermite subdivision schemes with linear-phase moments such that their basis vector functions are spline functions in $\dCH{m}$ and have linearly independent integer shifts.
As byproducts, our results resolve convergence, smoothness and existence of Lagrange, Hermite, or Birkhoff subdivision schemes.
Even in dimension one our results significantly generalize and extend many known results on extensively studied univariate Hermite subdivision schemes. To illustrate the theoretical results in this paper, we provide examples of convergent generalized Hermite subdivision schemes with symmetric matrix masks having short support and smooth basis vector functions with or without interpolation property.
\end{abstract}


\maketitle

\pagenumbering{arabic}

\section{Introduction and Main Results}

Hermite interpolation approximates a target function through an approximating function which matches the target function value and all its consecutive derivatives to a prescribed order at given data points.
As a generalization, Birkhoff interpolation allows the derivatives for interpolation to be not necessarily consecutive. Being an iterative locally averaging algorithm for fast computing a limit function and its consecutive derivatives,
a Hermite subdivision scheme is of interest for generating subdivision curves/surfaces in CAGD (e.g., see \cite{dl02,dm09,hyx05,mer92}), for building wavelets to develop multiscale numerical solvers of PDEs (\cite{dah97,dhjk00,hm21}), and for
implementing the fast wavelet transform to process data (\cite{han09,hanbook}),
due to their properties such as interpolation, high smoothness and approximation, and spline connections.
From the literature review in Subsection~\ref{subsec:review}, we shall see that only the univariate Hermite subdivision schemes have been extensively studied so far and there are barely systematic study in the multivariate setting in the literature.
In this paper we shall generalize and extend Hermite subdivision schemes to the multivariate setting, and then systematically analyze their mathematical properties. Such generalized Hermite subdivision schemes not only include all known Hermite and Lagrange subdivision schemes in the literature as special cases but also introduce new types of subdivision schemes
such as Birkhoff subdivision schemes with various properties.

\subsection{Characterize masks of generalized Hermite subdivision schemes}

To introduce the notion of generalized Hermite subdivision schemes, we first recall some necessary definitions and notations.
By $\dlrs{}{s}{r}$ we denote the linear space of all sequences $u=\{u(k)\}_{k\in \dZ}: \dZ \rightarrow \C^{s\times r}$.
Similarly,
$\dlrs{0}{s}{r}$ consists of all finitely supported sequences $u\in \dlrs{}{s}{r}$ with $\{k\in \dZ \setsp u(k)\ne 0\}$ being finite.
%For $u\in \dlrs{0}{r}{s}$ and $v\in \dlrs{0}{s}{t}$, their convolution is defined to be $u*v:=\sum_{k\in \dZ} u(\cdot-k)v(k)$. In terms of symbols, we have $\wh{u*v}(\xi)=\wh{u}(\xi)\wh{v}(\xi)$.
A multivariate matrix mask or filter is often given by $a=\{a(k)\}_{k\in \dZ}\in \dlrs{0}{r}{r}$ with $a(k)\in \C^{r\times r}$ for all $k\in \dZ$.
To state a generalized Hermite (dyadic) subdivision scheme,
the \emph{vector/matrix subdivision operator} $\sd_a: \dlrs{}{s}{r} \rightarrow \dlrs{}{s}{r}$ is defined to be
\be \label{sd}
(\sd_a v)(j):=2^d\sum_{k\in \dZ} v(k) a(j-2k),\qquad j\in \dZ
\ee
for $v=\{v(k)\}_{k\in \dZ}\in \dlrs{}{s}{r}$.
%From \eqref{sd} we observe
%$\wh{\sd_a v}(\xi)=2^d \wh{v}(2\xi)\wh{a}(\xi)$ for $v\in \dlrs{0}{s}{r}$ and
We often use initial data $v\in  \dlrs{}{s}{r}$ with either $s=1$ for sequences $v$  of row vectors or $s=r$ for sequences $v$ of $r\times r$ square matrices in this paper.
Any convergent subdivision scheme with a matrix mask $a\in \dlrs{0}{r}{r}$ naturally has a basis vector function $\phi=[\phi_1,\ldots,\phi_r]^\tp$ on $\dR$ which is defined in \eqref{def:phi}; and we shall see in \cref{thm:hsd:mask} that such a basis vector function $\phi$ must be a refinable vector function  satisfying the refinement equation:
\be \label{refeq:phi}
\phi=2^d\sum_{k\in \dZ} a(k)\phi(2\cdot-k),
\quad \mbox{or equivalently},\quad
\wh{\phi}(2\xi)=\wh{a}(\xi)\wh{\phi}(\xi),
\ee
where the Fourier transform is defined to be $\wh{f}(\xi):=\int_{\dR} f(x) e^{-ix \cdot \xi} dx, \xi\in \dR$ for $f\in \dLp{1}$ and can be naturally extended to tempered distributions through duality, and where
\[
\wh{a}(\xi):=\sum_{k\in \dZ} a(k) e^{-ik\cdot \xi},\qquad \xi \in \dR
\]
is the Fourier series of $a\in \dlrs{0}{r}{r}$, which is just an $r\times r$ matrix of $2\pi \dZ$-periodic trigonometric polynomials.
Note that an input sequence $v\in \dlrs{}{1}{r}$ of row vectors is multiplied to the left-hand side of the matrix mask $a$ in the definition of a subdivision operator in \eqref{sd}, while the column vector function $\phi$ is multiplied to the right-hand side of the matrix mask $a$ in the refinement equation \eqref{refeq:phi}.  To avoid potential confusion about notation differences,  we point out here that many references such as \cite{ch19,cmr14,ccs16,cmss19,dm06,dm09,dhmm05,dl95,mer92,ms11,ms17,ms19,md19,mhc20} studying univariate Hermite subdivision schemes adopt the following definition of a vector subdivision operator $\mathbf{S}_{\mathbf{A}}: (l(\dZ))^r\rightarrow (l(\dZ))^r$:
\[
\mathbf{S}_{\mathbf{A}} v:=\sum_{k\in \dZ} \mathbf{A}(j-2k) v(k) =
[\sd_a (v^\tp)]^\tp,\quad v\in (l(\dZ)^{r} \quad \mbox{with}\quad \mathbf{A}:=2^d a^\tp.
\]
Because we study subdivision schemes in this paper through frequently using their underlying refinable vector functions and vector cascade algorithms, to avoid transposes appearing in many places, it is very natural for us to follow the notations used in \cite{han03,hanbook} from the perspective of wavelet analysis.

Define $\dNN:=(\N\cup\{0\})^d$ to be the set of nonnegative multi-integers and $\mathbf{0}:=(0,,\ldots,0)^\tp\in \dNN$.
For $j=1,\ldots,d$, by $\partial_j$ we denote the partial derivative with respect to the $j$th coordinate of $\dR$.
For $\mu=(\mu_1,\ldots,\mu_d)^\tp\in \dNN$ and $x=(x_1,\ldots,x_d)^\tp \in \dR$, we define
\[
|\mu|:=\mu_1+\cdots+\mu_d, \qquad \mu!:=\mu_1!\cdots \mu_d!,\qquad
x^\mu:=x_1^{\mu_1}\cdots x_d^{\mu_d},\qquad  \partial^\mu:=\partial_1^{\mu_1}\cdots \partial_d^{\mu_d},
\]
where $\partial^\mu$ stands for the $\mu$th partial derivative. For a smooth $d$-variate function $f$ on $\dR$, for simplicity we shall use the notation $f^{(\mu)}:=\partial^\mu f$, the $\mu$th partial derivative of $f$.
For $m\in \NN$, we define
\be \label{indm}
\ind_{m}:=\{(\mu_1,\ldots,\mu_d)^\tp \in \dNN \setsp \mu_1+\cdots+
\mu_d \le m\},\qquad m\in \NN.
\ee
To introduce a generalized Hermite subdivision scheme, we shall use ordered multisets $\ind$ whose elements are
from $\ind_m$ in \eqref{indm} (or more generally $\dNN$) but can be repeated with multiplicity. Let $\#\ind$ be the cardinality of an ordered multiset $\ind$ counting multiplicities of elements in $\ind$.
An ordered multiset $\ind$ with cardinality $r:=\#\ind$ considered in this paper can be simply represented by
\be \label{indset}
\ind=\{\nu_1,\ldots,\nu_r\}\quad\mbox{with}\quad
\nu_1,\ldots,\nu_r\in \ind_{m}\quad \mbox{and}\quad \nu_1=\mathbf{0}\in \dNN.
\ee
We now define a generalized Hermite subdivision scheme of type $\ind$.
Let $r:=\#\ind\in \N$ and $a=\{a(k)\}_{k\in \dZ}\in \dlrs{0}{r}{r}$ be a finitely supported matrix mask on $\dZ$.
Let $w_0: \dZ \rightarrow \C^{1\times r}$ be a sequence on $\dZ$ (which is a sequence of  $1\times r$ row vectors) standing for a given input vector sequence/data on $\dZ$.
A sequence of generalized Hermite refinements $w_n: \dZ \rightarrow \C^{1\times r}$ for $n\in \N$ is obtained through recursively applying the vector subdivision operator $\sd_a$ on $w_0$ iteratively as follows:
\be \label{hsd:wn}
w_{n}:=(\sd_a^n w_0) \sD_\ind^{-n},\qquad n\in \N \quad \mbox{with}\quad
\sD_\ind:=\mbox{diag}(2^{-|\mu|})_{\mu\in \ind}=\mbox{diag}(2^{-|\nu_1|},\ldots, 2^{-|\nu_r|}).
\ee
Alternatively, we have the following equivalent form:
\be \label{wn}
w_n=(\sd_a (w_{n-1} \sD_\ind^{n-1}))\sD_\ind^{-n}=
\sd_{\sD^{n-1}_\ind a \sD_\ind^{-n}} w_{n-1},\qquad \mbox{for all } n\in \N.
\ee
The above iterative scheme for computing $\{w_n\}_{n\in \N}$ with the initial input data $w_0$ is called \emph{a generalized Hermite subdivision scheme of type $\ind$}.
For $\ind=\{\mathbf{0}\}$, it is just a scalar subdivision scheme (due to $\#\ind=1$ and a scalar mask $a\in \dlp{0}$) which has been extensively studied (e.g., see \cite{cdm91,dl02,hj98}).
For $m=0$ and $\ind=\{\mathbf{0},\ldots,\mathbf{0}\}$ in \eqref{indset}, it is simply a Lagrange subdivision scheme (also called a vector subdivision scheme in \cite{dm06,hkz09,ms98}).
For $\ind=\ind_m$, it is just a (standard) Hermite subdivision scheme of degree $m$, and a Birkhoff subdivision scheme if $\ind\subsetneq \ind_m$.
If $\ind$ is the $N$ copies of $\ind_m$ with $N\in \N$, as we shall discuss in \cref{sec:hsdint}, then a generalized Hermite subdivision scheme of type $\ind$ is closely related to multivariate interpolating refinable vector functions studied in \cite{hkz09,hz09}.

We now state the convergence of a generalized Hermite subdivision scheme of type $\ind$.

\begin{definition}\label{def:hsd}
{\rm
Let $m\in \NN$, $r\in \N$ and
$a\in \dlrs{0}{r}{r}$ be a finitely supported matrix mask on $\dZ$.
For an ordered multiset $\ind=\{\nu_1,\ldots,\nu_r\}\subseteq \ind_m$ with $\nu_1=\mathbf{0}\in \dNN$ as in \eqref{indset},
we say that \emph{a generalized Hermite subdivision scheme of type $\ind$ with mask $a$ is convergent with limit functions in $\dCH{m}$} if for every input initial vector sequence $w_0=\{w_0(k)\}_{k\in \dZ}: \dZ \rightarrow \C^{1\times r}$, i.e., $w_0\in \dlrs{}{1}{r}$, there exists a $\dCH{m}$ function $\eta: \dR\rightarrow \C$ such that for every constant $K>0$,
\be \label{hsd:converg}
\lim_{n\to \infty} \max_{k\in \dZ\cap [-2^nK, 2^n K]^d} |w_n(k) e_{\ell}- \eta^{(\nu_\ell)}(2^{-n}k)|=0,\qquad \forall\; \ell=1,\ldots,r,
\ee
where $e_\ell$ is the $\ell$th unit coordinate column vector of $\R^r$ and $w_n, n\in \N$ are the generalized Hermite refinement data which are defined in \eqref{hsd:wn} and can be recursively computed via \eqref{wn}.
}
\end{definition}

For convenience, throughout the paper we shall adopt the following big $\bo$ notion: For $m\in \NN$ and smooth functions $f$ and $g$,
\be \label{bo}
f(\xi)=g(\xi)+\bo(\|\xi\|^{m+1}),\quad \xi\to 0 \;\; \mbox{stands for}\;\;
f^{(\mu)}(0)=g^{(\mu)}(0)\qquad \forall\; \mu\in \dNN, |\mu|\le m.
\ee
By $\td$ we denote \emph{the Dirac sequence} such that $\td(0):=1$ and $\td(k):=0$ for all $k\in \dZ\bs\{0\}$.
Using the initial data $w_0=\td I_r\in \dlrs{0}{r}{r}$ in \eqref{hsd:wn} and the convention $\nu_1=\mathbf{0}\in \dNN$ in \eqref{indset} for the ordered multiset $\ind$,
we can define a vector function $\phi \in (\dCH{m})^r$ as the limit vector function for a convergent generalized Hermite subdivision scheme of type $\ind$ in \cref{def:hsd} through
\be \label{def:phi}
\lim_{n\to \infty} \sup_{k\in \dZ}\| (\sd_a^n (\td I_r))(k) e_1-\phi(2^{-n}k) \|=0\qquad \mbox{with}\quad
\phi=[\phi_1,\ldots,\phi_r]^\tp,
\ee
where $e_1$ is the first unit coordinate column vector in $\R^r$ and we used the fact $\sD_\ind^{-n} e_1=e_1$ due to $\nu_1=\mathbf{0}$.
The vector function $\phi$ in \eqref{def:phi} is called \emph{the basis vector function} of a generalized Hermite subdivision scheme in \cref{def:hsd}, because
the limit function $\eta$ in \cref{def:hsd} is given by $\eta=w_0*\phi:=\sum_{k\in \dZ} w_0(k) \phi(\cdot-k)$ for any initial sequence $w_0\in \lrs{}{1}{r}$.
To study the properties of a basis vector function $\phi$, we recall the definition of sum rules of a matrix mask $a\in \dlrs{0}{r}{r}$.
For a nonnegative integer $m\in \NN$, we say that a matrix mask \emph{$a\in \dlrs{0}{r}{r}$ has order $m+1$ sum rules with a matching filter $\vgu_a\in \dlrs{0}{1}{r}$} (e.g., see \cite[(1.6)]{han01} or \cite[(2.8) and (2.9)]{han03}) if $\wh{\vgu_a}(0)\ne 0$ and
\be \label{sr}
\wh{\vgu_a}(2\xi)\wh{a}(\xi)=\wh{\vgu_a}(\xi)+\bo(\|\xi\|^{m+1})
\; \mbox{ and }\;
\wh{\vgu_a}(2\xi)\wh{a}(\xi+\pi \omega)=\bo(\|\xi\|^{m+1}),\; \xi\to 0,\; \forall\;
\omega\in \Gamma\bs \{0\},
\ee
where $\Gamma:=[0,1]^d\cap \dZ$.
A filter $\vgu_a\in \dlrs{0}{1}{r}$ satisfying the first identity of \eqref{sr} and $\wh{\vgu_a}(0)\ne 0$ is often called \emph{an order $m+1$ matching filter} of the matrix mask $a$. By $\PL_m$ we shall denote the linear space of all $d$-variate polynomials of (total) degree no more than $m$.
The notion of sum rules plays a key role in studying subdivision schemes and wavelets  by linking it to polynomial reproduction and smoothness of its refinable vector function (e.g., \cite{han01,han03,hanbook} and references therein).

For a convergent generalized Hermite subdivision scheme of type $\ind$ with limit functions in $\dCH{m}$, we now show that its basis vector function $\phi$ must be a refinable vector function satisfying \eqref{refeq:phi}, which helps us to characterize its matrix mask through sum rules and the underlying matching filter. For improved readability, the proof of the following result is given in \cref{sec:hsdmask}.

\begin{theorem}\label{thm:hsd:mask}
Let $m\in \NN, r\in \N$ and an ordered multiset $\ind=\{\nu_1,\ldots,\nu_r\} \subseteq \ind_m$ as in \eqref{indset}.
%Let $a\in \dlrs{0}{r}{r}$ be a finitely supported matrix mask on $\dZ$.
Suppose that the generalized Hermite subdivision scheme of type $\ind$ with a mask $a\in \dlrs{0}{r}{r}$ is convergent with limit functions in $\dCH{m}$, as in the sense of \cref{def:hsd}. Then its basis vector function $\phi$ defined in \eqref{def:phi} is a compactly supported refinable vector function in $\dCH{m}$ satisfying
\be \label{refeq}
\wh{\phi}(\xi)=\lim_{n\to \infty} \Big(\Big[\prod_{j=1}^n \wh{a}(2^{-j}\xi)\Big] e_1\Big)
\qquad \mbox{and}\qquad \wh{\phi}(2\xi)=\wh{a}(\xi)\wh{\phi}(\xi)\quad \forall\, \xi\in \dR.
\ee
If in addition
\be \label{phi:cond0}
\mspan\{\wh{\phi}(\pi \omega+2\pi k) \setsp k\in \dZ\}=\C^r,\qquad \forall\; \omega\in [0,1]^d\cap \dZ,
\ee
then $\wh{\phi}(0)\ne 0$ and the following statements hold:
\begin{enumerate}
\item[(1)] The matrix matrix $\wh{a}(0):=\sum_{k\in \dZ} a(k)$ must satisfy that
\be \label{eig:one}
\mbox{$1$ is a simple eigenvalue of $\wh{a}(0)$ and all its other eigenvalues
are less than $2^{-m}$ in modulus.}
\ee
Consequently, there always exists a finitely supported sequence $\vgu_a\in \dlrs{0}{1}{r}$ whose values $\wh{\vgu_a}^{(\mu)}(0)$ for $\mu\in \ind_{m}$ are uniquely determined by
\be \label{vgua:0}
\wh{\vgu_a}(0)\wh{a}(0)=\wh{\vgu_a}(0) \quad \mbox{with}\quad \wh{\vgu_a}(0)\wh{\phi}(0)=1
\ee
and the following recursive formula:
\be \label{vgua:val}
\wh{\vgu_a}^{(\mu)}(0)=\sum_{0\le \beta \le \mu, \beta\ne \mu} \frac{2^{|\beta|} \mu!}{\beta! (\mu-\beta)!}
\wh{\vgu_a}^{(\beta)}(0)\wh{a}^{(\mu-\beta)}(0)[I_r-2^{|\mu|} \wh{a}(0)]^{-1},\qquad \mu\in \ind_{m}\bs\{\mathbf{0}\}
\ee
from $|\mu|=1$ to $|\mu|=m$, where $\beta\le \mu$ means that all the entries in $\mu-\beta$ are nonnegative.

\item[(2)] For any $\vgu_a\in \dlrs{0}{1}{r}$ satisfying \eqref{vgua:0} and \eqref{vgua:val} in item (1), the sequence $\vgu_a$ must satisfy
%
%\be \label{phi:poly}
\[
\wh{\vgu_a}(\xi)\wh{\phi}(\xi+2\pi k)=\td(k) +\bo(\|\xi\|^{m+1}),\quad \xi \to 0,\; \mbox{for all }\, k\in \dZ.
\]
Consequently, $\pp=\sum_{k\in \dZ} (\pp*\vgu_a)(k) \phi(\cdot-k)$ for all $\pp\in \PL_m$, where $\pp*\vgu_a:=\sum_{k\in \dZ} \pp(\cdot-k)\vgu_a(k)$.

\item[(3)]
The mask $a$ has order $m+1$ sum rules with a matching filter $\vgu_a\in \dlrs{0}{1}{r}$ satisfying
\be \label{vgua:hsd}
\wh{\vgu_a}(\xi)=
\left[(i\xi)^{\nu_1}+\bo(\|\xi\|^{|\nu_1|+1}), (i\xi)^{\nu_2}+\bo(\|\xi\|^{|\nu_2|+1}),\ldots,
(i\xi)^{\nu_r}+\bo(\|\xi\|^{|\nu_r|+1})\right], \quad \xi\to 0.
\ee

\item[(4)] For all $\mu\in \ind_m$, the vector polynomial $\vec{\pp}_\mu:=\frac{(\cdot)^\mu}{\mu!}*\vgu_a$ satisfies
$\sd_a \vec{\pp}_\mu=2^{-|\mu|}\vec{\pp}_\mu$ and
%
%\be \label{sd:pmu}
\[
\vec{\pp}_\mu(x) e_\ell=
\begin{cases}
\frac{x^{\mu-\nu_\ell}}{(\mu-\nu_\ell)!}+\sum_{0\le \beta\le \mu, |\beta|>|\nu_\ell|}
\frac{(-i)^{|\beta|}}{\beta!} \frac{x^{\mu-\beta}}{(\mu-\beta)!}
\wh{\vgu_a}^{(\beta)}(0)e_\ell, &\text{if $\nu_\ell\le \mu$},\\
\sum_{0\le \beta\le \mu, |\beta|>|\nu_\ell|}
\frac{(-i)^{|\beta|}}{\beta!} \frac{x^{\mu-\beta}}{(\mu-\beta)!}
\wh{\vgu_a}^{(\beta)}(0)e_\ell,
&\text{otherwise},
\end{cases}
\quad \forall\; \ell=1,\ldots,r,
\]
and $\sum_{k\in \dZ} \vec{\pp}_\mu(k)\phi(\cdot-k)=\frac{(\cdot)^\mu}{\mu!}$, where $\vgu_a\in \dlrs{0}{1}{r}$ satisfies \eqref{vgua:hsd}.
\end{enumerate}
\end{theorem}

A matrix mask $a\in \dlrs{0}{r}{r}$ is called \emph{a generalized Hermite mask of type $\ind$ having order $m+1$ sum rules} if  the mask $a$ has order $m+1$ sum rules with a matching filter $\vgu_a\in \dlrs{0}{1}{r}$ satisfying \eqref{vgua:hsd}, i.e., item (3) of \cref{thm:hsd:mask} holds.
By linking a basis vector function $\phi$ to a refinable vector function in \eqref{refeq}, item (3) of \cref{thm:hsd:mask} characterizes the matrix mask of a convergent generalized Hermite subdivision scheme of type $\ind$ with limit functions in $\dCH{m}$.
Moreover, its basis vector function $\phi$ and its vector subdivision operator $\sd_a$ must reproduce polynomials in $\PL_m$ in the special way as described in
items (2) and (4) of \cref{thm:hsd:mask}.

\subsection{Convergence and smoothness of generalized Hermite subdivision schemes}

To characterize convergence and smoothness of generalized Hermite subdivision schemes, built on \cref{thm:hsd:mask}
we now connect them to refinable vector functions and vector cascade algorithms.
For $1\le p\le \infty$, the \emph{vector refinement operator}
$\cd_a: (\dLp{p})^r \rightarrow (\dLp{p})^r$
is defined to be
\be \label{cd}
\cd_a f:=2^d \sum_{k\in \dZ} a(k) f(2\cdot-k),\qquad f\in (\dLp{p})^r.
\ee
Note that the refinement equation \eqref{refeq:phi} can be restated as $\cd_a \phi=\phi$, i.e., $\phi$ is a fixed point of the vector refinement operator $\cd_a$.
Similar to a subdivision scheme for discrete data, a cascade algorithm iteratively computes a sequence $\{\cd_a^n f\}_{n=1}^\infty$ of vector functions on $\R^d$ such that this sequence may converge to
the refinable vector function
$\phi$ in certain function spaces.
%such as the Sobolev space $W_p^m(\dR)$. \cite[Theorem~4.3]{han03} shows that a vector cascade algorithm with a mask $a\in \dlrs{0}{r}{r}$ converges in $W_p^m(\dR)$ if and only if $\sm_p(a)>m$, where the quantity $\sm_p(a)$ is defined in \cite[(4.3)]{han03} and is given in Subsection~\ref{subsec:sma}.
For $0<\tau\le 1$ and a function $f\in \dLp{p}$, we say that $f$ belongs to the Lipschitz space $\mbox{Lip}(\tau, \dLp{p})$ if there exists a positive constant $C$ such that $\|f-f(\cdot-t)\|_{\dLp{p}}\le C\|t\|^\tau$ for all $t\in \dR$. For convenience, we define $\mbox{Lip}(0, \dLp{p}):=\dLp{p}$.
The $L_p$ smoothness of a function $f\in \dLp{p}$ is measured by its \emph{$L_p$ critical exponent} $\sm_p(f)$ defined by
\[
\sm_p(f):=\sup\{ m+\tau \setsp  0\le \tau<1\; \mbox{ and }\; \partial^\mu f\in \mbox{Lip}(\tau, \dLp{p}) \; \mbox{ for all } \mu\in \ind_m, m\in \NN\}.
\]
If $\phi=[\phi_1,\ldots,\phi_r]^\tp$ is a vector function, then we define $\sm_p(\phi):=\min(\sm_p(\phi_1),\ldots,\sm_p(\phi_r))$.

For a column vector $\phi$ of compactly supported distributions on $\R^d$, we say that the integer shifts of $\phi$ are \emph{stable} if
\be \label{stable}
\mspan\{\wh{\phi}(\xi+2\pi k) \setsp k\in \dZ\}=\C^r, \qquad \forall\; \xi\in \dR.
\ee

We now characterize the convergence and smoothness of a generalized Hermite subdivision scheme in the following result, whose technical proof is deferred to \cref{sec:hsdconverg}.

\begin{theorem}\label{thm:converg}
Let $m\in \NN$, $r\in \N$, and $\ind=\{\nu_1,\ldots,\nu_r\}\subseteq \ind_m$ be an ordered multiset given in \eqref{indset}. Let $a\in \dlrs{0}{r}{r}$ be a generalized Hermite mask of type $\ind$ having order $m+1$ sum rules, i.e., the mask $a$ has order $m+1$ sum rules in \eqref{sr} with a matching filter $\vgu_a\in \dlrs{0}{1}{r}$ satisfying \eqref{vgua:hsd}.
Let $\phi=[\phi_1,\ldots,\phi_r]^\tp$ be a vector of compactly supported
distributions satisfying the refinement equation $\wh{\phi}(2\xi)=\wh{a}(\xi)\wh{\phi}(\xi)$ and the normalization condition $\wh{\vgu_a}(0)\wh{\phi}(0)=1$.
If
\begin{enumerate}
\item[(i)] $\sm_\infty(a)>m$ (see Subsection~\ref{subsec:sma} for definition), or equivalently, the vector cascade algorithm associated with the mask $a$ is convergent in $\dCH{m}$: For every compactly supported admissible initial vector function $f\in (\dCH{m})^r$ satisfying
\be \label{initialf}
\wh{\vgu_a}(0)\wh{f}(0)=1 \quad \mbox{and}\quad
\wh{\vgu_a}(\xi)\wh{f}(\xi+2\pi k)=\bo(\|\xi\|^{m+1}),\quad \xi \to 0, \forall\; k\in \dZ\bs\{0\},
\ee
the cascade sequence $\{\cd_a^n f\}_{n=1}^\infty$ converges to the refinable vector function $\phi$ in $(\dCH{m})^r$, i.e.,
\be \label{cd:converg}
\phi \in (\dCH{m})^r \quad \mbox{and}\quad \lim_{n\to \infty} \|\cd_a^n f-\phi\|_{(\dCH{m})^r}=0,
\ee
\end{enumerate}
then $\sm_\infty(\phi)\ge \sm_\infty(a)$ and
\begin{enumerate}
\item[(ii)] the generalized Hermite subdivision scheme of type $\ind$ with mask $a\in \dlrs{0}{r}{r}$ is convergent with limit functions in $\dCH{m}$, as in the sense of \cref{def:hsd}.
\end{enumerate}
Conversely, if item (ii) holds and
if the integer shifts of $\phi$ are stable,
then item (i) must hold and $\sm_\infty(\phi)=\sm_\infty(a)>m$.
\end{theorem}

It is known (e.g., \cite{jm91} and \cite[Theorem~5.3.4]{hanbook}) that
a compactly supported continuous vector function $\phi$ has stability as in \eqref{stable} if and only if there exist positive constants $C_1$ and $C_2$ such that
\be\label{stable:2}
C_1 \|w_0\|_{\infty}\le \Big \| \sum_{k\in \dZ} w_0(k) \phi(\cdot -k)\Big \|_{\mathscr{C}(\dR)} \le C_2 \|w_0\|_{\infty},\qquad \forall\; w_0\in \dlrs{\infty}{1}{r}.
\ee
Because $\phi$ has compact support, the upper bound in \eqref{stable:2} holds with $C_2=\sup_{x\in \dR} \sum_{k\in \dZ} \|\phi(x-k)\|_{l_1}<\infty$.
For a convergent generalized Hermite subdivision scheme of type $\ind$, the limit function $\eta$ in \cref{def:hsd} is given by $\eta=w_0*\phi=\sum_{k\in \dZ} w_0(k) \phi(\cdot-k)$. Therefore, the stability in \eqref{stable:2} guarantees that $C_1 \|w_0\|_{\infty}\le \|\eta\|_{\mathscr{C}(\dR)}\le C_2\|w_0\|_{\infty}$. The stability condition \eqref{stable} in \cref{thm:converg} cannot be removed even for univariate scalar subdivision schemes.
Suppose now that \eqref{stable} fails. One one hand, for any arbitrarily small $\varepsilon>0$, there exists $w_0\in \dlrs{\infty}{1}{r}$ with $\|w_0\|_{\infty}=1$ but the limit function $\eta:=w_0*\phi$ satisfies $\|\eta\|_{\mathscr{C}(\dR)}<\varepsilon$.
On the other hand, for any $m\in \NN$ and $\ind=\{0\}$ (i.e., the scalar case),
examples of univariate scalar masks $a\in l_0(\Z)$ are given at the end of \cite[Section~5]{han20}
such that $\sm_\infty(a)=0$
%(even $\sm_\infty(a)<0$ by modifying these examples)
but its scalar refinable function $\phi\in \mathscr{C}^m(\R)$.
For more details, see \cite{han06} and references therein for studying refinable vector functions and vector cascade algorithms without the stability condition in \eqref{stable}. \cref{thm:converg} generalizes \cite[Theorem~1.3]{han20} even in dimension one.

Let $m\in \NN$ and $\ind$ be an arbitrarily given ordered multiset in \eqref{indset}.
To construct a convergent generalized Hermite subdivision scheme of type $\ind$ with limit functions in $\dCH{m}$, by \cref{thm:hsd:mask,thm:converg},
the task is to construct a matrix mask $a\in \dlrs{0}{r}{r}$ such that $\sm_\infty(a)>m$ and $a$ is a generalized Hermite mask of type $\ind$ having order $m+1$ sum rules, i.e., item (3) of \cref{thm:hsd:mask} holds.
Moreover, using \cref{thm:hsd:mask,thm:converg}, we shall prove in
\cref{thm:existence} that
there always exists a generalized Hermite mask $a\in \dlrs{0}{r}{r}$ such that $\sm_\infty(a)>m$ and its basis vector function $\phi$ is a spline refinable vector function in $\dCH{m}$. Hence, by \cref{thm:converg}, the generalized Hermite subdivision scheme of type $\ind$ with such a mask $a$ is convergent with limit functions in $\dCH{m}$.

\subsection{Characterize refinable generalized Hermite interpolants}\label{subsec:int}

We now study a special family of generalized Hermite subdivision schemes whose basis vector functions have the interpolation property.
A vector function $\phi=[\phi_1,\ldots,\phi_r]^\tp$ in $\dCH{m}$ is \emph{a generalized Hermite interpolant of type $\ind$ and translation $T$} if it satisfies the following generalized Hermite interpolation property:
\be \label{int:phi}
\phi^{(\nu_\ell)}(k+\tau_\ell)=\td(k) e_{\ell},
\quad \mbox{i.e.},\quad
\phi^{(\nu_\ell)}_j(k+\tau_\ell)=\td(k) \td(\ell-j), \quad \forall\; k\in \dZ \mbox{ and } j, \ell=1,\ldots,r,
\ee
for some ordered multiset $T:=\{\tau_1,\ldots,\tau_r\}\subseteq \dR$ for translation, which is also related to dual Hermite subdivision schemes (see \cref{sec:hsdint} for details). For $f\in \dCH{m}$, we deduce from \eqref{int:phi} that the approximating function $g(x):=\sum_{\ell=1}^r \sum_{k\in \dZ} f^{(\nu_\ell)}(k+\tau_\ell)\phi_\ell(x-k)$ interpolates $f$ such that
$g^{(\nu_\ell)}(k+\tau_\ell)=f^{(\nu_\ell)}(k+\tau_\ell)$ for all $k\in \dZ$ and $\ell=1,\ldots,r$.
Consequently, a generalized Hermite interpolant $\phi$ of type $\ind$ and translation $T$ satisfying \eqref{int:phi} obviously has linearly independent integer shifts, that is, $\sum_{k\in \dZ} v(k) \phi(\cdot-k)=0$ for a sequence $v\in \dlrs{}{1}{r}$ always implies $v(k)=0$ for all $k\in \dZ$.
The generalized Hermite interpolation property in \eqref{int:phi} is just the Lagrange interpolation property for $\ind=\{\mathbf{0},\ldots,\mathbf{0}\}$, the (standard) Hermite interpolation property of degree $m$ for $\ind=\ind_m$ and $T=\{0,\ldots,0\}$, and the Birkhoff interpolation property for $\ind\subsetneq \ind_m$.
%If $\ind$ is the $N$ copies of $\ind_m$ with $N\in \N$, then the generalized Hermite interpolation property in \eqref{int:phi} has been extensively studied in \cite{hz09} for multivariate interpolating refinable vector functions.
For a generalized Hermite subdivision scheme, here we are particularly interested in establishing criteria such that its basis vector function $\phi$
has the generalized Hermite interpolation property in \eqref{int:phi}.

\begin{definition}\label{def:int}
\normalfont
Let $\ind=\{\nu_1,\ldots,\nu_r\}$ be an ordered multiset as in \eqref{indset} for derivatives and $T=\{\tau_1,\ldots, \tau_r\}\subseteq \R^d$ be an ordered multiset for translation. Let $\theta$ be a mapping on $\{1,\ldots,r\}$ satisfying
\be \label{theta}
\theta: \{1,\ldots,r\}\rightarrow \{1,\ldots,r\},
\ell \mapsto \theta(\ell)
\quad\mbox{such that}\quad
\nu_{\theta(\ell)}=\nu_\ell \quad \mbox{and}\quad
2\tau_\ell-\tau_{\theta(\ell)}\in \dZ.
\ee
We say that a generalized Hermite subdivision scheme of type $\ind$ with mask $a\in \dlrs{0}{r}{r}$ is \emph{an interpolatory generalized Hermite subdivision scheme of type $\ind$ and translation $T$} if for any initial data $w_0: \dZ\rightarrow \C^{1\times r}$, the refinement data $\{w_n\}_{n\in \N}$ in \eqref{hsd:wn} satisfy the following interpolation condition:
\be \label{interpolation}
w_n(2k+2\tau_\ell-\tau_{\theta(\ell)})e_{\theta(\ell)}=w_{n-1}(k) e_{\ell},\qquad \forall\; k\in \dZ, n\in \N, \ell=1,\ldots,r.
\ee
\end{definition}

We shall explain the role of the mapping $\theta$ in \cref{sec:hsdint}. The condition \eqref{theta} on the mapping $\theta$ is mainly due to the refinability of generalized Hermite interpolants. For interpolatory generalized Hermite subdivision schemes, we have the following result, whose proof is given in \cref{sec:hsdint}.

\begin{theorem}\label{thm:int}
Let $m\in \NN$ and $r\in \N$.
Take an ordered multiset $\ind=\{\nu_1,\ldots,\nu_r\}\subseteq \ind_m$ with $\nu_1=\mathbf{0}\in \dNN$ as in \eqref{indset} for derivatives,
an ordered multiset $T=\{\tau_1,\ldots, \tau_r\}\subseteq \R^d$ for translation, and  a mapping $\theta: \{1,\ldots,r\}\rightarrow \{1,\ldots,r\}$ satisfying \eqref{theta}. Let $a\in \dlrs{0}{r}{r}$ be a finitely supported matrix mask.
If the following conditions are satisfied:
\begin{enumerate}
\item[(i)] the generalized Hermite subdivision scheme of type $\ind$ with mask $a$
is convergent with limit functions in $\dCH{m}$ as in the sense of \cref{def:hsd};
\item[(ii)] the generalized Hermite subdivision scheme of type $\ind$ with mask $a$ is interpolatory of type $\ind$ and translation $T$ as in the sense of \cref{def:int};

\item[(iii)] the basis vector function $\phi$ in \eqref{def:phi} associated with mask $a$ satisfies the condition in \eqref{phi:cond0},
\end{enumerate}
then all the following statements hold:
\begin{enumerate}
\item[(1)] The matrix mask $a$ satisfies the following interpolation condition:
\be \label{int:mask}
a(2k+2\tau_\ell-\tau_{\theta(\ell)})e_{\theta(\ell)}=
2^{-d-|\nu_\ell|} \td(k) e_{\ell},\quad \forall\; k\in \dZ, \ell=1,\ldots,r.
\ee

\item[(2)] The vector function $\phi$ is a refinable generalized Hermite interpolant of type $\ind$ and translation $T$ satisfying the refinement equation \eqref{refeq:phi} and the interpolation condition \eqref{int:phi}.

\item[(3)]
Let $M\in \NN$ be the largest possible integer such that the mask $a$ has order $M+1$ sum rules with a matching filter $\vgu_a\in \dlrs{0}{1}{r}$ satisfying $\wh{\vgu_a}(0)\wh{\phi}(0)=1$. Then $M\ge m$ and
\be \label{int:vgua}
\wh{\vgu_a}(\xi)=
\left[(i\xi)^{\nu_1}e^{i\tau_1\cdot\xi},\ldots,
(i\xi)^{\nu_r}e^{i\tau_r\cdot \xi}\right]+\bo(\|\xi\|^{M+1}),
\quad \xi \to 0.
\ee

\item[(4)] For all $\mu\in \ind_M$, the vector polynomial $\vec{\pp}_\mu:=\frac{(\cdot)^\mu}{\mu!}*\vgu_a$ satisfies $\sd_a \vec{\pp}_\mu=2^{-|\mu|} \vec{\pp}_\mu$, $\sum_{k\in \dZ} \vec{\pp}_\mu (k)\phi(\cdot-k)=\frac{(\cdot)^\mu}{\mu!}$ and
\be\label{int:pmu}
\vec{\pp}_\mu(x) e_\ell=
\partial^{\nu_\ell} \left(\frac{(x+\tau_\ell)^\mu}{\mu!}\right)=
\begin{cases}
\frac{(x+\tau_\ell)^{\mu-\nu_\ell}}{(\mu-\nu_\ell)!},
&\text{if $\nu_\ell\le \mu$},\\
0, &\text{otherwise},
\end{cases}
\qquad \forall\; \ell=1,\ldots,r.
\ee
\end{enumerate}
Conversely, if a refinable vector function $\phi$ associated with mask $a\in \dlrs{0}{r}{r}$ is a generalized Hermite interpolant of type $\ind$ and translation $T$ (i.e., $\phi$ satisfies the generalized Hermite interpolation condition in \eqref{int:phi}, or equivalently, item (2) holds), then all the items (i)--(iii) and (1)--(4) must hold.
\end{theorem}

A matrix mask $a\in \dlrs{0}{r}{r}$ is \emph{interpolatory of type $\ind$ and translation $T$} if \eqref{int:mask} is satisfied. As we shall see in the proof of \cref{thm:int} in Section~\ref{sec:hsdint}, if a mask $a\in \dlrs{0}{r}{r}$ is interpolatory of type $\ind$ and translation $T$ and if the mask $a$ has order $M+1$ sum rules with a matching filter $\vgu_a\in \dlrs{0}{1}{r}$, then the matching filter $\vgu_a$ must satisfy \eqref{int:vgua} in item (3) of \cref{thm:int}. To construct a refinable generalized Hermite interpolant in $\dCH{m}$ and a convergent interpolatory generalized Hermite subdivision scheme of type $\ind$ and translation $T$,
by \cref{thm:converg,thm:int}, the task is to construct a matrix mask $a\in \dlrs{0}{r}{r}$ such that $\sm_\infty(a)>m$ and $a$ is an interpolatory generalized Hermite mask of type $\ind$ and translation $T$ having order $M+1$ sum rules with a matching filter $\vgu_a\in \dlrs{0}{1}{r}$ satisfying \eqref{int:vgua} with some $M\ge m$.

\subsection{Literature review, motivations and contributions}
\label{subsec:review}

We now review related works on Hermite subdivision schemes.
\emph{A (standard) Hermite subdivision scheme of degree $m$} is just a generalized Hermite subdivision scheme of type $\ind_m$. \cref{def:hsd} with $d=1$ for univariate Hermite subdivision schemes is given in \cite{hyx05,dm09,ms19} and references therein.
Moreover, \cref{def:int} with $\ind=\ind_m$ and $T=\{0,\ldots,0\}$
becomes the standard interpolatory Hermite subdivision schemes.
For $\ind=\{0,1\}$, the univariate interpolatory Hermite subdivision schemes of degree $1$ were studied in Merrien
\cite{mer92} and in
Dyn and Levin \cite{dl95}. \cite{zhou00} studied univariate Hermite interpolants. A family of univariate interpolatory Hermite masks is presented in \cite[Theorems~4.2 and 4.3]{han01}.
Convergence of univariate Hermite subdivision schemes have been extensively studied by many researchers through factorization of matrix masks and their derived subdivision schemes in \cite{cmr14,ccs16,cmss19,dm06,dm09,ju02,ms11,ms17,ms19,mhc20} and references therein. Recently, \cite{han20} characterized univariate Hermite masks and the convergence of univariate Hermite subdivision schemes without factorizing univariate Hermite masks. Further related topics have been addressed in \cite{ch19,jy19,md19} and references therein.
For the ordered multiset $\ind=\{\mathbf{0},\ldots,\mathbf{0}\}$ with $\mathbf{0}\in \dNN$ being repeated $r$ times,
due to $\sD_\ind=I_{r}$ in \eqref{hsd:wn}, \cref{def:hsd} agrees with the definition of convergence of Lagrange subdivision schemes in \cite{hyx05} and is just the convergence of vector subdivision schemes in \cite[Definition~2.1]{ms98} and references therein. Lagrange subdivision schemes were studied in \cite{dm06,ms98} and univariate refinable Lagrange interpolants were studied in \cite{hkz09} and references therein.
If $\ind$ is a proper subset of $\{0,1,\ldots,r\}$, e.g., $\ind=\{0,2\}$, \cref{def:hsd} is related to Birkhoff interpolation and even in dimension one we are not aware of any work on univariate Birkhoff subdivision schemes.

For multiple dimensions,
multivariate refinable Hermite interpolants have been characterized in \cite[Corollary~5.3]{han03} and constructed in \cite{dhmm05,hm07}.
The interpolating refinable vector functions in \cite{hz09}
are just generalized Hermite interpolants of type $\ind$ with translation $T$ such that $\ind$ contains $N$ copies of $\ind_m$. Such particular generalized Hermite interpolants have been analyzed and constructed in \cite{hz09}. However, in these papers, their associated multivariate interpolatory generalized Hermite subdivision schemes have not been addressed there, because these papers focus on the interpolation property of the underlying refinable vector functions through studying the convergence of vector cascade algorithms.
Built on the work in \cite{han03},
multivariate interpolatory Hermite subdivision schemes have been studied in \cite{hyp04}. Convergence of multivariate vector subdivision schemes was studied in \cite{ccs05} through the factorization of multivariate masks.
To our best knowledge, \cite{hyx05} is the only known paper discussing multivariate Hermite subdivision schemes by providing a restrictive sufficient condition in \cite[Theorem~2.2]{hyx05}.
As a byproduct of \cite{hyx05}, face-based (i.e., dual) Hermite subdivision schemes were addressed in \cite{hy06}.
However, we are not aware of any other work addressing multivariate Hermite subdivision schemes and multivariate generalized Hermite subdivision schemes using an ordered multiset $\ind$.
%Due to the importance of generalized Hermite subdivision schemes and generalized Hermite interpolants in mathematics,

This paper is also motivated by an example of refinable spline spaces known in approximation theory  and numerical PDEs.
%\begin{example}\label{expl:spline} \normalfont
For $m\in \NN$ and $N\in \N$,
the $\CH{m}$ piecewise polynomial space of degree at most $(m+1)N+m$ is given by
\[
S_{m,N}:=\{ f\in \CH{m} \setsp
f|_{(k, k+1]} \mbox{ is a polynomial of degree } \le (m+1)N+m \mbox{ for all } k\in \Z\}.
\]
Let $\ind$ and $T$ be the $N$ copies of $\{0,1,\ldots,m\}$ and $\{\frac{0}{N},\frac{1}{N},\ldots,\frac{N-1}{N}\}$, respectively. Explicitly, $\ind=\{\nu_1,\ldots,\nu_r\}$ and $T=\{\tau_1,\ldots,\tau_r\}$ with $r:=(m+1)N$, where $\nu_\ell:=(\ell-1) \mod (m+1)$ and $\tau_\ell:={\lfloor \frac{\ell-1}{m+1}\rfloor}/{N}$ for $\ell=1,\ldots,r$.
Define
$\pq_\ell$ to be the unique polynomial of degree $m$ satisfying
\[
\pq_\ell (x+\tau_\ell)
=\frac{x^{\nu_\ell}}
{\nu_\ell !\pp_\ell(x+\tau_\ell)}+\bo(x^{m+1}), \quad x\to 0
\quad\mbox{with}\quad
\pp_\ell(x):=\prod_{k\in \{0,\ldots,N\}\bs \{N\tau_\ell\}} \left(x-\frac{k}{N}\right)^{m+1}.
\]
Now we can construct a basis vector function $\phi=[\phi_1,\ldots,\phi_r]^\tp$ for the above spline space by defining
\be \label{phi:spline}
\phi_\ell:=
\begin{cases}
\pp_\ell\pq_\ell\chi_{[0,1]}+
(-1)^{\nu_\ell} \pp_\ell(-\cdot)\pq_\ell(-\cdot)\chi_{[-1,0)},
&\ell=1,\ldots,m+1,\\
\pp_\ell \pq_\ell\chi_{[0,1]}, &\ell=m+2,\ldots, (m+1)N.
\end{cases}
\ee
Note that $\deg(\pp_\ell \pq_\ell)\le (m+1)N+m$ and
the Hermite cubic splines correspond to $\phi$ in \eqref{phi:spline} with $m=1$ and $N=1$.
For $\ell=m+2,\ldots, (m+1)N$, we observe $0<\tau_\ell<1$ and hence, $x^{m+1}(x-1)^{m+1}\mid \pp_\ell(x)$, from which we conclude $\phi_\ell\in S_{m,N}$. For $\ell=1,\ldots,m+1$, we have $\tau_\ell=0$ and hence $(x-1)^{m+1}\mid \pp_\ell(x)$ and $\phi_\ell(x)=\frac{x^{\nu_\ell}}{\nu_\ell!}+\bo(x^{m+1})$ as $\xi\to 0$, from which we have $\phi_\ell\in S_{m,N}$. By calculation, we can directly check that $\phi$ in \eqref{phi:spline} is indeed a generalized Hermite interpolant of type $\ind$ and translation $T$ satisfying \eqref{int:phi}.
Now we can conclude that the integer shifts of $\phi$ in \eqref{phi:spline} generate the spline space $S_{m,N}$.
Moreover, $\phi$ is a refinable vector function satisfying the refinement equation in \eqref{refeq:phi} with a mask $a\in \lrs{0}{r}{r}$ supported on $\{-1,0,1\}$ and given by
\be \label{spline:a}
a(k)=\left[2^{-1-\nu_1} \phi^{(\nu_1)}(2^{-1}(k+\tau_1)),\ldots,
2^{-1-\nu_r}\phi^{(\nu_r)}(2^{-1}(k+\tau_r))
\right],\qquad k\in \Z
\ee
which can be obtained using the interpolation property in \eqref{int:phi}. Define a mapping on $\{1,\ldots,r\}$ by
\[
\theta(\ell):=\ell+(m+1)N_\ell,\quad \ell=1,\ldots,r
\quad \mbox{with}\quad
N_\ell:=\begin{cases}
N\tau_\ell &\text{if $\ell\le (m+1)N (1-\tau_\ell)$},\\
N\tau_\ell-N, &\text{if $\ell>(m+1)N (1-\tau_\ell)$}.
\end{cases}
\]
%Note that the above choice of $N_\ell\in \Z$ guarantees $\theta(\ell)\in \{1,\ldots,r\}$.
Because $\nu_\ell:=(\ell-1) \mod (m+1)$ and $N_\ell\in \Z$, we have $\nu_{\theta(\ell)}=\nu_\ell$ and $2\tau_\ell-\tau_{\theta(\ell)}=\tau_\ell-\frac{N_\ell}{N}\in \{0,1\}$ by the choice of $N_\ell$ for all $\ell=1,\ldots,r$.
Hence, the above mapping $\theta$ indeed satisfies \eqref{theta}.
By \cref{thm:converg,thm:int}, the interpolatory generalized Hermite subdivision scheme of type $\ind$ and translation $T$ with such a matrix mask $a\in \lrs{0}{r}{r}$ is convergent with limit functions in $\CH{m}$.
%Many other known multivariate spline spaces in mathematics are generated by generalized Hermite interpolants, see \cite{cdm91,cj03,hm07,hyp04,hy06,hz09} and references therein.
%In dimension $d$, the tensor product spline space $\otimes^d S_{m,N}$ has a generalized Hermite interpolant generator $\otimes^d \phi$ with ordered multisets $\otimes^d \ind$ and $\otimes^d T$, where $\otimes^d\phi$ stands for the tensor product of $d$ copies of $\phi$.
%Note that the set $\ind_m$ does not have the tensor product structure while $\ind$ in this paper covers the tensor product of $\ind_m$.

To apply wavelets to the boundary value problems in numerical PDEs, one key task is to adapt wavelets from the Euclidean spaces to bounded domains (\cite{dah97}). Recently, \cite{hm21} developed a systematic direct approach which can adapt any univariate compactly supported biorthogonal multiwavelet from the real line $\R$ to the interval $[0,1]$ with or without any homogeneous boundary conditions. The special spline generalized Hermite interpolants in \eqref{phi:spline} motivate us to search for new refinable vector functions, which are associated with generalized Hermite subdivision schemes and have some desired properties for numerical PDEs.
Motivated by our study of univariate Hermite subdivision schemes in \cite{han20} and references therein, in this paper we shall comprehensively study various properties of multivariate generalized Hermite subdivision schemes in \cref{def:hsd}  with or without the interpolation property and/or the spline connections. Further motivated by \cref{thm:int} on refinable generalized Hermite interpolants, we shall
introduce the notion of linear-phase moments for the polynomial-interpolation property.
%\cref{def:hsd} can be extended to a general dilation matrix but the presentation is much more involved. Consequently we shall restrict ourselves for the dyadic case only in this paper.
Using ordered multisets $\ind$, our results in this paper on generalized Hermite subdivision schemes significantly advance the theory on standard multivariate Hermite subdivision schemes in several aspects as byproducts.
The contributions of this paper are as follows:
\begin{enumerate}
\item Introduce the unifying notion of generalized Hermite subdivision schemes, which include all known variants of Hermite-like subdivision schemes as special cases and introduce new types of subdivision schemes such as general Lagrange or Birkhoff subdivision schemes.
\item Provide a systematic analysis by resolving several key problems on generalized Hermite subdivision schemes (including multivariate standard Hermite or Lagrange subdivision schemes as special cases)
    such as convergence and smoothness in \cref{thm:converg}, polynomial reproduction in \cref{thm:sd:sr}, and characterization of generalized Hermite masks in \cref{thm:hsd:mask}.
%As byproducts, our results systematically analyze standard multivariate Hermite subdivision schemes which are rarely studied yet in the literature and even in dimension one, our results significantly extend known results on univariate Hermite and Lagrange subdivision schemes.

\item Prove constructively in \cref{thm:existence} that for any $m\in \NN$ and any ordered multiset $\ind\subseteq \ind_m$, there always exists a generalized Hermite mask $a\in \dlrs{0}{r}{r}$ with $r:=\#\ind$ such that the generalized Hermite subdivision scheme of type $\ind$ with mask $a$ is convergent with limit functions in $\dCH{m}$ and its basis vector function is a spline refinable vector function in $\dCH{m}$, whose integer shifts are linearly independent (and therefore, stable).

\item Characterize interpolatory generalized Hermite subdivision schemes in \cref{thm:int}
    and linear-phase moments in \cref{thm:lpm} for the polynomial-interpolation property.
    These results extend all known interpolatory Hermite/Lagrange subdivision schemes to a wide class of generalized Hermite subdivision schemes (e.g., interpolatory Birkhoff subdivision schemes) with desired interpolation properties. %which are of interest in applied mathematics.

\item Develop new techniques for analyzing convergence of generalized Hermite subdivision schemes of type $\ind$, because some arguments in \cite{han20} for univariate Hermite subdivision schemes (i.e., $\ind=\ind_m$) are no longer applicable even in dimension one for ordered multisets $\ind$.
\end{enumerate}

The structure of the paper is as follows. In \cref{sec:hsdmask} we study in \cref{thm:sd:sr} how the vector subdivision operator acts on vector polynomial spaces and investigate the structure of eigenvalues and eigenvectors of the subdivision operator.
%We shall prove in \cref{thm:sd:sr} that the subdivision operator is invariant under a vector polynomial space if and only if its mask must satisfy certain order of sum rules, see \eqref{sr} for the definition of sum rules.
Built on \cref{thm:sd:sr}, we shall prove \cref{thm:hsd:mask} in \cref{sec:hsdmask} for characterizing masks of convergent generalized Hermite subdivision schemes.
Using the notion of the normal form of a matrix mask and employing \cref{thm:sd:sr,thm:hsd:mask},
we prove \cref{thm:converg} in \cref{sec:hsdconverg} for characterizing convergence and smoothness of generalized Hermite subdivision schemes, by linking them with vector cascade algorithms and refinable vector functions.
In \cref{sec:hsdint} we prove  \cref{thm:int} for characterizing interpolatory generalized Hermite subdivision schemes. Our result in \cref{thm:int} not only covers all known interpolatory Hermite subdivision schemes as special cases but also extends to new types of interpolatory subdivision schemes. Then we introduce the notion of linear-phase moments for the polynomial-interpolation property in \cref{thm:lpm}.
In \cref{sec:ex} we prove that
every multivariate subdivision scheme with a matrix mask can be turned into a generalized Hermite subdivision scheme and even into a standard Hermite subdivision scheme, as long as their masks have the same multiplicity. Consequently,
for any $m\in \NN$ and ordered multiset $\ind\subseteq \ind_m$,
we constructively prove in \cref{thm:existence} that there always exists a convergent generalized Hermite subdivision scheme of type $\ind$ whose limit functions are spline functions in $\dCH{m}$ and have linearly independent integer shifts.
Finally, we provide several examples of generalized Hermite subdivision schemes with the symmetry property to illustrate the theoretical results in this paper.

\section{Characterize Masks of Generalized Hermite Subdivision Schemes}
\label{sec:hsdmask}

In this section, we shall prove \cref{thm:hsd:mask} for characterizing the matrix mask of a convergent generalized Hermite
subdivision scheme, by connecting it to refinable vector functions in \eqref{refeq:phi} and vector subdivision operators in \eqref{sd} acting on vector polynomials.

Let us first discuss how a vector subdivision operator acting on vector polynomials.
For $m\in \NN$ and $r\in \N$,
%by $\PL_m$ we denote the space of all $d$-variate polynomials of (total) degree no more than $m$. By
by $(\PL_m)^{1\times r}$ we denote the space of all $1\times r$ (row) vectors of $d$-variate polynomials of (total) degree no more than $m$. Note that a polynomial $\pp$ is uniquely determined by its restriction $\pp|_{\dZ}=\{\pp(k)\}_{k\in \dZ}$ and hence throughout the paper we do not distinguish a polynomial defined on $\dR$ or regarded as a polynomial sequence restricted on $\dZ$.
For $\mu\in \dNN$, recall that $\pp^{(\mu)}(x):=(\partial^\mu \pp)(x)$, the $\mu$th partial derivative of $\pp$.
By \cite[Lemma~3.1]{han13},
the convolution
$\pp*v$ of a polynomial $\pp$ and a sequence $v\in \dlrs{0}{1}{r}$
is always a vector polynomial and is given by
\be \label{conv:poly}
\pp*v:=\sum_{k\in \dZ} \pp(\cdot-k) v(k)=
\sum_{\mu\in \dNN} \frac{(-i)^{|\mu|}}{\mu!} \pp^{(\mu)}(\cdot) \wh{v}^{(\mu)}(0).
\ee
%
%where $\wh{v}(\xi):=\sum_{k\in \dZ} v(k) e^{-ik\cdot \xi}$ is the Fourier series of $v$, which is a row vector of $2\pi \dZ$-periodic trigonometric polynomials.
%
%For $u,v\in \dlrs{0}{1}{r}$ satisfying $\wh{u}(\xi)=\wh{v}(\xi)+\bo(\|\xi\|^{m+1})$ as $\xi \to 0$, we observe from \eqref{conv:poly} and \eqref{bo} that $\pp*u=\pp*v$ for all $\pp\in \PL_m$.
Throughout the paper, we shall use the following notation:
\be \label{Gamma}
\Gamma:=[0,1]^d\cap \dZ=\{ (k_1,\ldots,k_d)^\tp \in \dZ \setsp
k_1,\ldots,k_d\in \{0,1\}\}.
\ee
That is, $\Gamma$ is a particular complete set of all distinct representatives of cosets in the quotient group $\dZ/(2\dZ)$. For $\gamma\in \dZ$, the $\gamma$-coset $a^{[\gamma]}$ of a mask $a$ is defined to be
$a^{[\gamma]}(k):=a(\gamma+2k)$ for $k\in \dZ$.
%Since $\wh{a^{[\gamma]}}(\xi)=\sum_{k\in \dZ} a(\gamma+2k)e^{-ik\cdot \xi}$, it is easy to verify that
%\[
%\wh{a}(\xi)=\sum_{\gamma\in \Gamma} \sum_{k\in \dZ} a(\gamma+2k)e^{-i(\gamma+2k)\cdot \xi}
%=\sum_{\gamma\in \Gamma} e^{-i\gamma\cdot \xi}\wh{a^{[\gamma]}}(2\xi).
%\]

Using a similar idea as in \cite[Proposition~2.2]{han01} and \cite[Theorem~3.4]{han13}, we have
the following result on how a subdivision operator acts on a single vector polynomial.

\begin{lemma}\label{lem:sd:poly}
Let $a\in \dlrs{0}{r}{r}$ be a finitely supported matrix mask and
let $\vec{\pp}=[\pp_1,\ldots,\pp_r]$ be a row vector of $d$-variate polynomials.
Then the following statements are equivalent to each other:
\begin{enumerate}
\item[(1)] $\sd_a \vec{\pp}$ is a vector polynomial sequence, i.e., all the entries of $\sd_a \vec{\pp}$ are polynomial sequences.
\item[(2)] $\sum_{k\in \dZ} \vec{\pp}^{(\mu)}(-2^{-1}\gamma-k)a(\gamma+2k)=\sum_{k\in \dZ} \vec{\pp}^{(\mu)}(-k) a(2k)$ for all $\mu\in \dNN$ and $\gamma\in \Gamma$.
\item[(3)] $[\vec{\pp}^{(\mu)}(-2^{-1}\gamma-i\partial) \wh{a^{[\gamma]}}(\xi)]|_{\xi=0}=[\vec{\pp}^{(\mu)}(-i\partial) \wh{a^{[0]}}(\xi)]|_{\xi=0}$ for all $\mu\in \dNN$ and $\gamma\in \Gamma$, where we define
$\vec{\pp}(x-i\partial) f(\xi):=
\sum_{\mu\in \dNN}
\frac{(-i)^{|\mu|}}{\mu!} \vec{\pp}^{(\mu)}(x) f^{(\mu)}(\xi)$ with $\partial:=(\partial_1,\ldots,\partial_d)^\tp$ only acting on $\xi$.

\item[(4)] $[\vec{\pp}^{(\mu)} (-i\partial)(e^{-i2^{-1}\gamma\cdot\xi} \wh{a^{[\gamma]}}(\xi))]|_{\xi=0}=
    [\vec{\pp}^{(\mu)} (-i\partial) \wh{a^{[0]}}(\xi))]|_{\xi=0}$ for all $\mu\in \dNN$ and $\gamma\in \Gamma$.

\item[(5)] $[\vec{\pp}^{(\mu)}(-i\partial) \wh{a}(2^{-1}\xi+\pi \omega)]|_{\xi=0}=0$ for all $\mu\in \dNN$ and $\omega\in \Gamma\bs\{0\}$.
\end{enumerate}
Moreover, the above items (1)--(5) imply $\deg(\sd_a \vec{\pp})\le \deg(\vec{\pp}):=\max(\deg(\pp_1),\ldots,\deg(\pp_r))$ and
\be \label{sd:poly:single}
[\sd_a \vec{\pp}]^{(\mu)}=
2^{-|\mu|} \sd_a (\vec{\pp}^{(\mu)})=
2^{-|\mu|}(\vec{\pp}^{(\mu)}(2^{-1}\cdot))*a,\qquad \forall\; \mu\in \dNN.
\ee
%
%and consequently $\sd_a (\vec{\pp}(\cdot-t))=(\vec{\pp}(2^{-1}\cdot-t))*a=[\sd_a \vec{\pp}](\cdot-2t)$ for all $t\in \dR$.
\end{lemma}

For every univariate vector polynomial $\vec{\pp}\in (\PL_m)^{1\times r}$, one can show (e.g., see \cite[Lemma~2.2]{han20}) that $\vec{\pp}=\frac{(\cdot)^m}{m!}*v$ for some $v\in \lrs{0}{1}{r}$. Moreover, from \eqref{conv:poly}, we see that $\wh{v}(0)\ne 0$ if and only if $\deg(\vec{\pp})=m$. For the multivariate case, things are more complicated and we have

\begin{lemma}\label{lem:p:v}
For a vector $d$-variate polynomial $\vec{\pp}$ and a $d$-variate polynomial $\pq$, there exists  a finitely supported sequence $v\in \dlrs{0}{1}{r}$ satisfying $\vec{\pp}=\pq*v$ if and only if
\be \label{poly:v:cond}
\{\pp_1,\ldots,\pp_r\}\subseteq
\mspan\{\pq^{(\mu)} \setsp \mu\in \dNN\}
\quad \mbox{with}\quad [\pp_1,\ldots,\pp_r]:=\vec{\pp}.
\ee
%
%Consequently, there always exist a $d$-variate polynomial $\pq$ and a sequence $v\in \dlrs{0}{1}{r}$ such that $\vec{\pp}=\pq*v$.
%In addition, $\wh{v}(0)\ne 0$ if and only if $\deg(\pq)=\deg(\vec{\pp}):=\max(\deg(\pp_1),\ldots,\deg(\pp_r))$.
\end{lemma}

\bp If $\vec{\pp}=\pq*v$, then it follows from \eqref{conv:poly} that
\eqref{poly:v:cond} must hold.
Conversely, we deduce from \eqref{poly:v:cond} that
$\vec{\pp}=\sum_{\mu\in \dNN, |\mu|\le \deg(\pq)} c_\mu \pq^{(\mu)}$ for some $c_\mu\in \C^{1\times r}$. Note that there always exists $v\in \dlrs{0}{1}{r}$ such that
$\wh{v}^{(\mu)}(0)=\frac{\mu!}{(-i)^{|\mu|}} c_\mu$ for all $\mu\in \dNN$ with $|\mu|\le \deg(\pq)$.
It follows from \eqref{conv:poly} that
$\pq*v=\vec{\pp}$.
\ep

For a general vector polynomial $\vec{\pp}$,
\eqref{poly:v:cond} is satisfied by taking $\pq(x_1,\ldots,x_d)=x_1^{m}\cdots x_d^{m}$ with $m:=\deg(\vec{\pp})$. Hence, \cref{lem:p:v} tells us that $\vec{\pp}=\pq*v$ for some $v\in \dlrs{0}{1}{r}$.
However, $\wh{v}(0)\ne 0$ cannot be guaranteed for $d\ge 2$.
Consider $d=2$ and $\vec{\pp}=[x,y]$. Suppose $\vec{\pp}=\pq*v$ with $\wh{v}(0)\ne 0$. Then $\deg(\pq)=\deg(\vec{\pp})=1$ and we can write $\pq(x,y)=c_0+c_1x+c_2y$. By
\eqref{conv:poly}, we must have
\[
[x,y]=\vec{\pp}(x,y)=(\pq*v)(x,y)=(c_0+c_1x+c_2y)\wh{v}(0)-ic_1 \wh{v}^{(1,0)}(0)-ic_2 \wh{v}^{(0,1)}(0),
\]
which forces $c_1 \wh{v}(0)=[1,0]$ and $c_2\wh{v}(0)=[0,1]$, leading to a contradiction $[\frac{1}{c_1},0]=\wh{v}(0)=[0,\frac{1}{c_2}]$.

The condition $\wh{v}(0)\ne 0$ is important for studying generalized Hermite subdivision schemes by using special vector polynomials $\pq*v$ with $\pq\in \PL_m$ and $v\in \lrs{}{1}{r}$.
For $m\in \NN$ and $v\in \dlrs{0}{1}{r}$, following \cite{han01,han03,hanbook}, we define a vector polynomial subspace $\PR_{m,v}$ of $(\PL_m)^{1\times r}$ by
\be \label{Pvm}
\PR_{m,v}:=\{\pp*v \setsp \pp\in \PL_m\},\qquad m\in \NN, v\in \dlrs{0}{1}{r}.
\ee
Note that $\PR_{m,v}\subseteq (\PL_m)^{1\times r}$.
Then $\wh{v}(0)\ne 0$ if and only if $\dim(\PR_{m,v})=\dim(\PL_m)$.
The polynomial space $\PR_{m,v}$ is closely linked to the notion of sum rules and serves as the same role for generalized Hermite subdivision schemes as the polynomial space $\PL_m$ for scalar subdivision schemes.

The relations between subdivision operators and sum rules have been initially studied in \cite{han01,han03,han13} and \cite[Theorem~1.1]{han20} for special cases.  Now we have the following result providing a complete picture on relations between sum rules and vector subdivision operators acting on vector polynomial spaces.

\begin{theorem}\label{thm:sd:sr}
Let $r\in \N$ and $a\in \dlrs{0}{r}{r}$ be a finitely supported matrix mask and $v\in \dlrs{0}{1}{r}$ be a finitely supported sequence with $\wh{v}(0)\ne 0$. For $m\in \NN$, the following are equivalent to each other:
\begin{enumerate}
\item[(1)] $\sd_a \PR_{m,v}=\PR_{m,v}$ and $\sd_a (\wh{v}(0))=\wh{v}(0)$, where $\wh{v}(0)$ is regarded as a constant sequence on $\dZ$.

\item[(2)] There exists a finitely supported sequence $c\in \dlp{0}$ with $\wh{c}(0)=1$ such that
\be \label{sr:2}
\wh{v}(2\xi)\wh{a}(\xi)=\wh{c}(\xi)\wh{v}(\xi)+\bo(\|\xi\|^{m+1})
\quad \mbox{and}\quad
\wh{v}(2\xi)\wh{a}(\xi+\pi \omega)=\bo(\|\xi\|^{m+1}),\quad \xi\to 0,
\omega\in \Gamma\bs \{0\}.
\ee

\item[(3)] There exists a finitely supported sequence $b\in \dlp{0}$ with $\wh{b}(0)=1$ such that the sequence $\vgu_a:=b*v$ satisfies \eqref{sr} for order $m+1$ sum rules with the matching filter $\vgu_a\in \dlrs{0}{1}{r}$.

\item[(4)] $\sd_a \pp_\mu=2^{-|\mu|}\pp_\mu$ for all $\mu\in \dNN$ with $|\mu|\le m$, where
\be \label{pmu}
\pp_\mu(x):=\left(\frac{(\cdot)^\mu}{\mu!}*\vgu_a\right)(x)=
\sum_{0\le \nu\le \mu} \frac{(-i)^{|\nu|}}{\nu! (\mu-\nu)!}
x^{\mu-\nu} \wh{\vgu_a}^{(\nu)}(0),\qquad \mu\in \dNN,
\ee
where $\vgu_a$ is given in item (3) and $\nu\le \mu$ means that all entries of $\mu-\nu$ are nonnegative.

\item[(5)] For every $\mu\in \dNN$ with $|\mu|=m$,
there exist $\pq_\mu\in (\PL_m)^{1\times r}$ such that $\sd_a \pq_\mu=2^{-|\mu|}\pq_\mu$ and
\be \label{poly:cond}
\mspan\{ \pq_\mu^{(\nu)} \setsp \mu,\nu \in \dNN, |\mu|=m \}=\PR_{m,v}.
\ee
\item[(6)] $\sd_a \wh{v}(0)=\wh{v}(0)$ and $\sd_a(\frac{(\cdot)^\mu}{\mu!}*v) \in \PR_{m,v}$ for all $\mu\in \dNN$ with $|\mu|=m$.
\end{enumerate}
Moreover, any of the above items (1)--(6) implies that for all $\pp\in \PL_m$,
\be \label{sd:poly:vva}
\sd_a (\pp*v)=((\pp*v)(2^{-1}\cdot))*a=(\pp(2^{-1}\cdot))*\mathring{v},
\quad
\sd_a (\pp*\vgu_a) =((\pp*\vgu_a)(2^{-1}\cdot))*a=(\pp(2^{-1}\cdot))*\vgu_a,
\ee
with $\wh{\mathring{v}}(\xi):=\wh{v}(2\xi)\wh{a}(\xi)$, and
\be \label{sd:poly:v:deriv}
[\sd_a(\pp*v)]^{(\nu)}
=2^{-|\nu|} \sd_a (\pp^{(\nu)}*v),\quad
[\sd_a(\pp*\vgu_a)]^{(\nu)}
=2^{-|\nu|} \sd_a (\pp^{(\nu)}*\vgu_a),\qquad \forall\; \pp\in \PL_m, \nu\in \dNN.
\ee
\end{theorem}

\bp
(1)$\imply$(2). For any smooth function $f(\xi)$, we have
\be \label{diff:qv}
[(\pq*v)(x-i\partial) f(\xi)]|_{\xi=0}=[\pq(x-i\partial) ( \wh{v}(\xi) f(\xi))]|_{\xi=0},
\ee
where $\partial:=(\partial_1,\ldots,\partial_d)^\tp$ only acts on the variable $\xi$.
Since $\sd_a \PR_{m,v}\subseteq \PR_{m,v}$,
for every $\pq\in \PL_m$, obviously $\vec{\pp}:=\pq*v\in \PR_{m,v}$ and we deduce from \eqref{diff:qv} and \cref{lem:sd:poly} that
\[
[\pq^{(\mu)}(-i\partial)(\wh{v}(\xi)\wh{a}(2^{-1}\xi+\pi \omega))] |_{\xi=0}=
[\vec{\pp}^{(\mu)}(-i\partial) \wh{a}(2^{-1}\xi+\pi \omega)]|_{\xi=0}=0,\qquad \forall\, \omega\in \Gamma\bs\{0\}
\]
by applying \eqref{diff:qv} with $\pq$ being replaced by $\pq^{(\mu)}$.
Consequently, the second identity in \eqref{sr:2} must hold, because the above identity holds for all $\pq\in \PL_m$. Moreover, we must have $\sd_a \PR_{m,v}=\PR_{m,\mathring{v}}$ by \eqref{sd:poly:single} of \cref{lem:sd:poly} and \eqref{diff:qv}, where $\wh{\mathring{v}}(\xi)=\wh{v}(2\xi)\wh{a}(\xi)$.
By \cite[Lemma~3.3]{han03} (also c.f. \cref{lem:p:v}) and $\PR_{m,\mathring{v}}=\sd_a \PR_{m,v}=\PR_{m,v}$, there must exist $c\in \dlp{0}$ such that $\wh{c}(0)\ne 0$ and $\wh{\mathring{v}}(\xi)=\wh{c}(\xi)\wh{v}(\xi)+\bo(\|\xi\|^{m+1})$ as $\xi \to 0$, which proves the first identity in \eqref{sr:2}.
Note that $\wh{v}(0)=1*v\in \PR_{m,v}$. By assumption $\sd_a \wh{v}(0)=\wh{v}(0)$ in item (1), it follows from $\sd_a \wh{v}(0)=\sd_a(1*v)=1*\mathring{v}= \wh{\mathring{v}}(0)$ that $\wh{\mathring{v}}(0)=\wh{v}(0)$. Since $\wh{\mathring{v}}(0)= \wh{c}(0)\wh{v}(0)$ and $\wh{\mathring{v}}(0)=\wh{v}(0)\ne 0$, we must have $\wh{c}(0)=1$.
This proves (1)$\imply$(2).

(2)$\imply$(3). Since $\wh{c}(0)=1$, we can define $\wh{\varphi}(\xi):=\prod_{j=1}^\infty \wh{c}(2^{-j}\xi)$ for $\xi\in \dR$, which satisfies $\wh{\varphi}(2\xi)=\wh{c}(\xi)\wh{\varphi}(\xi)$ and $\wh{\varphi}(0)=1$. Take $b\in \dlp{0}$ such that $\wh{b}(\xi)=1/\wh{\varphi}(\xi)+\bo(\|\xi\|^{m+1})$ as $\xi \to 0$. Then $\frac{\wh{b}(\xi)}{\wh{b}(2\xi)}=\wh{c}(\xi)+\bo(\|\xi\|^{m+1})$ as $\xi \to 0$ and $\wh{b}(0)=1$. Now it follows directly from \eqref{sr:2} that \eqref{sr} holds with $\wh{\vgu_a}(\xi)=\wh{b}(\xi)\wh{v}(\xi)$.
This proves (2)$\imply$(3).

(3)$\imply$(4). By \eqref{sd:poly:single} of \cref{lem:sd:poly} and \eqref{diff:qv}, \eqref{sr} implies
$\sd_a (\frac{(\cdot)^\mu}{\mu!}*\vgu_a)\in \PR_{m,\vgu_a}$ for all $|\mu|\le m$ and
\[
\sd_a \pp_\mu=\sd_a \left(\frac{(\cdot)^\mu}{\mu!} *\vgu_a\right)
=\frac{(2^{-1}\cdot)^\mu}{\mu!} *\vgu_a=2^{-|\mu|} \frac{(\cdot)^\mu}{\mu!} *\vgu_a
=2^{-|\mu|} \pp_\mu.
\]
This proves (3)$\imply$(4).
Moreover, \eqref{sd:poly:vva} and \eqref{sd:poly:v:deriv} follow directly from \eqref{sd:poly:single} and \eqref{diff:qv}.

(4)$\imply$(5). Take $\pq_\mu=\frac{(\cdot)^\mu}{\mu!}*\vgu_a$ for $\mu\in \dNN$ with $|\mu|=m$, which satisfies $\sd_a \pq_\mu=2^{-|\mu|} \pq_\mu$ by item (4). By the definition of $\PR_{m,\vgu_a}$ in \eqref{Pvm}, we trivially have $\mspan\{ \pq_\mu^{(\nu)} \setsp \mu,\nu \in \dNN, |\mu|=m \}=\PR_{m,\vgu_a}$, due to $\pq_\mu^{(\nu)}=\frac{(\cdot)^{\mu-\nu}}{(\mu-\nu)!}*\vgu_a$ for all $0\le \nu\le \mu$. Now we conclude that (4)$\imply$(5) by noting $\PR_{m,\vgu_a}=\PR_{m,v}$, due to $\wh{\vgu_a}(\xi)=\wh{b}(\xi)\wh{v}(\xi)$ and $\wh{b}(0)\ne 0$.

(5)$\imply$(6). By item (5) and \eqref{sd:poly:single}, we have $\sd_a (\pq_\mu^{(\nu)})=2^{|\nu|} [\sd_a \pq_\mu]^{(\nu)}=2^{|\nu|-|\mu|} \pq_\mu^{(\nu)}\in \PR_{m,v}$.
For $\mu\in \dNN$ with $|\mu|=m$,
noting that $\pp_\mu:=\frac{(\cdot)^\mu}{\mu!}*v \in \PR_{m,v}$ and using \eqref{poly:cond}, we can write
$\pp_\mu=\sum_{|\eta|=m, \nu\in \dNN}
c_{\eta,\nu} \pq_\eta^{(\nu)}$
for some coefficients $c_{\eta,\nu}$.
Consequently, we have $\sd_a \pp_\mu =
\sum_{|\eta|=m, \nu\in \dNN}
c_{\eta,\nu} \sd_a (\pq_\eta^{(\nu)})\in \PR_{m,v}$.
This proves $\sd_a \pp_\mu\in \PR_{m,v}$ for all $|\mu|=m$.
On the other hand, we must have $\deg(\pq_\mu)=m$ for some $|\mu|=m$. Suppose not.
Then $\deg(\pq_\mu)<m$ for all $|\mu|=m$ and consequently, we conclude from \eqref{poly:cond} that any vector polynomial in $\PR_{m,v}$ has degree less than $m$. This contradicts our assumption $\wh{v}(0)\ne 0$, which implies $\deg(\frac{(\cdot)^\mu}{\mu!}*v)=m$ for all $|\mu|=m$.
Thus, there must exist $\pq_\mu$ in item (5) satisfying $\deg(\pq_\mu)=m$ for some $|\mu|=m$. Then there exists $\nu\in \dNN$ with $|\nu|=m$ such that $\pp_\mu^{(\nu)}$ is a nonzero constant vector and $\pq_\mu^{(\nu)}\in \PR_{m,v}$.
Since $\wh{v}(0)\ne 0$, we observe from \eqref{conv:poly} that any constant vectors in $\PR_{m,v}$ must be a multiple of $\wh{v}(0)$. Hence, $\pq_\mu^{(\nu)}=\alpha \wh{v}(0)$ for some $\alpha\ne 0$. Now by \eqref{sd:poly:single}, we have
\[
\alpha \sd_a (\wh{v}(0))=
\sd_a (\pp_\mu^{(\nu)})
=2^{|\nu|} [\sd_a \pp_\mu]^{(\nu)}=
2^{|\nu|-|\mu|} \pp_\mu^{(\nu)}=
\alpha \wh{v}(0).
\]
Since $\alpha\ne 0$, the above identity proves $\sd_a (\wh{v}(0))=\wh{v}(0)$. This proves (5)$\imply$(6).

(6)$\imply$(1). Define $\pp_\mu:=\frac{(\cdot)^\mu}{\mu!}*v$ with $\mu\in \dNN$ and $|\mu|=m$.
By \eqref{sd:poly:single}, we have $\sd_a \pp_\mu^{(\nu)} \in \PR_{m,v}$ for all $\nu\in \dNN$ and $|\mu|=m$. This proves $\sd_a \PR_{m,v}\subseteq \PR_{m,v}$. By $\sd_a (\wh{v}(0))=\wh{v}(0)$ and \eqref{conv:poly}, we can easily prove that $\sd_a: \PR_{m,v}\rightarrow \PR_{m,v}$ must be injective.
Consequently, since $\PR_{m,v}$ is a finite dimensional space, we must have $\sd_a \PR_{m,v}=\PR_{m,v}$. This proves (6)$\imply$(1).
\ep

Note that $\lim_{|\xi|\to \infty} \wh{f}(\xi)=0$ for $f\in \dLp{1}$ by the Riemann-Lebesgue lemma.
Using \cref{thm:sd:sr} and linking a generalized Hermite subdivision scheme with a refinable vector function satisfying the refinement equation in \eqref{refeq:phi},
we now prove \cref{thm:hsd:mask} for characterizing a matrix mask and its associated refinable vector function of a convergent generalized Hermite subdivision scheme.

\bp[Proof of \cref{thm:hsd:mask}]
In terms of Fourier series, \eqref{hsd:wn} can be equivalently rewritten as
\be \label{hsd:wn:F}
\wh{w_n}(\xi)=2^{dn} \wh{w_0}(2^n\xi)\wh{a_n}(\xi) \sD_\ind^{-n}
\quad \mbox{with}\quad
\wh{a_n}(\xi):=\wh{a}(2^{n-1}\xi)\cdots \wh{a}(2\xi)\wh{a}(\xi).
\ee
Since the mask $a$ is finitely supported, without loss of generality, we assume that the mask $a$ is supported inside $[-N,N]^d$ for some $N\in \N$.
Now by the definition of the subdivision operator $\sd_a$ in \eqref{sd}, the sequence $a_n$ must be supported inside $[-2^nN, 2^nN]^d$ for all $n\in \N$.
Since $\wh{w_0}(\xi)=I_r$ by $w_0=\td I_r$, the above identity in \eqref{hsd:wn:F} becomes $\wh{w_n}(\xi)=2^{dn} \wh{a_n}(\xi) \sD_\ind^{-n}$, that is, the refinement data are given by $w_n=2^{dn} a_n \sD_\ind^{-n}$ for $n\in \N$.
From the definition $\sD_\ind=\mbox{diag}(2^{-|\nu_1|},\ldots, 2^{-|\nu_r|})$  in \eqref{hsd:wn}, we have $\sD_\ind^{-1} e_1=e_1$ by $\nu_1=\mathbf{0}$. Hence, by $w_n=2^{dn} a_n \sD_\ind^{-n}$, \eqref{hsd:converg} with $K=N$ becomes \eqref{def:phi}, that is,
\be \label{def:phi:2}
\lim_{n\to \infty} \sup_{k\in \dZ} \| 2^{dn} a_n(k)e_1-\phi(2^{-n}k)\|=0.
\ee
Because $a_n$ is supported inside $[-2^nN,2^n N]^d$, it follows directly from \eqref{def:phi:2} that  $\phi$ must be supported inside $[-N,N]^d$.
Noting that $\phi$ is continuous and using Riemann sums and \eqref{def:phi:2}, we have
\begin{align*}
\wh{\phi}(\xi)
&=\int_{\dR} \phi(x) e^{-ix \cdot \xi} dx=
\lim_{n\to \infty}
2^{-dn}\sum_{k\in \dZ\cap [-2^nN,2^nN]^d} \phi(2^{-n}k)e^{-i 2^{-n}k \cdot \xi}\\
&=\lim_{n\to \infty}
\sum_{k\in \dZ\cap [-2^nN,2^nN]^d} a_n(k) e^{-i 2^{-n}k \cdot \xi} e_1
=\lim_{n\to \infty}
\wh{a_n}(2^{-n}\xi) e_1=\lim_{n\to \infty} \Big(\Big [\prod_{j=1}^n \wh{a}(2^{-j}\xi)\Big]e_1\Big),
\end{align*}
where we used the identity
$\sum_{k\in \dZ} a_n(k) e^{-i2^{-n} k \cdot \xi}=\wh{a_n}(2^{-n}\xi)=\prod_{j=1}^n a(2^{-j}\xi)$ and the inequality
\begin{align*}
&2^{-dn}\sum_{k\in \dZ\cap [-2^nN,2^nN]^d} \| \phi(2^{-n}k)e^{-i 2^{-n}k \cdot \xi}
-2^{dn} a_n(k) e^{-i 2^{-n}k \cdot \xi} e_1\|\\
&\qquad \le 2^{-dn} (2^{n+1}N+1)^d  \sup_{k\in \dZ}
\|\phi(2^{-n}k)-2^{dn} a_n(k) e_1\|
\le (2N+1)^d \sup_{k\in \dZ}
\|\phi(2^{-n}k)-2^{dn} a_n(k) e_1\|,
\end{align*}
which goes to $0$ as $n\to \infty$ by \eqref{def:phi:2}. This proves the first part of \eqref{refeq}, from which we get
\[
\wh{\phi}(2\xi)=
\lim_{n\to \infty} \Big(\Big [\prod_{j=1}^n \wh{a}(2^{1-j}\xi)\Big]e_1\Big)
=\wh{a}(\xi) \lim_{n\to \infty}
\Big(\Big [\prod_{j=1}^{n-1} \wh{a}(2^{-j}\xi)\Big]e_1\Big)
=\wh{a}(\xi)\wh{\phi}(\xi).
\]
That is, $\phi=2^d \sum_{k\in \dZ} a(k) \phi(2\cdot-k)$.
This proves the second part of \eqref{refeq}. Hence, we proved \eqref{refeq}.

We now prove $\wh{\phi}(0)\ne 0$ and items (1)--(4) under the condition in \eqref{phi:cond0}. We first prove $\wh{\phi}(0)\ne 0$.
Note that $\phi\in (\dCH{m})^r$ has compact support. For $\mu\in \ind_{m}$, we have $\phi^{(\mu)}\in \dLp{1}$ and $\wh{\phi^{(\mu)}}(\xi)=(i\xi)^\mu \wh{\phi}(\xi)$. Since \eqref{refeq} implies $\wh{\phi}(2^n \xi)=\wh{a_n}(\xi)\wh{\phi}(\xi)$ and
$\wh{a_n}$ is $2\pi\dZ$-periodic, we obtain
\[
[\wh{a}(0)]^n \wh{\phi}(2\pi k)
=\wh{a_n}(0)\wh{\phi}(2\pi k)=
\wh{a_n}(2\pi k)\wh{\phi}(2\pi k)=\wh{\phi}(2^{n+1}\pi k),\qquad k\in \dZ.
\]
Applying the Riemann-Lebesgue lemma to $\wh{\phi^{(\mu)}}(2^{n+1}\pi k)$, we have
\[
\lim_{n\to \infty}  (i2^{n+1}\pi k)^\mu [\wh{a}(0)]^n \wh{\phi}(2\pi k)=
\lim_{n\to \infty} \wh{\phi^{(\mu)}}(2^{n+1}\pi k)=0,\qquad \forall\, k\in \dZ\bs\{0\}, \mu\in \ind_{m}.
\]
Since the above identity holds for all $|\mu|=m$, we conclude from the above identity that
\be \label{an:phi}
\lim_{n\to\infty} 2^{mn} [\wh{a}(0)]^n \wh{\phi}(2\pi k)=0,\qquad \forall\; k\in \dZ\bs\{0\}.
\ee
We use proof by contradiction to show $\wh{\phi}(0)\ne 0$. Suppose not. Then
$\wh{\phi}(0)=0$. From \eqref{an:phi} and \eqref{phi:cond0}  with $\omega=0$, we must have
$\lim_{n\to \infty} 2^{mn} [\wh{a}(0)]^n=0$.
In particular, $\lim_{n\to \infty} [\wh{a}(0)]^n=0$ due to $m\ge 0$. However, by $\wh{\phi}(\xi)=\lim_{n\to \infty} ([\prod_{j=1}^{n} \wh{a}(2^{-j}\xi)]e_1)$ in \eqref{refeq} and $\lim_{j\to \infty} \wh{a}(2^{-j}\xi)=\wh{a}(0)$, we must end up with $\wh{\phi}(\xi)=\lim_{n\to \infty}([\prod_{j=1}^{n} \wh{a}(2^{-j}\xi)]e_1)=0$ for every $\xi\in \dR$. That is, $\phi$ must be identically zero, which is a contradiction to our assumption in \eqref{phi:cond0} with $\omega=0$. Hence, $\wh{\phi}(0)\ne 0$.

By $\wh{\phi}(0)=\wh{a}(0)\wh{\phi}(0)$ in \eqref{refeq} and $\wh{\phi}(0)\ne 0$, we see that $1$ must be an eigenvalue of $\wh{a}(0)$. Now by \eqref{phi:cond0} with $\omega=0$ and \eqref{an:phi}, all the other eigenvalues of $\wh{a}(0)$ must be less than $2^{-m}$ in modulus. The existence and uniqueness of the values $\wh{\vgu_a}^{(\mu)}(0)$ for $\mu\in \ind_{m}$ in \eqref{vgua:0} and \eqref{vgua:val} are guaranteed, because $1$ is a simple eigenvalue of $\wh{a}(0)$ and the matrices $I_r-2^{|\mu|} \wh{a}(0)$ are invertible for all $\mu\in \ind_{m}\bs\{0\}$.
This proves item (1).

We now prove item (2).
Applying the Leibniz differentiation formula to the following identity
\be \label{vgua:00}
\wh{\vgu_a}(2\xi)\wh{a}(\xi)=\wh{\vgu_a}(\xi)+\bo(\|\xi\|^{m+1}),\qquad \xi\to 0,
\ee
we observe that
the above identity together with $\wh{\vgu_a}(0)\wh{\phi}(0)=1$ is equivalent to \eqref{vgua:0} and \eqref{vgua:val}.
Consequently, \eqref{vgua:0} and \eqref{vgua:val} together imply \eqref{vgua:00}, which is just the first part of \eqref{sr}.
Define $g:=\vgu_a*\phi=\sum_{k\in \dZ} \vgu_a(k) \phi(\cdot-k)$ and $g_\mu:=(-ix)^\mu g(x)$ for $\mu\in \ind_{m}$. Then both $g$ and $g_\mu$ are compactly supported functions in $\dCH{m}$.
Note that $\wh{g_\mu^{(\nu)}}(\xi)=(i\xi)^\nu \wh{g_\mu}(\xi)=(i\xi)^\nu \wh{g}^{(\mu)}(\xi)$. Since $g_\mu^{(\nu)}\in \dLp{1}$ by $g_\mu \in \dCH{m}$ and $|\nu|\le m$, we conclude from the Riemann-Lebesgue lemma that
\[
\lim_{n\to\infty} 2^{n|\nu|}(i\xi)^\nu \wh{g}^{(\mu)}(2^n\xi)=\lim_{n\to \infty} \wh{g_\mu^{(\nu)}}(2^n\xi)=0
\qquad \forall\, \xi\in \dR\bs\{0\}, \mu,\nu\in \ind_{m}.
\]
Considering $\nu\in \ind_{m}$ with $|\nu|=m$ in the above identity, we have
\be \label{fto0}
\lim_{n\to \infty} 2^{mn} \wh{g}^{(\mu)}(2^n 2\pi k)=0 \qquad
\forall\; \mu\in \ind_{m}, k\in \dZ\bs\{0\}
\quad \mbox{with}\quad g=\vgu_a*\phi.
\ee
Using $\wh{g}(\xi)=\wh{\vgu_a}(\xi)\wh{\phi}(\xi)
=\wh{\vgu_a}(\xi) \wh{a_n}(2^{-n}\xi)\wh{\phi}(2^{-n}\xi)$ and \eqref{vgua:00}, we deduce that for $k\in \dZ$,
\begin{align*}
\wh{g}(2^n(\xi+2\pi k))
&=\wh{\vgu_a}(2^n(\xi+2\pi k))\wh{a_n}(\xi+2\pi k)\wh{\phi}(\xi+2\pi k)\\
&=\wh{\vgu_a}(2^n\xi)\wh{a}(2^{n-1}\xi)\wh{a}(2^{n-2}\xi)\cdots \wh{a}(2\xi)\wh{a}(\xi)\wh{\phi}(\xi+2\pi k)\\
&=\wh{\vgu_a}(2^{n-1}\xi)\wh{a}(2^{n-2}\xi)\cdots \wh{a}(2\xi)\wh{a}(\xi)\wh{\phi}(\xi+2\pi k)+\bo(\|\xi\|^{m+1})\\
%&=\wh{\vgu_a}(2\xi)\wh{a}(\xi)\wh{\phi}(\xi+2\pi k)+\bo(\|\xi\|^{m+1})\\
&=\wh{\vgu_a}(\xi)\wh{\phi}(\xi+2\pi k)+\bo(\|\xi\|^{m+1})\\
&=\wh{g}(\xi+2\pi k)+\bo(\|\xi\|^{m+1}),
\end{align*}
as $\xi \to 0$. On one hand, the above identity with $k=0$ becomes $\wh{g}(2^n\xi)=\wh{g}(\xi)+\bo(\|\xi\|^{m+1})$ as $\xi\to 0$, which forces $\wh{g}(\xi)=\wh{g}(0)+\bo(\|\xi\|^{m+1})$ as $\xi \to 0$ by employing the Taylor expansion of $\wh{g}$ at $\xi=0$.
On the other hand,
the above identity
for all  $k\in \dZ\bs\{0\}$ and $\mu\in \ind_{m}$ is equivalent to
\[
\wh{g}^{(\mu)}(2\pi k)=[ \wh{g}(2^n(\xi+2\pi k))]^{(\mu)}|_{\xi=0}=
2^{n|\mu|} \wh{g}^{(\mu)}(2^n2\pi k)\to 0 \quad \mbox{as}\;\; n\to \infty,
\]
where we used \eqref{fto0} in the last identity.
Consequently, we proved $\wh{\vgu_a}(\xi)\wh{\phi}(\xi+2\pi k)=\wh{g}(\xi+2\pi k)=\wh{g}(0)\td(k)+\bo(\|\xi\|^{m+1})$ as $\xi \to 0$. Since $\wh{g}(0)=\wh{\vgu_a}(0)\wh{\phi}(0)=1$,
this proves item (2).

Now we prove item (3).
Since \eqref{vgua:00} holds and is just the first part of \eqref{sr},
it suffices to prove the second part of \eqref{sr}.
For $\omega\in \Gamma\bs\{0\}$ and $k\in \dZ$,
noting that $\wh{a}(\xi+\pi \omega)\wh{\phi}(\xi+\pi \omega+2\pi k)=\wh{a}(\xi+\pi \omega+2\pi k)\wh{\phi}(\xi+\pi \omega+2\pi k)=\wh{\phi}(2\xi+2\pi \omega+4\pi k)$,
we have
\[
\wh{\vgu_a}(2\xi)\wh{a}(\xi+\pi \omega)\wh{\phi}(\xi+\pi \omega+2\pi k)=
\wh{\vgu_a}(2\xi) \wh{\phi}(2\xi+2\pi \omega+4\pi k)=\wh{g}(2\xi+2\pi \omega+4\pi k)=\bo(\|\xi\|^{m+1})
\]
as $\xi \to 0$, where we used item (2) in the last identity and $\omega+2k\in \dZ\bs\{0\}$. By our assumption in \eqref{phi:cond0}, the above identity yields the second part of \eqref{sr}. This shows that the mask $a$ must have order $m+1$ sum rules with a matching filter $\vgu_a$ satisfying \eqref{vgua:0} and \eqref{vgua:val}. We now further show that $\vgu_a$ must satisfy \eqref{vgua:hsd}.
By \eqref{refeq} and \eqref{vgua:0}, we have
\[
1=\wh{\vgu_a}(0)\wh{\phi}(0)=
\lim_{n\to \infty} \wh{\vgu_a}(0) [\wh{a}(0)]^n e_1=\lim_{n\to \infty} \wh{\vgu_a}(0)e_1=\wh{\vgu_a}(0)e_1.
\]
Hence, $\wh{\vgu_a}(0)e_1=1$.
Consider an initial data $w_0=\vec{\pp}_\mu:=\frac{(\cdot)^\mu}{\mu!}*\vgu_a$ with $\mu\in \ind_{m}$. By \cref{thm:sd:sr} and item (3),
we have $\sd_a \vec{\pp}_\mu =2^{-|\mu|}\vec{\pp}_\mu$.
Then $\{w_n\}_{n=1}^\infty$ in \eqref{hsd:wn} must satisfy
\[
w_n:=\sd_a^n (\vec{\pp}_\mu)\sD_\ind^{-n}=
2^{-|\mu|n} \vec{\pp}_\mu \sD_\ind^{-n}.
\]
From the definition of $\sD_\ind$ in \eqref{hsd:wn}, we have
$\sD_\ind^{-n} e_\ell=2^{|\nu_\ell| n} e_\ell$ for $\ell=1,\ldots,r$.
Therefore, $w_n e_\ell=2^{(|\nu_\ell|-|\mu|)n}\vec{\pp}_\mu e_\ell$.
Since the generalized Hermite subdivision scheme of type $\ind$ is convergent with limit functions in $\dCH{m}$, by \cref{def:hsd}, there exists a function $\eta_\mu\in \dCH{m}$ such that \eqref{hsd:converg} holds. That is, for any constant $K>0$,
\be \label{hsd:converg:2}
\lim_{n\to \infty} \sup_{k\in \dZ\cap [-2^n K, 2^n K]^d}
|2^{(|\nu_\ell|-|\mu|)n} \vec{\pp}_\mu(k) e_\ell-\eta^{(\nu_\ell)}_\mu(2^{-n}k)|=0,\qquad \forall\; \mu\in \ind_{m}, \ell=1,\ldots,r.
\ee
By definition, $\vec{\pp}_\mu=\frac{(\cdot)^\mu}{\mu!}*\vgu_a$ and
therefore, using \eqref{conv:poly} we get
\be \label{Hmuelln}
2^{(|\nu_\ell|-|\mu|)n} \vec{\pp}_\mu(x)e_\ell
=
\sum_{0\le \beta \le \mu}
\frac{(-i)^{|\beta|}}{\beta!}
\frac{(2^{-n}x)^{\mu-\beta}}{(\mu-\beta)!}
\wh{\vgu_a}^{(\beta)}(0) e_\ell 2^{(|\nu_\ell|-|\beta|)n}
:=H_{\mu,\ell,n}(2^{-n} x).
\ee
That is, $2^{(|\nu_\ell|-|\mu|)n}\vec{\pp}_\mu (k) e_\ell=H_{\mu,\ell,n}(2^{-n}k)$ for all $k\in \dZ$ and $n\in \N$.
Now \eqref{hsd:converg:2} becomes
\[
\lim_{n\to \infty} \sup_{k\in \dZ\cap [-2^nK,2^nK]^d} |H_{\mu,\ell,n}(2^{-n}k)-\eta_\mu^{(\nu_\ell)}(2^{-n}k)|=0,\qquad \forall\; \mu\in \ind_{m}, \ell=1,\ldots,r.
\]
Since all involved functions are continuous, the above identity implies
$\lim_{n\to \infty} H_{\mu,\ell,n}(x)=\eta^{(\nu_\ell)}_\mu(x)$ on any compact set of $\dR$.
By the above identity and the definition of $H_{\mu,\ell,n}$ in \eqref{Hmuelln}, we have
\be \label{etaH}
\eta^{(\nu_\ell)}_\mu(x)
=\lim_{n\to \infty} H_{\mu,\ell,n}(x)=
\sum_{0\le \beta\le \mu} \frac{(-i)^{|\beta|}}{\beta!}
\frac{x^{\mu-\beta}}{(\mu-\beta)!}
\wh{\vgu_a}^{(\beta)}(0) e_\ell \lim_{n\to \infty}2^{(|\nu_\ell|-|\beta|)n},
\ee
for all $\mu\in \ind_m$ and $\ell=1,\ldots,r$.
Since $\lim_{n\to \infty} 2^{(|\nu_\ell|-|\beta|)n}=\infty$ for $|\beta|<|\nu_\ell|$, the above limit $\lim_{n\to \infty} H_{\mu,\ell,n}(x)$ in \eqref{etaH} exists only if
\[
\wh{\vgu_a}^{(\beta)}(0)
e_\ell=0\quad \mbox{for all}\; \mu\in \ind_{m}, 0\le \beta\le \mu, \ell=1,\ldots,r
\quad \mbox{satisfying}\quad
|\beta|<|\nu_\ell|.
\]
For every $\nu_\ell\in \ind\subseteq\ind_m$ and $\beta\in \ind_m$ satisfying $|\beta|<|\nu_\ell|$, there always exists $\mu\in \ind_{m}$ such that $\beta\le \mu$ (for example, we can simply take $\mu=\beta$).
Hence, we conclude from the above identities that
\be \label{vgua:relation}
\wh{\vgu_a}^{(\beta)}(0)
e_\ell=0,\qquad \forall\; \beta\in \ind_m, |\beta|<|\nu_\ell|, \ell=1,\ldots,r.
\ee
Applying the above identities in \eqref{vgua:relation} to \eqref{etaH} and noting $\lim_{n\to \infty} 2^{(|\nu_\ell|-|\beta|)n}=0$ for $|\beta|>|\nu_\ell|$,
we end up with
\be \label{etamu1}
\eta_\mu^{(\nu_\ell)}(x)
%&=\sum_{0\le \beta\le \mu, |\beta|=|\nu_\ell|} \frac{(-i)^{|\beta|}}{\beta!}
%\frac{x^{\mu-\beta}}{(\mu-\beta)!}
%\wh{\vgu_a}^{(\beta)}(0) e_\ell\\
=\begin{cases}
\sum_{0\le \beta\le \mu, |\beta|=|\nu_\ell|} \frac{(-i)^{|\beta|}}{\beta!} \frac{x^{\mu-\beta}}{(\mu-\beta)!} \wh{\vgu_a}^{(\beta)}(0) e_\ell, &\text{if $|\nu_\ell|\le |\mu|$},\\
0, &\text{if $|\nu_\ell|>|\mu|$}.
\end{cases}
\ee
Taking $\ell=1$ in the above first identity and noting $\nu_1=\mathbf{0}$ by our convention in \eqref{indset},
 we see that
\[
\eta_\mu(x)=\frac{x^\mu}{\mu!} \wh{\vgu_a}(0) e_1=
\frac{x^\mu}{\mu!},\qquad \forall\; \mu\in \ind_{m},
\]
where we used the proved identity $\wh{\vgu_a}(0)e_1=1$.
From the above identity we trivially have
\be \label{etamu2}
\eta_\mu^{(\nu_\ell)}(x)=
\begin{cases}
\frac{x^{\mu-\nu_\ell}}{(\mu-\nu_\ell)!}, &\text{if $\nu_\ell \le \mu$},\\
0, &\text{otherwise},
\end{cases}
\qquad \forall\; \mu\in \ind_{m}, \ell=1,\ldots,r.
\ee
Note that $\nu_\ell \in \ind \subseteq \ind_{m}$ and hence there exists $\mu\in \ind_{m}$ such that $\nu_\ell \le \mu$ (i.e., take $\mu=\nu_\ell$).
Using \eqref{vgua:relation} and comparing the two expressions for $\eta_\mu^{(\nu_\ell)}$ in \eqref{etamu1} and \eqref{etamu2} with all $\mu\in \ind_{m}$, we have
\be \label{vgua:hsd:cond}
\frac{(-i)^{|\beta|}}{\beta!} \wh{\vgu_a}^{(\beta)}(0) e_\ell=
\begin{cases}
1, &\text{if $\beta=\nu_\ell$},\\
0, &\text{if $|\beta|\le |\nu_\ell|, \beta\ne \nu_\ell$, $\beta\in \ind_m$},
\end{cases}\qquad \forall\; \ell=1,\ldots,r,
\ee
which is just equivalent to \eqref{vgua:hsd} by using the definition of the big $\bo$ notation in \eqref{bo}. This proves item (3).

By \cref{thm:sd:sr} and item (3), we must have
$\sd_a \vec{\pp}_\mu=2^{-|\mu|} \vec{\pp}_\mu$.
The identity $\sum_{k\in \dZ} \vec{\pp}_\mu(k)\phi(\cdot-k)=\frac{(\cdot)^\mu}{\mu!}$ follows directly from item (2) (also see \cite[Theorem~2.4]{han01}).
Since $\vec{\pp}_\mu=\frac{(\cdot)^\mu}{\mu!}*\vgu_a$, the expression in item (4) for $\vec{\pp}_\mu (x) e_\ell$ follows
directly from \eqref{conv:poly} and  \eqref{vgua:hsd}.
This proves item (4).
\ep

\section{Convergence and Smoothness of Generalized Hermite Subdivision Schemes}
\label{sec:hsdconverg}

With the help of the normal form of a matrix mask,
in this section we shall prove \cref{thm:converg} for characterizing convergence and smoothness of generalized Hermite subdivision schemes by connecting them to vector cascade algorithms and refinable vector functions.

\subsection{Normal form of a matrix mask}
The study of vector cascade algorithms and refinable vector functions is often much more difficult and complicated than their scalar counterparts (i.e., $r=1$), largely because its underlying mask $a$ is a matrix mask.
Convergence of vector cascade algorithms in Sobolev spaces have been systematically studied in \cite{han03} and many references therein (e.g., \cite{hj98,hj06,jrz98}). In comparison with previous work on vector cascade algorithms in the literature, the approach for studying vector cascade algorithms in \cite{han03} is relatively easy, thanks to the notion of the normal (or canonical) form of a matrix mask.
To introduce the normal form of a matrix mask, let us recall some definitions. A finitely supported sequence $U\in \dlrs{0}{r}{r}$ is said to be \emph{strongly invertible} if $\wh{U}(\xi):=\sum_{k\in \dZ} U(k) e^{-ik\cdot \xi}$ is invertible for every $\xi\in \dR$ and $\wh{U}(\xi)^{-1}$ is an $r\times r$ matrix of $2\pi \dZ$-periodic trigonometric polynomials. Thus, for a strongly invertible $U\in \dlrs{0}{r}{r}$, we can define another sequence $U^{-1}$ through $\wh{U^{-1}}(\xi):=\wh{U}(\xi)^{-1}$. Then $U^{-1}$ must be a finitely supported sequence and $U*U^{-1}=U^{-1}*U=\td I_r$.
The key idea of the normal form is to transform the original matrix mask $a$ and its refinable vector function $\phi$ into a new mask and a new refinable vector function with relatively simple structures.
Using a strongly invertible sequence $U\in \dlrs{0}{r}{r}$, we define
\be \label{newa}
\wh{\mathring{\phi}}(\xi):=\wh{U}(\xi)^{-1}\wh{\phi}(\xi),
\quad \wh{\mathring{a}}(\xi):=\wh{U}(2\xi)^{-1} \wh{a}(\xi) \wh{U}(\xi),
\quad \mbox{and}\quad
\wh{\mathring{\vgu}_a}(\xi):=\wh{\vgu_a}(\xi)\wh{U}(\xi).
\ee
Since $U$ is strongly invertible,
one can easily verify the following claims:
\begin{enumerate}
\item[(1)] The mask $a$ is finitely supported if and only if the new mask $\mathring{a}$ is finitely supported;

\item[(2)] The vector function $\phi$ is a compactly supported refinable vector function associated with the mask $a$ satisfying $\wh{\phi}(2\xi)=\wh{a}(\xi) \wh{\phi}(\xi)$ if and only if the new vector function $\mathring{\phi}$ is a compactly supported refinable vector function associated with the new mask $\mathring{a}$ satisfying $\wh{\mathring{\phi}}(2\xi)=\wh{\mathring{a}}(\xi)
    \wh{\mathring{\phi}}(\xi)$;

\item[(3)] The mask $a$ has order $m+1$ sum rules with a matching filter $\vgu_a\in \dlrs{0}{1}{r}$ if and only if the new mask $\mathring{a}$ has order $m+1$ sum rules with the new matching filter $\mathring{\vgu}_a \in \dlrs{0}{1}{r}$.
\end{enumerate}

The key idea of the normal form of a matrix mask is
to use a suitable strongly invertible sequence $U\in \dlrs{0}{r}{r}$ such that
the new mask $\mathring{a}$ behaves essentially like a scalar mask by employing a matching filter $\mathring{\vgu}_a$
with a simple special structure.
To do so, let us recall a result from \cite[Lemma~2.3 and Theorem~5.1]{han10mc}, also c.f.
\cite[Theorem~2.1]{han09} and
\cite[Proposition~2.4]{han03}. Due to its importance,
we shall sketch the proof given in \cite[Lemma~2.3]{han10mc}.

\begin{lemma}\label{lem:vgu}
Let $r\in \N$ with $r>1$. Let $u, v\in \dlrs{0}{1}{r}$ be two sequences such that $\wh{u}(0)\ne 0$ and $\wh{v}(0)\ne 0$. For each $m\in \NN$, there exists a strongly invertible sequence $U\in \dlrs{0}{r}{r}$ such that
\be \label{uvU}
\wh{u}(\xi)\wh{U}(\xi)=\wh{v}(\xi)+\bo(\|\xi\|^{m+1}),\qquad \xi\to 0.
\ee
\end{lemma}

\bp Let us first prove \eqref{uvU}
for the special case $\wh{v}(\xi)=[1,0,\ldots,0]$.
Write $[u_1,\ldots,u_r]:=u$.
Since $\wh{u}(0)\ne 0$, without loss of generality, we can assume $\wh{u_1}(0)\ne 0$; otherwise, we can perform a permutation on $\wh{u}$ first.
Since $\wh{u_1}(0)\ne 0$,
there always exist sequences $c_1,\ldots,c_r\in \dlp{0}$ satisfying
\[
\wh{c_1}(\xi)=1/\wh{u_1}(\xi)+\bo(\|\xi\|^{m+1})
\quad \mbox{and}\quad
\wh{c_\ell}(\xi)=\wh{u_\ell}(\xi)/\wh{u_1}(\xi)+\bo(\|\xi\|^{m+1}),\qquad \xi\to 0, \ell=2,\ldots,r.
\]
Using binomial expansion and $\wh{c_1}(0)\ne 0$, we can write
$(1-\wh{c_1}(\xi)/\wh{c_1}(0))^{2m+2}=1-\wh{c_1}(\xi) \wh{g}(\xi)$ for a unique $2\pi\dZ$-periodic trigonometric polynomial $g$.
Define two sequences $U_1$ and $U_2$ in $\dlrs{0}{r}{r}$ by
\[
\wh{U_1}(\xi):=\left[ \begin{matrix}
1 &-\wh{c_2} &\cdots &-\wh{c_r}(\xi)\\
0 &1 &\cdots &0\\
\vdots &\vdots &\ddots &\vdots\\
0 &0 &0 &1\end{matrix}\right],\quad
\wh{U_2}(\xi):=\left[ \begin{matrix}
\wh{c_1}(\xi) &-(1-\wh{c_1}(\xi)/\wh{c_1}(0))^{m+1} &0\\
(1-\wh{c_1}(\xi)/\wh{c_1}(0))^{m+1} &g(\xi) &0\\
0 &0 &I_{r-2}\end{matrix}\right].
\]
Note that only the definition of $U_2$ requires the condition $r>1$.
Since $\det(\wh{U_1}(\xi))=\det(\wh{U_2}(\xi))=1$, the sequences $U_1$ and $U_2$ are obviously strongly invertible. Since $1-\wh{c_1}(\xi)/\wh{c_1}(0)=\bo(\|\xi\|)$ as $\xi \to 0$, we can directly verify that
\[
\wh{u}(\xi) \wh{U_1}(\xi)=[\wh{u_1}(\xi),0,\ldots,0]+\bo(\|\xi\|^{m+1})
\quad \mbox{and}\quad
\wh{u}(\xi)\wh{U_1}(\xi)\wh{U_2}(\xi)=
[1,0,\ldots,0]+\bo(\|\xi\|^{m+1})
\]
as $\xi \to 0$. This proves \eqref{uvU} for the special case $\wh{v}=[1,0,\ldots,0]$. For the general case, there exist strongly invertible sequences $U_u, U_v\in \dlrs{0}{r}{r}$ such that $\wh{u}(\xi) \wh{U_u}(\xi)=[1,0,\ldots,0]+\bo(\|\xi\|^{m+1})=
\wh{v}(\xi) \wh{U_v}(\xi)$ as $\xi \to 0$. Then \eqref{uvU} holds by taking $\wh{U}(\xi)=\wh{U_u}(\xi)\wh{U_v}(\xi)^{-1}$.
\ep

From now on we assume that a matrix mask $a\in \dlrs{0}{r}{r}$ has order $m+1$ sum rules with a matching filter $\vgu_a\in \dlrs{0}{1}{r}$ satisfying $\wh{\vgu_a}(0)\ne 0$.
Then \cref{lem:vgu} tells us that there always exists a strongly invertible sequence $U\in \dlrs{0}{r}{r}$ satisfying
\be \label{vgu:normal}
\wh{\vgu_a}(\xi)\wh{U}(\xi)=\wh{\mathring{\vgu}_a}(\xi)
+\bo(\|\xi\|^{m+1}), \quad \xi \to 0
\quad \mbox{with}\quad
\mathring{\vgu}_a:=[c,0,\ldots,0],
\quad c\in \dlp{0}, \wh{c}(0)=1.
\ee
For $r>1$, we can further choose $c=\td$.
Then the new mask $\mathring{a}$ in \eqref{newa} must have order $m+1$ sum rules with the new matching filter $\mathring{\vgu}_a=[c,0,\ldots,0]$. Hence, from the definition of sum rules in \eqref{sr} for mask $\mathring{a}$,
we observe that the new mask $\mathring{a}$ with a special matching filter $\mathring{\vgu}_a$ in \eqref{vgu:normal} must satisfy
\be \label{normalform}
\mathring{a}=
\left[ \begin{matrix} \mathring{a}_{1,1} &\mathring{a}_{1,2}\\
\mathring{a}_{2,1} &\mathring{a}_{2,2}\end{matrix}
\right]
\quad \mbox{with}\quad
\wh{\mathring{a}_{1,1}}(0)=1,\quad
\wh{\mathring{a}_{1,1}}(\xi+ \pi \omega)=\bo(\|\xi\|^{m+1}),\quad \xi\to 0,\; \forall\; \omega\in \Gamma\bs \{0\}
\ee
and
\be \label{a12}
\wh{\mathring{a}_{1,2}}(\xi+\pi \omega)=\bo(\|\xi\|^{m+1}),\qquad \xi \to 0,\; \forall\; \omega\in \Gamma,
\ee
where $\mathring{a}_{1,1}\in \dlp{0}, \mathring{a}_{1,2}\in \dlrs{0}{1}{(r-1)}$, $\mathring{a}_{2,1}\in \dlrs{0}{(r-1)}{r}$ and $\mathring{a}_{2,2}\in \dlrs{0}{(r-1)}{(r-1)}$.
If in addition $c=\td$ in \eqref{vgu:normal}, then we further have $\wh{\mathring{a}_{1,1}}(\xi)=1+\bo(\|\xi\|^{m+1})$ as $\xi \to 0$.
A matrix mask $\mathring{a}$ satisfying \eqref{normalform} and \eqref{a12} is called \emph{a normal (or canonical) form} of the original matrix mask $a$. Obviously, $\mathring{a}_{1,1}$ now behaves exactly like a scalar mask with order $m+1$ sum rules, while $\mathring{a}_{1,2}$ satisfying \eqref{a12} behaves essentially like the trivial zero filter.
The normal form of a matrix mask makes the study of vector cascade algorithms and refinable vector functions relatively easy, since many techniques for studying scalar masks and scalar refinable functions can be applied now without much difficulty to the new mask $\mathring{a}$ satisfying \eqref{normalform} and \eqref{a12}.
For more details,
see \cite[Proposition~2.4]{han03},
\cite[Theorem~2.1]{han09}, and
\cite[Theorem~5.1]{han10mc} for the normal form of matrix masks and its applications.

\subsection{The smoothness quantity $\sm_p(a)$}\label{subsec:sma}
To study the convergence of a generalized Hermite subdivision scheme,
we need to recall a few definitions.
For $m\in \NN$ and $v\in \dlrs{0}{1}{r}$, we define
\be \label{Vvm}
\PV_{m,v}:=\{u\in (\dlp{0})^r \setsp \wh{v}(\xi)\wh{u}(\xi)=\bo(\|\xi\|^{m+1}),\quad \xi\to 0\}.
\ee
In fact, $\PV_{m,v}=\{u \in (\dlp{0})^r \setsp \vec{\pp}*u=0\; \forall\, \vec{\pp}\in \PR_{m,v}\}$,
which is perpendicular to the vector polynomial space $\PR_{m,v}$ in \eqref{Pvm}.
For $v\in \dlrs{0}{s}{r}$ and $k\in \dZ$, we define $\nabla_k v:=v-v(\cdot-k)$. In particular,  we define
\[
\nabla^\mu v:=\nabla_{e_1}^{\mu_1}\cdots \nabla_{e_d}^{\mu_d} v,\qquad \mu=(\mu_1,\ldots,\mu_d)^\tp\in \dNN, v\in \dlrs{0}{s}{r}.
\]
For $\mathring{\vgu}_a$ in \eqref{vgu:normal},
we define a set $\PB_{m,\mathring{\vgu}_a}\subseteq (\dlp{0})^r$ consisting of the elements in $(\dlp{0})^r$ as follows:
\[
[\nabla^\mu \td,0,\ldots,0]^\tp, \quad
\forall\, \mu \in \dNN, |\mu|=m+1 \quad \mbox{and}\quad
\td e_j, \quad j=2,\ldots,r.
\]
It is trivial to check that
the integer shifts of elements in $\PB_{m,\mathring{\vgu}_a}$ span
the linear space $\PV_{m,\mathring{\vgu}_a}$.
Consequently, the integer shifts of elements in $\PB_{m,\vgu_a}:=\{U*u \setsp u\in \PB_{m,\mathring{\vgu}_a}\}$ span the linear space $\PV_{m,\vgu_a}$.
For $1\le p\le \infty$, as in \cite[(4.3)]{han03} and \cite[(7.2.2)]{hanbook}, we define a key smoothness quantity $\sm_p(a)$ by
\be \label{sma}
\sm_p(a):=\frac{d}{p}-\log_2 \rho_{m+1}(a,\vgu_a)_p
\ee
where
\be \label{rho}
\rho_{m+1}(a,\vgu_a)_p:=2^d\max\Big\{\limsup_{n\to \infty} \|a_n*u \|_{(\dlp{p})^r}^{1/n} \setsp u\in \PB_{m,\vgu_a}\Big\},
\ee
where $a_n:=2^{-dn} \sd_a^n (\td I_r)$ in \eqref{hsd:wn:F}, i.e., $\wh{a_n}(\xi):=\wh{a}(2^{n-1}\xi)\cdots\wh{a}(2\xi)\wh{a}(\xi)$.
In the definition of $\sm_p(a)$ we always take $m$ to be the largest possible integer for the mask $a$ having order $m+1$ sum rules with a matching filter $\vgu_a\in \lrs{0}{1}{r}$. A technical argument using joint spectral radius (\cite{hj98,han03,hanbook}) guarantees that $\sm_p(a)$ is independent of the choice of a matching filter $\vgu_a$ in \eqref{sma}.
The set $\PB_{m,\vgu_a}$ in \eqref{rho} can be also replaced by $\PV_{m,\vgu_a}$. For $m\in \NN$ and $1\le p\le \infty$, the Sobolev space $W^m_p(\dR)$ consists of all functions $f\in L_p(\dR)$ such that all the (weak) derivatives $f^{(\mu)}\in L_p(\dR)$ for all $\mu\in \dNN$ with $|\mu|\le m$. Let $\phi$ be a refinable vector function associated with the mask $a\in \dlrs{0}{r}{r}$ with the normalization condition $\wh{\vgu_a}(0)\wh{\phi}(0)=1$. Recall that the vector cascade operator $\cd_a$ is defined in \eqref{cd}.
We say that the vector cascade algorithm with mask $a\in \dlrs{0}{r}{r}$ is convergent in the Sobolev space $W_p^m(\dR)$ if the cascade sequence $\{\cd_a^n f\}_{n=1}^\infty$ converges to $\phi$ in $W^m_p(\dR)$ for all compactly supported admissible initial functions $f\in W^m_p(\dR)$ satisfying \eqref{initialf}. \cite[Theorem~4.3]{han03} shows that a vector cascade algorithm with a mask $a\in \dlrs{0}{r}{r}$ is convergent in $W^m_p(\dR)$ if and only if $\sm_p(a)>m$ (The same conclusion also holds if $W^m_p(\dR)$ with $p=\infty$ is replaced by $\dCH{m}$). Moreover,  $\sm_p(\phi)\ge \sm_p(a)$ always holds, see \cite{han03,hanbook} and references therein for more details.
As we already see in \cref{thm:converg},
the quantity $\sm_\infty(a)$ plays a key role in characterizing convergence and smoothness of generalized Hermite subdivision schemes as well. Moreover, $\sm_\infty(a)\ge \sm_2(a)-\frac{d}{2}$ always holds (see \cite[(4.7)]{han03}) and $\sm_2(a)$ can be computed by finding the eigenvalues of a finite matrix associated with the matrix mask $a\in \dlrs{0}{r}{r}$, e.g., see \cite[Theorem~1.1]{jj03}, \cite[Theorem~7.1]{han03} and \cite[Theorem~2.4]{han03smaa}.

\subsection{Admissible initial vector functions for vector cascade algorithms}
Note that \eqref{vgua:hsd} is satisfied for any sequence $\vgu_a\in \dlrs{0}{1}{r}$ satisfying
\be \label{special:vgua}
\wh{\vgu_a}(\xi)=\left[(i\xi)^{\nu_1}\wh{c_1}(\xi),\ldots,
(i\xi)^{\nu_r}\wh{c_r}(\xi)\right]+\bo(\|\xi|^{m+1}),\quad \xi\to 0
\ee
for some sequences $c_1,\ldots,c_r\in \dlp{0}$ with $\wh{c_1}(0)=\cdots=\wh{c_r}(0)=1$. For dimension $d=1$, it is also trivial to observe that \eqref{vgua:hsd} implies \eqref{special:vgua}. However, for dimension $d>1$, this is no longer true and \eqref{vgua:hsd} may hold while \eqref{special:vgua} fails.
Moreover, even in dimension one, elements in an ordered multiset $\ind$ in this paper are allowed to repeat.
Consequently, the technique in \cite{han20} for studying convergence of univariate Hermite subdivision schemes is no longer applicable to the multivariate setting. We need to circumvent this difficulty by developing other more general techniques.

To prove \cref{thm:converg}, we need the following auxiliary result on admissible initial vector functions for vector cascade algorithms,
which is of interest in its own right.

\begin{lemma}\label{lem:initialfunc}
Let $m\in \NN$ and $r\in \N$. Let $\ind=\{\nu_1,\ldots,\nu_r\}\subseteq \ind_m$ with $\nu_1=\mathbf{0}\in \dNN$ be an ordered multiset as in \eqref{indset}.
Let $\vgu_H\in \dlrs{0}{1}{r}$ satisfy
\be \label{vguH}
\wh{\vgu_H}(\xi)=\left[
(i\xi)^{\nu_1},(i\xi)^{\nu_2},\ldots,
(i\xi)^{\nu_r}\right]+\bo(\|\xi\|^{m+1}),\quad \xi \to 0.
\ee
Pick any subset $J_\ind\subseteq \{1,\ldots,r\}$ such that elements in $\{\nu_\ell \setsp \ell\in J_\ind\}$ are distinct and there is $\ell_0\in J_\ind$ such that $\nu_{\ell_0}=\mathbf{0}\in \dNN$.
Then there always exists a compactly supported spline vector function $h=[h_1,\ldots,h_r]^\tp\in (\dCH{m})^r$ such that
\be \label{initial:cond}
\wh{\vgu_H}(\xi)\wh{h}(\xi+2\pi k)=\td(k)+\bo(\|\xi\|^{m+1}),\quad \xi \to 0, \; \forall\; k\in \dZ
\ee
and
\be \label{intial:hermite}
h_\ell=0 \quad \forall\; \ell\in \{1,\ldots,r\}\bs J_\ind
\quad \mbox{and}\quad
h^{(\nu_\ell)}(k)=\td(k)e_\ell, \qquad \forall\; k\in \dZ, \ell\in J_\ind.
\ee
Moreover, let $\vgu_a\in \lrs{0}{1}{r}$ be any given sequence with $\wh{\vgu_a}(0)\ne 0$. For $r=1$, the function $f=h\in \dCH{m}$ is obviously an admissible initial vector function satisfying \eqref{initialf}. For $r>1$, \cref{lem:vgu} guarantees that there exists
a strongly invertible sequence $U\in \dlrs{0}{r}{r}$ satisfying
\be \label{U:vgua}
\wh{\vgu_a}(\xi)\wh{U}(\xi)=\wh{\vgu_H}(\xi)+\bo(\|\xi\|^{m+1}),\quad \xi \to 0
\ee
and then the following compactly supported spline vector function
\be \label{f:initial}
f:=U*h=\sum_{k\in \dZ} U(k) h(\cdot-k)
\ee
must be an admissible initial vector function in $\dCH{m}$ satisfying a stronger version than \eqref{initialf}:
\be \label{initialf:strong}
\wh{\vgu_a}(\xi)\wh{f}(\xi)=1+\bo(\|\xi\|^{m+1})
\quad \mbox{and}\quad
\wh{\vgu_a}(\xi)\wh{f}(\xi+2\pi k)=\bo(\|\xi\|^{m+1}),\quad \xi\to 0,\; \forall\; k\in \dZ\bs\{0\}.
\ee
\end{lemma}

\bp We first recall a refinable Hermite interpolant
(e.g., see \cite[Proposition~6.2.2]{hanbook}):
\[
\varphi_{\ell}(x):=
\begin{cases}
(1-x)^{m+1} \frac{x^{\ell}}{\ell!}\sum_{j=0}^{m-\ell}
\frac{(m+j)!}{m!j!} x^j, &\text{$x\in [0,1]$},\\
(1+x)^{m+1} \frac{x^{\ell}}{\ell!}\sum_{j=0}^{m-\ell}
\frac{(m+j)!}{m!j!} (-x)^j, &\text{$x\in [-1,0)$},\\
0, &\text{$x\in \R\bs [-1,1]$},
\end{cases}
\]
for $\ell=0,\ldots,m$. Then $\varphi\in (\CH{m})^{m+1}$ is a compactly supported vector function possessing the Hermite interpolation property:
\be \label{theta:hint}
\varphi_{\ell}^{(j)}(k)=\td(\ell-j)\td(k),\qquad \forall\; \ell,j=0,\ldots,m\quad \mbox{and}\quad k\in \Z,
\ee
and by \cite[Corollary~5.2]{han03} or \cite[Theorem~6.2.3]{hanbook}, we have
\be \label{theta:approx}
\left[1,i\xi,\ldots, (i\xi)^m \right] \wh{\varphi}(\xi+2\pi k)=\td(k)+\bo(|\xi|^{m+1}),\qquad \xi\to0, k\in \Z.
\ee
For $\mu=(\mu_1,\ldots,\mu_d)^\tp \in \dNN$ with $|\mu|_\infty:=\max(\mu_1,\ldots,\mu_d)\le m$, using the tensor product, we define
\[
\Phi_\mu(x_1,\ldots,x_d):=\varphi_{\mu_1}(x_1)\cdots
\varphi_{\mu_d}(x_d),\qquad (x_1,\ldots,x_d)^\tp \in \dR.
\]
Define $h_\ell:=0$ for all $\ell\in \{1,\ldots,r\}\bs J_\ind$,
$h_\ell:=\Phi_{\nu_\ell}$ for $\ell\in J_{\ind}\bs\{\ell_0\}$, and
\[
h_{\ell_0}:=\Phi_{\nu_{\ell_0}}+
\sum_{|\mu|_\infty\le m, \mu \not\in \{\nu_\ell \setsp \ell\in J_\ind\}}
c_\mu*\Phi_\mu,
\]
where $c_\mu\in \dlp{0}$ are chosen to satisfy $\wh{c_\mu}(\xi)=(i\xi)^\mu+\bo(\|\xi\|^{m+1})$ as $\xi \to 0$. Note that $\nu_{\ell_0}=\mathbf{0}\in \dNN$.
Using \eqref{theta:hint}, we can directly check that \eqref{intial:hermite} holds. Because $\Phi_\mu$ is constructed through tensor product, we trivially deduce from \eqref{theta:approx} that
\be \label{tensorhermite}
\sum_{\mu\in \dNN, |\mu|_\infty\le m}
(i\xi)^\mu \wh{\Phi_\mu}(\xi+2\pi k)=\td(k)+\bo(\|\xi\|^{m+1}),\quad \xi \to 0, k\in \dZ.
\ee
Now by the definition of the vector function $h$, for every $k\in \dZ$, we have
\begin{align*}
\wh{\vgu_H}(\xi)\wh{h}(\xi+2\pi k)
&=\wh{h_{\ell_0}}(\xi+2\pi k)
+\sum_{\ell\in J_\ind\bs \{\ell_0\}}
(i\xi)^{\nu_\ell} \wh{h_\ell}(\xi+2\pi k)+\bo(\|\xi\|^{m+1})\\
&=\sum_{\mu\in \dNN, |\mu|_\infty\le m}
(i\xi)^\mu \wh{\Phi_\mu}(\xi+2\pi k)+\bo(\|\xi\|^{m+1})
=\td(k)+\bo(\|\xi\|^{m+1}),\quad \xi \to0,
\end{align*}
where we used \eqref{tensorhermite}.
This proves \eqref{initial:cond}.

By definition of $f$ in \eqref{f:initial},
$f\in (\dCH{m})^r$ is a compactly supported vector function and $\wh{f}(\xi)=\wh{U}(\xi)\wh{h}(\xi)$.
For $r=1$, it is obvious that $f=h$ is an admissible initial function satisfying \eqref{initialf}.
For $r>1$, we deduce from \eqref{initial:cond} that
\[
\wh{\vgu_a}(\xi)\wh{f}(\xi+2\pi k)
=\wh{\vgu_a}(\xi)\wh{U}(\xi+2\pi k)\wh{h}(\xi+2\pi k)=
\wh{\vgu_H}(\xi)\wh{h}(\xi+2\pi k)+\bo(\|\xi\|^{m+1})
=\td(k)+\bo(\|\xi\|^{m+1})
\]
as $\xi \to 0$. That is, the vector function $f$ satisfies \eqref{initialf:strong} and hence is an admissible initial vector function satisfying a stronger version in \eqref{initialf:strong} than \eqref{initialf}.
\ep

\subsection{Proof of \cref{thm:converg}}
To prove \cref{thm:converg}, we make a remark about \cref{def:hsd}. Let $\phi$ be the underlying
basis vector function defined in \eqref{def:phi}.
It is easy to see that the limit function $\eta$ in \eqref{hsd:converg} in \cref{def:hsd} is given by $\eta=w_0*\phi=\sum_{k\in \dZ} w_0(k)\phi(\cdot-k)$. In other words, \eqref{hsd:converg} in \cref{def:hsd} can be equivalently rewritten as: for any constant $K>0$,
\[
\lim_{n\to \infty} \max_{k\in \dZ\cap [-2^nK, 2^n K]^d} |w_n(k) e_{\ell}- [w_0*\phi]^{(\nu_\ell)}(2^{-n}k)|=0,\qquad \forall\; \ell=1,\ldots,r.
\]
Therefore, a generalized Hermite subdivision scheme of type $\ind$ with mask $a$ is convergent with limit functions in $\dCH{m}$ as stated in \cref{def:hsd} if and only if
there exists a compactly supported vector function $\phi\in (\dCH{m})^r$ satisfying
\be \label{hsd:converg:phi}
\lim_{n\to \infty} \max_{k\in \dZ}
\|2^{(d+|\nu_\ell|)n}a_n(k)
e_{\ell}- \phi^{(\nu_\ell)}(2^{-n}k)\|=0,\quad \forall\; \ell=1,\ldots,r\quad \mbox{with}\quad a_n:=2^{-dn} \sd_a^n (\td I_r).
\ee
That is, $\wh{a_n}(\xi):=\wh{a}(2^{n-1}\xi)\cdots \wh{a}(2\xi)\wh{a}(\xi)$.
We are now ready to prove \cref{thm:converg}.

\bp[Proof of \cref{thm:converg}] We first prove the easy direction (ii)$\imply$(i) under the condition in \eqref{stable}.
Since \eqref{stable} implies \eqref{phi:cond0},
by \cref{thm:hsd:mask}, $\phi$ must be its basis vector function in $(\dCH{m})$ satisfying \eqref{hsd:converg:phi} and item (3) of \cref{thm:hsd:mask}.
Therefore, we deduce from item (3) of \cref{thm:hsd:mask} that $\phi$ trivially satisfies the condition in \eqref{initialf} with $f=\phi$ and hence, $\phi$ is an admissible initial vector function. Obviously, we have $\cd_a \phi=\phi$ by \eqref{refeq} and consequently, we trivially have $\lim_{n\to \infty} \|\cd_a^n \phi-\phi\|_{(\dCH{m})^r}=0$.
Because the integer shifts of $\phi$ are stable due to our assumption in \eqref{stable}, it follows from item~(2) of \cite[Theorem~4.3]{han03} that $\sm_\infty(a)>m$. This proves (ii)$\imply$(i).

We now prove the key step: (i)$\imply$(ii).
For the scalar case $r=1$, the connections between scalar subdivision schemes and scalar cascade algorithms are well studied, e.g., see \cite[Theorem~2.1]{hj06}.
In particular, it follows from
\cite[Theorem~2.1]{hj06} that (i)$\imply$(ii) for $r=1$. Thus, we assume $r>1$.

Arbitrarily take $\ell_*\in \{1,\ldots,r\}$ and then fix it. If $\nu_{\ell_*}=\mathbf{0}\in \dNN$, then we just take $J_\ind=\{\ell_*\}$; otherwise, we take $J_\ind=\{1,\ell_*\}$. Let a vector function $h$ as in \cref{lem:initialfunc} for $J_\ind$.
By \cref{lem:vgu}, there exists a strongly invertible sequence $U\in \dlrs{0}{r}{r}$ such that \eqref{U:vgua} holds.
Define a vector function $f$ as in \eqref{f:initial}. By \cref{lem:initialfunc}, the vector function $f$ satisfies \eqref{initialf:strong} and hence is an admissible initial vector function. Since item (i) holds, \eqref{cd:converg} holds.
Define $f_n:=\cd_a^n f$ and $b_n:=a_n*U$. By $f=U*h$,
\[
f_n=2^{dn}\sum_{k\in \dZ} a_n(k) f(2^n\cdot-k)=2^{dn} \sum_{k\in \dZ} b_n(k) h(2^n\cdot-k).
\]
The above identity can be also proved by noting
\[
\wh{f_n}(\xi)=\wh{a_n}(2^{-n}\xi) \wh{f}(2^{-n}\xi)= \wh{a_n}(2^{-n}\xi)\wh{U}(2^{-n}\xi)\wh{h}(2^{-n}\xi)
=\wh{b_n}(2^{-n}\xi)\wh{h}(2^{-n}\xi).
\]
Since $h\in (\dCH{m})^r$ satisfies \eqref{intial:hermite}, we deduce from $\ell_*\in J_\ind$ and \eqref{intial:hermite} that
\[
f_n^{(\nu_{\ell_*})}(2^{-n}k)
=2^{dn} \sum_{j\in \dZ} b_n(j) h^{(\nu_{\ell_*})}(k-j)
2^{|\nu_{\ell_*}|n}
=2^{(d+|\nu_{\ell_*}|)n} \sum_{j\in \dZ} b_n(j) \td(k-j)e_{\ell_*}
=2^{(d+|\nu_{\ell_*}|)n} b_n(k)e_{\ell_*}.
\]
By item (i) and $|\nu_{\ell_*}|\le m$, $\lim_{n\to \infty}\|f_n^{(\nu_{\ell_*})}-
\phi^{(\nu_{\ell_*})}\|_{(\dCH{})^r}=0$. Hence, by $b_n=a_n*U$ and the above identity,
\be \label{anU}
\lim_{n\to \infty} \sup_{k\in \dZ}\|
2^{(d+|\nu_{\ell_*}|)n} (a_n*U e_{\ell_*})(k)-\phi^{(\nu_{\ell_*})}(2^{-n}k)\|=0.
\ee
Define a sequence $u:=\td e_{\ell_*}-U e_{\ell_*}$.
It follows directly from \eqref{U:vgua} and \eqref{vgua:hsd} that
\[
\wh{\vgu_a}(\xi) \wh{u}(\xi)=
\wh{\vgu_a}(\xi) e_{\ell_*}-
\wh{\vgu_a}(\xi)\wh{U}(\xi) e_{\ell_*}
=\left(\wh{\vgu_a}(\xi)-\wh{\vgu_H}(\xi)\right) e_{\ell_*}+\bo(\|\xi\|^{m+1})
=\bo(\|\xi\|^{|\nu_{\ell_*}|+1}),\quad \xi \to 0.
\]
This proves $u\in \PV_{|\nu_{\ell_*}|, \vgu_a}$.
Since $|\nu_{\ell_*}|\le m$,
it is known in \cite[Theorem~5.7.6]{hanbook} for dimension one and \cite[Theorem~3.1]{han03smaa} for high dimensions that
\[
\rho_{|\nu_{\ell_*}|+1}(a,\vgu_a)_\infty=
\max(2^{-|\nu_{\ell_*}|-1}, \rho_{m+1}(a,\vgu_a)_\infty)=
\max(2^{-|\nu_{\ell_*}|-1}, 2^{-\sm_\infty(a)}),
\]
because  $\rho_{m+1}(a,\vgu_a)_\infty=
2^{-\sm_\infty(a)}$ and $|\nu_{\ell_*}|\le m$.
Since $u\in \PV_{|\nu_{\ell_*}|, \vgu_a}$ and the integer shifts of $\PB_{|\nu_{\ell_*}|, \vgu_a}$ span the linear space $\PV_{|\nu_{\ell_*}|, \vgu_a}$, we conclude from the definition of $\rho_{|\nu_{\ell_*}|+1}(a,\vgu_a)_\infty$ in \eqref{rho}
that
\[
\limsup_{n\to \infty} 2^{d+|\nu_{\ell_*}|}\|a_n*u\|_{(\dlp{p})^r}^{1/n}
\le 2^{|\nu_{\ell_*}|} \rho_{|\nu_{\ell_*}|+1}(a,\vgu_a)_\infty
=\max\left(2^{-1},  2^{-(\sm_\infty(a)-|\nu_{\ell_*}|)}\right)<1,
\]
since $\sm_\infty(a)>m\ge |\nu_{\ell_*}|$,
where the above inequality also follows from \cite[item (5) of Theorem~4.3]{han03}.
Therefore, we conclude from the above inequality that
\be \label{anuell}
\lim_{n\to \infty} 2^{(d+|\nu_{\ell_*}|)n} \|a_n*u\|_{(\dlp{\infty})^r}=0.
\ee
Since $Ue_{\ell_*}=\td e_{\ell_*}-u$, we conclude from \eqref{anU} and \eqref{anuell} that
$\lim_{n\to \infty} \sup_{k\in \dZ}\|
2^{(d+|\nu_{\ell_*}|)n} a_n(k) e_{\ell_*}-\phi^{(\nu_{\ell_*})}(2^{-n}k)\|=0$. This proves \eqref{hsd:converg:phi} for $\ell=\ell_*$. Because $\ell_*\in \{1,\ldots,r\}$ can be chosen arbitrarily, this completes the proof of \eqref{hsd:converg:phi}. Therefore, we proved (i)$\imply$(ii).
\ep

\section{Generalized Hermite Subdivision Schemes with Interpolation Properties}
\label{sec:hsdint}

In this section, we shall prove \cref{thm:int}
for the convergence of interpolatory generalized Hermite subdivision schemes, which consist of a special family of multivariate generalized Hermite subdivision schemes. Motivated by \cref{thm:int},
we shall introduce the notion of linear-phase moments for a generalized Hermite subdivision scheme to have the polynomial-interpolation property.

%We shall see that such interpolatory generalized Hermite subdivision schemes include all different types of interpolatory subdivision schemes known in the literature as special cases and consist of a much wider class of interpolatory subdivision schemes.
%For interpolatory generalized Hermite subdivision schemes, we have some natural compatibility condition on the ordered multiset $T$ with $\ind$.
%We assume that there is a mapping $\theta$ satisfying

To study interpolatory generalized Hermite subdivision schemes, we need the following result.

\begin{lemma}\label{lem:pv}
Let $M\in \NN$, $\nu\in \ind_M$ and $\tau\in \dR$. For $v\in \dlp{0}$, $\pp*v=\pp^{(\nu)}(\cdot+\tau)$ for all $\pp\in \PL_M$ if and only if
$\wh{v}(\xi)=(i\xi)^\nu e^{i\tau \cdot \xi}+\bo(\|\xi\|^{M+1})$ as $\xi \to 0$.
\end{lemma}

\bp For $y\in \dR$, we take $v_y\in \dlp{0}$ satisfying $\wh{v_{y}}(\xi)=e^{i y\cdot\xi}+\bo(\|\xi\|^{M+1})$ as $\xi \to 0$.
For any $\pp\in \PL_M$, we deduce from \eqref{conv:poly} that
\[
\pp*v_{y}=\sum_{\mu\in \dNN}\frac{(-i)^{|\mu|}}{\mu!} \pp^{(\mu)}(\cdot) (iy)^\mu=
\sum_{\mu\in \dNN}\frac{1}{\mu!} \pp^{(\mu)}(\cdot) y^\mu
=\pp(\cdot+y),
\]
where we used the Taylor expansion of $\pp(\cdot+y)$ in the last identity.
Consequently, we have
\be \label{poly:vy}
\pp*v_{y}=\pp(\cdot+y),\quad \forall\, \pp\in \PL_M, y\in \dR\quad \mbox{with}\quad
\wh{v_{y}}(\xi)=e^{iy \cdot \xi}+\bo(\|\xi\|^{M+1}),\quad \xi \to 0.
\ee
By \eqref{poly:vy}, we have $\pp=\pp(\cdot+\tau-\tau)=(\pp(\cdot+\tau))*v_{-\tau}$ and hence,
$\pp*v=((\pp(\cdot+\tau))*v_{-\tau})*v=
(\pp(\cdot+\tau))*(v_{-\tau}*v)$.
Now by \eqref{conv:poly}, we have
\be \label{pv}
\pp*v=
\sum_{\alpha \in \ind_M}
\frac{(-i)^{|\alpha|}}{\alpha!}
\pp^{(\alpha)}(\cdot+\tau) [e^{-i\tau \cdot \xi}\wh{v}(\xi)]^{(\alpha)}(0),\qquad \pp\in \PL_M.
\ee
Suppose that $\wh{v}(\xi)=(i\xi)^\nu e^{i\tau \cdot \xi}+\bo(\|\xi\|^{M+1})$ as $\xi \to 0$. Then $e^{-i\tau \cdot \xi} \wh{v}(\xi)=(i\xi)^\nu +\bo(\|\xi\|^{M+1})$ as $\xi \to 0$. It follows from \eqref{pv} that
$\pp*v=
\sum_{\alpha \in \ind_M}
\frac{(-i)^{|\alpha|}}{\alpha!}
\pp^{(\alpha)}(\cdot+\tau)
i^{|\nu|} \nu! \td(\alpha-\nu)
=\pp^{(\nu)}(\cdot+\tau)$ for all $\pp\in \PL_M$.

Conversely, suppose that $\pp*v=\pp^{(\nu)}(\cdot+\tau)$ for all $\pp\in \PL_M$. By \eqref{pv}, we have
\[
\pp^{(\nu)}(\cdot+\tau)=\pp*v=
\sum_{\alpha \in \ind_M}
\frac{(-i)^{|\alpha|}}{\alpha!}
\pp^{(\alpha)}(\cdot+\tau) [e^{-i\tau \cdot \xi}\wh{v}(\xi)]^{(\alpha)}(0),
\]
from which we must have $[e^{-i\tau \cdot \xi}\wh{v}(\xi)]^{(\alpha)}=i^{|\nu|}\nu! \td(\alpha-\nu)$ for all $\alpha\in \ind_M$. Therefore, we conclude that
$\wh{v}(\xi)=(i\xi)^\nu e^{i\tau \cdot \xi}+\bo(\|\xi\|^{M+1})$ as $\xi \to 0$.
\ep

We are now ready to prove \cref{thm:int}  for convergence of interpolatory generalized Hermite subdivision schemes of type $\ind$ and translation $T$.

\bp[Proof of \cref{thm:int}] To prove item (1), by the recursive formula in \eqref{wn}, we have
\[
w_n=\sd_{\sD_\ind^{n-1} a\sD_\ind^{-n}} w_{n-1}=
2^d\sum_{j\in \dZ} w_{n-1}(j) \sD_\ind^{n-1} a(\cdot-2j)\sD_\ind^{-n}.
\]
Define $\beta_\ell:=2\tau_\ell-\tau_{\theta(\ell)}\in \dZ$ by the condition in \eqref{theta}.
For $k\in \dZ$ and $\ell=1,\ldots,r$, we have
\be \label{wn:2level}
\begin{split}
w_n(2k+\beta_\ell)
&=2^d\sum_{j\in \dZ} w_{n-1}(j) \sD_\ind^{n-1} a(2k+\beta_\ell-2j)\sD_\ind^{-n} \\
&=2^d\sum_{j\in \dZ} w_{n-1}(k-j) \sD_\ind^{n-1} a(2j+\beta_\ell )\sD_\ind^{-n}.
\end{split}
\ee
Now we deduce from the interpolation condition in \eqref{interpolation} that
for all $k\in \dZ$, $n\in \N$ and $\ell=1,\ldots,r$,
\[
2^d\sum_{j\in \dZ} w_{n-1}(k-j) \sD_\ind^{n-1} a(2j+\beta_\ell )\sD_\ind^{-n}e_{\theta(\ell)}=
w_n(2k+\beta_\ell)e_{\theta(\ell)}=
w_{n-1}(k)e_{\ell}.
\]
Since $w_0\in \dlrs{}{1}{r}$ is arbitrary, we conclude from the above identity with $n=1$ that
\be \label{int:cond}
2^d  a(\beta_\ell)2^{|\nu_\ell|} e_{\theta(\ell)}
=e_{\ell}
\quad \mbox{and}\quad
a(2k+\beta_\ell)2^{|\nu_\ell|} e_{\theta(\ell)}=0,\qquad \forall\; k\in \dZ\bs\{0\}, \ell=1,\ldots,r,
\ee
where we used $\nu_{\theta(\ell)}=\nu_\ell$ and
$\sD_\ind^{-1} e_{\theta(\ell)}=2^{|\nu_{\theta(\ell)}|} e_{\theta(\ell)}=
2^{|\nu_\ell|}e_{\theta(\ell)}$ by the definition of $\sD_\ind$ in \eqref{hsd:wn}.
Therefore, \eqref{int:cond} obviously implies \eqref{int:mask}.
This proves item (1).

We now prove item (2).
Let $w_0=\td I_r$ be the initial data and $\{w_n\}_{n=1}^\infty$ be the refinement data defined in \eqref{hsd:wn}.
Now by the interpolation condition in \eqref{interpolation} and using induction, we have
\begin{align*}
&w_n\left(2^n k+2^{n-1}\beta_\ell+\cdots +2\beta_{\theta^{n-2}(\ell)}+\beta_{\theta^{n-1}(\ell)}\right)e_{\theta^n(\ell)}\\
&\qquad =w_{n-1}\left(2^{n-1} k+2^{n-2}\beta_\ell+\cdots +2\beta_{\theta^{n-3}(\ell)}+\beta_{\theta^{n-2}(\ell)}\right)
e_{\theta^{n-1}(\ell)}
%\\&\qquad
=w_1(2k+\beta_\ell)e_{\theta(\ell)}=w_0(k) e_{\ell}.
\end{align*}
Since $w_0(k)=\td(k) I_r$, from the above identity we must have
\be \label{wn:int}
w_n\left(2^n k+2^{n-1}\beta_\ell+\cdots +2\beta_{\theta^{n-2}(\ell)}+
\beta_{\theta^{n-1}(\ell)}\right)
e_{\theta^n(\ell)}
=\td(k)e_{\ell},\quad k\in \dZ, \ell=1,\ldots,r.
\ee
By \eqref{theta}, we have $\beta_{\theta^j(\ell)}=2\tau_{\theta^j(\ell)}-\tau_{
\theta^{j+1}(\ell)}$ for all $\ell=1,\ldots,r$ and $j\in \N$.
Hence, we have
\begin{align*}
2^{n-1}\beta_\ell+2^{n-2}\beta_{\theta(\ell)}+
\cdots+\beta_{\theta^{n-1}(\ell)}
&=\sum_{j=1}^n 2^{n-j} \beta_{\theta^{j-1}(\ell)}=
\sum_{j=1}^n 2^{n-j}[ 2\tau_{\theta^{j-1}(\ell)}-\tau_{\theta^j(\ell)}]\\
&=\sum_{j=1}^n [2^{n-(j-1)} \tau_{\theta^{j-1}(\ell)}-2^{n-j}\tau_{\theta^j(\ell)}]
=2^n \tau_\ell-\tau_{\theta^n(\ell)}.
\end{align*}
Hence, \eqref{wn:int} can be rewritten as
\be \label{wn:int:2}
w_n\left(2^n k+2^n\tau_\ell-\tau_{\theta^n(\ell)}\right)
e_{\theta^n(\ell)}
=\td(k)e_{\ell},\quad k\in \dZ, n\in \N, \ell=1,\ldots,r.
\ee
Note that $\nu_{\theta^n(\ell)}=\nu_\ell$ by \eqref{theta}.
By \cref{def:hsd} with the initial data $w_0=\td(k)I_r$, we must have $\eta=w_0*\phi=\phi$ and \eqref{hsd:converg} is equivalent to
\be \label{hsd:converg:phi:0}
\lim_{n\to \infty} \max_{k\in \dZ}
\| w_n(k) e_{\ell}-\phi^{(\nu_\ell)}(2^{-n}k)\|=0,\qquad \forall\; \ell=1,\ldots,r.
\ee
Now by \eqref{wn:int:2} and \eqref{hsd:converg:phi:0} and noting that
$\nu_{\theta^n(\ell)}=\nu_\ell$,
we have
\begin{align*}
\td(k) e_{\ell}
&=\lim_{n\to \infty} w_n(2^n k+2^n \tau_\ell-\tau_{\theta^n(\ell)}) e_{\theta^n(\ell)}
=
\lim_{n\to \infty} \phi^{(\nu_{\theta^n(\ell)})}
\left(2^{-n}(2^n k+2^n \tau_\ell-\tau_{\theta^n(\ell)})\right)\\
&=\lim_{n\to \infty}
\phi^{(\nu_\ell)}(k+ \tau_\ell-2^{-n}\tau_{\theta^n(\ell)})
=\phi^{(\nu_\ell)}(k+\tau_\ell),
\end{align*}
where we used the fact that $\phi^{(\nu_\ell)}$ is uniformly continuous on $\dR$, since $\phi^{(\nu_\ell)}$ is a compactly supported continuous function.
This proves \eqref{int:phi} and hence item (2).

Next, we prove item (3).
Since all the conditions in \cref{thm:hsd:mask} are satisfied, all the claims in items (1)--(4) of \cref{thm:hsd:mask} hold. In particular, by item (3) of \cref{thm:hsd:mask},
the mask $a$ must have order $m+1$ sum rules in \eqref{sr} with a matching filter $\vgu_a\in \dlrs{0}{1}{r}$ satisfying \eqref{vgua:hsd}.
Using coset sequences $\wh{a^{[\gamma]}}(\xi)=\sum_{k\in \R^d} a(\gamma+2k)e^{-ik\cdot \xi}$,
one can easily deduce that \eqref{sr} is equivalent to
\be \label{sr:coset}
e^{-i\gamma\cdot \xi} \wh{\vgu_a}(2\xi) \wh{a^{[\gamma]}}(2\xi)=2^{-d}
\wh{\vgu_a}(\xi)+\bo(\|\xi\|^{m+1}),\quad \xi \to 0, \; \forall\; \gamma \in \Gamma.
\ee
By the definition of the integer $M$ in item (3), we must have $M\ge m$ and in terms of coset sequences, \eqref{sr:coset} must hold with $m$ being replaced by $M$. In particular, we have
\be \label{sr:coset:nu}
e^{-i\beta_\ell \cdot \xi} \wh{\vgu_a}(2\xi)\wh{a^{[\beta_\ell]}}(2\xi)
=2^{-d}\wh{\vgu_a}(\xi)+\bo(\|\xi\|^{M+1}),\qquad \xi\to 0, \ell=1,\ldots,r.
\ee
By the proved item (1), the mask $a$ is interpolatory of type $\ind$ and translation $T$ and satisfies \eqref{int:mask}. Hence, we must have
$\wh{a^{[\beta_\ell]}}(2\xi)e_{\theta(\ell)}=
a(\beta_\ell )e_{\theta(\ell)}=2^{-d-|\nu_\ell|} e_{\ell}$. Now we deduce from \eqref{sr:coset:nu} that
\[
e^{-i\beta_\ell \cdot \xi}
2^{-d-|\nu_\ell|}
\wh{\vgu_a}(2\xi)e_{\ell}=
e^{-i\beta_\ell \cdot \xi} \wh{\vgu_a}(2\xi)\wh{a^{[\beta_\ell]}}(2\xi)
e_{\theta(\ell)}=
2^{-d}[\wh{\vgu_a}(\xi)]_{\theta(\ell)}+\bo(\|\xi\|^{M+1}),
\]
as $\xi \to 0$.
Because $\beta_\ell=2\tau_\ell-\tau_{\theta(\ell)}$ by \eqref{theta},
we deduce from the above identity that
\[
e^{-i\tau_{\theta(\ell)}\cdot \xi} [\wh{\vgu_a}(\xi)]_{\theta(\ell)}=
e^{-i \tau_\ell \cdot 2 \xi}
2^{-|\nu_\ell|}
[\wh{\vgu_a}(2\xi)]_{\ell}
+\bo(\|\xi\|^{M+1}),
\quad \xi\to 0,\ell=1,\ldots,r.
\]
Since $\nu_{\theta^n(\ell)}=\nu_\ell$ for all $n\in \N$, it follows directly from the above identity that
\be \label{vgua:theta}
e^{-i\tau_{\theta^n(\ell)}\cdot \xi} [\wh{\vgu_a}(\xi)]_{\theta^n(\ell)}=
2^{-n |\nu_\ell|} e^{-i \tau_\ell \cdot 2^n \xi}[\wh{\vgu_a}(2^n\xi)]_{\ell}
+\bo(\|\xi\|^{M+1}),
\quad \xi\to 0, n\in \N, \ell=1,\ldots,r.
\ee
For every fixed $\ell=1,\ldots,r$,
since $\theta$ is a mapping on the finite set $\{1,\ldots,r\}$, we must have $\theta^{n_1}(\ell)=\theta^{n_2} (\ell)$ for some integers $n_1<n_2$.
Consequently, it follows from \eqref{vgua:theta} that
\begin{align*}
2^{-n_1 |\nu_\ell|}
e^{-i \tau_{\ell} \cdot 2^{n_1} \xi}
[\wh{\vgu_a}(2^{n_1}\xi)]_\ell
&=e^{-i\tau_{\theta^{n_1}(\ell)}\cdot\xi}
[\wh{\vgu_a}(\xi)]_{\theta^{n_1}(\ell)}
+\bo(\|\xi\|^{M+1})
=e^{-i\tau_{\theta^{n_2}(\ell)}\cdot\xi}
[\wh{\vgu_a}(\xi)]_{\theta^{n_2}(\ell)}+\bo(\|\xi\|^{M+1})\\
&=2^{-n_2 |\nu_\ell|}
e^{-i \tau_{\ell} \cdot 2^{n_2} \xi}
[\wh{\vgu_a}(2^{n_2}\xi)]_\ell+\bo(\|\xi\|^{M+1}),\quad \xi \to 0.
\end{align*}
That is, for $n:=n_2-n_1>1$, we conclude from the above identity that
\be \label{vgua:theta:n}
e^{-i\tau_\ell\cdot \xi} [\wh{\vgu_a}(\xi)]_\ell=
2^{-n |\nu_\ell|} e^{-i\tau_\ell\cdot 2^n \xi} [\wh{\vgu_a}(2^n\xi)]_\ell+\bo(\|\xi\|^{M+1}),\quad \xi \to 0.
\ee
Since $\wh{\vgu_a}$ must satisfy \eqref{vgua:hsd} by item (3) of \cref{thm:hsd:mask}, we deduce from \eqref{vgua:hsd} that $e^{-i\tau_\ell\cdot \xi} [\wh{\vgu_a}(\xi)]_{\ell} =(i\xi)^{\nu_\ell}
+\bo(\|\xi\|^{|\nu_\ell|+1})$ as $\xi \to 0$.
Now using the Taylor expansion for $e^{-i\tau_\ell\cdot \xi} [\wh{\vgu_a}(\xi)]_{\ell}$ at $\xi=0$, we conclude from \eqref{vgua:theta:n} that $e^{-i \tau_\ell \cdot \xi}
[\wh{\vgu_a}(\xi)]_{\ell}
=(i\xi)^{\nu_\ell} +\bo(\|\xi\|^{M+1})$ as $\xi\to 0$.
This proves \eqref{int:vgua} and item (3).

Finally, we prove that item (3) implies item (4).
By item (3), we conclude from \cref{thm:sd:sr} that $\sd_a \vec{\pp}_\mu=2^{-|\mu|} \vec{\pp}_\mu$.
By \eqref{int:vgua}, \eqref{int:pmu} is a direct consequence of \cref{lem:pv}.
The identity $\sum_{k\in \dZ} \pp_\mu(k)\phi(\cdot-k)=\sum_{k\in \Z} (\frac{(\cdot)^\mu}{\mu!}*\vgu_a)(k)\phi(\cdot-k)=
\frac{(\cdot)^\mu}{\mu!}$ for all $\mu\in \ind_M$ follows directly from the same proof as in item (2) of \cref{thm:hsd:mask} (also see \cite{han03} for details).
This proves item (4).

Conversely, suppose that item (2) holds, i.e., \eqref{int:phi} holds. The generalized Hermite interpolation property in \eqref{int:phi} obviously implies that the integer shifts of $\phi$ must be linearly independent
(see Subsection~\cref{lem:vgu} for details).
It is known (e.g., see \cite[Theorem~5.2.1]{hanbook} and references therein) that
the integer shifts of $\phi$ are linearly independent if and only if $\mspan\{\wh{\phi}(\xi+2\pi k) \setsp k\in \dZ\}=\C^r$ for all $\xi\in \C^d$, which automatically implies the conditions in \eqref{stable} and \eqref{phi:cond0}. This proves item (ii).
By \cite[Corollary~5.1]{han03}, we must have $\sm_\infty(a)>m$ and hence by \cref{thm:converg}, the generalized Hermite subdivision scheme of type $\ind$ with mask $a$ is convergent with limit functions in $\dCH{m}$. This proves item (i).
Using \eqref{int:phi} and \eqref{theta}, we deduce from the refinement equation  \eqref{refeq:phi} that for all $k\in \dZ$ and $\ell=1,\ldots,r$,
\[
\td(k) e_\ell=\phi^{(\nu_\ell)}(k+\tau_\ell)=
2^{|\nu_\ell|+d}\sum_{j\in \dZ}
a(j) \phi^{(\nu_{\theta(\ell)})}((2k-j+\beta_\ell)+\tau_{\theta(\ell)})
=2^{|\nu_\ell|+d} a(2k+\beta_\ell) e_{\theta(\ell)}.
\]
This proves \eqref{int:mask}, i.e., item (1) holds. By the definition of $w_n$ in \eqref{wn}, we deduce from \eqref{int:mask} that \eqref{interpolation} holds.
Indeed, \eqref{wn} implies \eqref{wn:2level}, from which and $\nu_{\theta(\ell)}=\nu_\ell$ we have
\begin{align*}
w_n(2k+\beta_\ell)
&=2^d \sum_{j\in \dZ} w_{n-1}(k-j) \sD_\ind^{n-1} a(2j+\beta_\ell) e_{\theta(\ell)} 2^{|\nu_{\theta(\ell)}|n} \\
&=2^{|\nu_{\theta(\ell)}|n-|\nu_\ell|}
\sum_{j\in \dZ} w_{n-1}(k-j) \td(j) \sD_\ind^{n-1} e_\ell\\
&=w_{n-1}(k)e_\ell 2^{|\nu_{\theta(\ell)}|n-|\nu_\ell|-|\nu_\ell|(n-1)}=
w_{n-1}(k) e_\ell.
\end{align*}
Hence, item (ii) holds. Items (3) and (4) are direct consequences of \cref{thm:hsd:mask}.
\ep

We now explain the role played by the mapping $\theta$ in \eqref{theta}.
We observe from the refinement equation in \eqref{refeq:phi} that
$\phi^{(\nu_\ell)}=2^{|\nu_\ell|+d}\sum_{j\in \dZ} a(j) \phi^{(\nu_\ell)}(2\cdot-j)$.
If $\nu_{\theta(\ell)}=\nu_\ell$ for some $\ell=1,\ldots,r$, then
for all $k\in \dZ$ we have
\[
\phi^{(\nu_\ell)}(k+\tau_\ell)
=2^{|\nu_\ell|+d}\sum_{j\in \dZ} a(j)\phi^{(\nu_{\theta(\ell)})}
((2k-j)+\tau_{\theta(\ell)}+
(2\tau_\ell-\tau_{\theta(\ell)})).
\]
Hence, it is natural to require $2\tau_\ell-\tau_{\theta(\ell)}\in \dZ$ in \eqref{theta} whenever $\nu_{\theta(\ell)}=\nu_\ell$. For any integers $s_1,\ldots, s_r\in \dZ$, we observe that $\{\tau_1-s_1,\ldots,\tau_r-s_r\}$ is the translation multiset for the new basis vector function $[\phi_1(\cdot+s_1),\ldots,\phi_r(\cdot+s_r)]^\tp$, which is essentially the same as the original basis vector function $\phi$ through integer shifts. Hence, we often translate elements in $T$ by integers and only need to consider $T\subseteq [0,1)^d$.
%After integer shifting elements in $T$ to $[0,1)^d$, the mapping $\theta$ in \eqref{theta} is often unique.
If elements in $\ind$ do not repeat (e.g., $\ind=\{\mathbf{0}\}$ for scalar subdivision schemes, or $\ind=\ind_m$ for standard Hermite subdivision schemes), to satisfy \eqref{theta}, then obviously $\theta$ must be the identity mapping and hence $T\subseteq \dZ$ (which becomes essentially $T=\{0,\ldots,0\}$).
If $\nu_\ell=\nu_{\ell'}$  with $\ell\ne \ell'$ in $\ind$, then one can easily deduce that the definition of a generalized Hermite interpolant in \eqref{int:phi} forces $\tau_\ell-\tau_{\ell'}\not\in \dZ$.
Let $\ind\subseteq \ind_m$ be an ordered multiset with some elements repeated. For each $\nu\in \ind_m$, we define $J_\nu:=\{\ell\in \{1,\ldots,r\} \setsp \nu_\ell=\nu\}$, i.e., all the indices of the repeated element $\nu$ in $\ind$. Then $\{1,\ldots,r\}$ is a disjoint union of $J_\nu, \nu\in \ind_m$ and $\theta$ maps $J_\nu$ into $J_\nu$.
For each nonempty set $J_\nu$, we can pick up a $d\times d$ integer matrix $N_\nu$ satisfying $\det (N_\nu)=\#J_\nu$. Now we can choose a multiset
$T$ with $\{\tau_\ell \setsp \ell\in J_\nu\}=\Omega_{N_\nu}:=[N_\nu^{-1}\dZ]\cap [0,1)^d$. For every  $\omega\in \Omega_{N_\nu}$, there exists a unique integer $\beta_\omega\in \dZ$ such that $2\omega-\beta_\omega \in \Omega_{N_\mu}$. This induces a natural mapping $\theta$ in \eqref{theta}.
If $\ind$ consists of $N$ copies of $\ind_m$ (i.e., every $\mu\in \ind_m$ has multiplicity $N$ in $\ind$), then the generalized Hermite interpolants $\phi$ in \cref{thm:int} become the interpolating refinable vector functions in \cite{hkz09,hz09}.
Hence, \cref{thm:int} not only covers all interpolatory Hermite subdivision schemes known in the literature but also generalizes them to a much wider class of interpolatory generalized Hermite subdivision schemes including Birkhoff interpolation sets $\ind$ and more general interpolation multisets $T$.

For a convergent generalized Hermite subdivision scheme, the interpolation property in \cref{def:int} imposes a stringent condition on its mask in \eqref{int:mask}. In many applications, such interpolation property can be weakened and one may only require interpolation property for polynomials.
This leads to the notion of linear-phase moments which was first explicitly introduced in \cite{han10} for scalar complex orthogonal wavelets.
Such polynomial-interpolation property is of interest for nearly shape preservation subdivision schemes in CAGD.
%(e.g., see \cite{fhs22,hanbook}).
Dual (also called face-based) Hermite subdivision schemes were discussed in \cite{hy06} by attaching the data $w_n(k) e_\ell$ in \eqref{hsd:converg} at the position $2^{-n}(k+\tau_\ell)$ with $\tau_\ell\in \dR$ such that not all $\tau_\ell$ are zero.
For an ordered multiset $T=\{\tau_1,\ldots,\tau_r\}\subseteq \dR$, we say that \emph{a generalized Hermite subdivision scheme of type $\ind$ and translation $T$ with mask $a$ is convergent with limit functions in $\dCH{m}$} if \cref{def:hsd} holds with \eqref{hsd:converg} being replaced by
\be \label{hsd:converg:face}
\lim_{n\to \infty} \max_{k\in \dZ\cap [-2^nK, 2^n K]^d} |w_n(k) e_{\ell}- \eta^{(\nu_\ell)}(2^{-n}(k+\tau_\ell))|=0,\qquad \forall\; \ell=1,\ldots,r.
\ee
Obviously, \eqref{hsd:converg} is a special case of \eqref{hsd:converg:face} with $\tau_1=\cdots=\tau_r=0$. Moreover, due to the uniform continuity of $\eta^{(\nu_\ell)}$ on compact sets, it is easy to observe that \eqref{hsd:converg} holds if and only if \eqref{hsd:converg:face} holds. Consequently, the face-based or dual version does not affect the convergence property in \cref{def:hsd} at all. However, it will make a difference for the polynomial-interpolation property below.

\begin{theorem}\label{thm:lpm}
Let $m\in \NN$ and $r\in \N$.
Take ordered multisets $\ind=\{\nu_1,\ldots,\nu_r\}\subseteq \ind_m$ with $\nu_1=\mathbf{0}$ as in \eqref{indset} and $T=\{\tau_1,\ldots,\tau_r\}\subseteq \dR$ for translation.
Suppose that the generalized Hermite subdivision scheme of type $\ind$ with mask $a\in \dlrs{0}{r}{r}$ is
convergent with limit functions in $\dCH{m}$. Let $\phi$ be its basis vector function.
For any integer $M\ge m$, the following are equivalent to each other:
\begin{enumerate}
\item[(1)]
%The matrix mask $a$ has \emph{order $M+1$ linear-phase moments of type $\ind$ and translation $T$}, i.e.,
    The mask $a$ has order $M+1$ sum rules with a matching filter $\vgu_a\in \dlrs{0}{1}{r}$ satisfying \eqref{int:vgua}.
\item[(2)] The polynomial-interpolation property of order $M+1$ holds: For all $\pp\in \PL_M$ and $n\in \NN$,
\be \label{polyinterpolation}
\sd_a^n [\pp^{(\nu_1)}(\cdot+\tau_1),\ldots,
\pp^{(\nu_r)}(\cdot+\tau_r)]
=[[\pp(2^{-n}\cdot)]^{(\nu_1)}(\cdot+\tau_1),\ldots,
[\pp(2^{-n}\cdot)]^{(\nu_r)}(\cdot+\tau_r)].
\ee
\item[(3)] For every initial sequence $w_0:=[\pp^{(\nu_1)}(\cdot+\tau_1),\ldots,
\pp^{(\nu_r)}(\cdot+\tau_r)]$,
\[
w_n:=(\sd_a^n w_0)\sD_\ind^{-n}
=[\pp^{(\nu_1)}(2^{-n}(\cdot+\tau_1)),\ldots,
\pp^{(\nu_r)}(2^{-n}(\cdot+\tau_r))],\qquad
\forall\; \pp\in \PL_M, n\in \NN.
\]
\end{enumerate}
Moreover, any of the above items (1)--(3) implies
\begin{enumerate}
\item[(i)] The limit in \eqref{hsd:converg:face} trivially holds, that is, $w_n(k)e_\ell=\eta^{(\nu_\ell)}(2^{-n}(k+\tau_\ell))$ with the initial sequence $w_0:=[\pp^{(\nu_1)}(\cdot+\tau_1),\ldots,
\pp^{(\nu_r)}(\cdot+\tau_r)]$ for all $\ell=1,\ldots,r$, $k\in \dZ$, and $\pp\in \PL_M$, where $\eta:=w_0*\phi$.
\item[(ii)] For any mapping  $\theta: \{1,\ldots,r\}\rightarrow \{1,\ldots,r\}$ satisfying $\nu_{\theta(\ell)}=\nu_\ell$ for all $\ell=1,\ldots,r$,
the identities in \eqref{interpolation} hold with $w_0:=[\pp^{(\nu_1)}(\cdot+\tau_1),\ldots,
\pp^{(\nu_r)}(\cdot+\tau_r)]$ for all $\pp\in \PL_M$.
\end{enumerate}
\end{theorem}

\bp (1)$\imply$(2).  For $\pp\in \PL_M$,
we conclude from \cref{lem:pv} and \eqref{int:vgua} that \eqref{int:pmu} holds and hence $w_0:=[\pp^{(\nu_1)}(\cdot+\tau_1),\ldots,
\pp^{(\nu_r)}(\cdot+\tau_r)]=\pp*\vgu_a$.
By \cref{thm:sd:sr} and \eqref{sd:poly:vva}, we have
\[
\sd_a^n (\pp*\vgu_a)=(\pp(2^{-n}\cdot))*\vgu_a
=[[\pp(2^{-n}\cdot)]^{(\nu_1)}(\cdot+\tau_1),\ldots,
[\pp(2^{-n}\cdot)]^{(\nu_r)}(\cdot+\tau_r)],
\]
where we used \cref{lem:pv} and \eqref{int:vgua} again in the last identity.
This proves \eqref{polyinterpolation} and item (2).

(2)$\imply$(1). Let $\vgu_a\in \dlrs{0}{1}{r}$ satisfy \eqref{int:vgua} as in item (1). By \cref{lem:pv} and \eqref{int:vgua},
the identity \eqref{polyinterpolation} in item (2) can be equivalently expressed as
\be \label{sdan:lpm}
\sd_a^n (\pp*\vgu_a)=(\pp(2^{-n}\cdot))*\vgu_a\in \PR_{M,\vgu_a},\qquad
\forall\, \pp\in \PL_M.
\ee
Hence, item (1) of \cref{thm:sd:sr} is satisfied with $v=\vgu_a$ by noting $1*\vgu_a=\wh{\vgu_a}(0)$. Hence, items (3) and (4) of \cref{thm:sd:sr} hold, that is,
the matrix mask $a$ has order $M+1$ sum rules for some matching filter $\tilde{\vgu}_a:=b*\vgu_a
\in \dlrs{0}{1}{r}$ with $b\in \dlp{0}$ satisfying $\wh{b}(0)=1$, and $\sd_a (\pq_\mu*\tilde{\vgu}_a)=2^{-|\mu|} \pq_\mu*\tilde{\vgu}_a$ for all $\mu\in \ind_M$, where $\pq_\mu(x):=\frac{x^\mu}{\mu!}$. Hence, we conclude from \eqref{sdan:lpm} with $n=1$ that for all $\mu\in \ind_M$,
\[
\pq_\mu*(b*\vgu_a)=\pq_\mu*\tilde{\vgu}_a=2^{|\mu|} \sd_a((\pq_\mu*b)*\vgu_a)=
2^{|\mu|} ((\pq_\mu*b)(2^{-1}\cdot))*\vgu_a=\pq_\mu*(b_2*\vgu_a),
\]
where $\wh{b_2}(\xi):=\wh{b}(2\xi)$ and we used $2^{|\mu|}\pq_\mu(2^{-1}x)=\frac{x^\mu}{\mu!}=\pq_\mu(x)$. We conclude from the above identity and \eqref{conv:poly} that $\wh{b}(\xi)\wh{\vgu_a}(\xi)=\wh{b}(2\xi)\wh{\vgu_a}(\xi)+\bo(\|\xi\|^{M+1})$ as $\xi\to 0$. Because $\nu_1=\mathbf{0}$, the first entry of $\wh{\vgu_a}(0)$ is $1$ and hence, we must have $\wh{b}(\xi)=\wh{b}(2\xi)+\bo(\|\xi\|^{M+1})$ as $\xi \to 0$. Using the Taylor expansion of $\wh{b}$ at $\xi=0$ and noting $\wh{b}(0)=1$, we conclude that $\wh{b}(\xi)=1+\bo(\|\xi\|^{M+1})$ as $\xi\to 0$. This proves
$\wh{\tilde{\vgu}_a}(\xi)=\wh{b}(\xi)\wh{\vgu_a}(\xi)=\wh{\vgu_a}(\xi)+\bo(\|\xi\|^{M+1})$ as $\xi\to 0$.
Hence, we proved item (1).

(2)$\iff$(3) is straightforward by observing $[\pp(2^{-n}\cdot)]^{(\nu_\ell)}=2^{-|\nu_\ell|n} \pp^{(\nu_\ell)}(2^{-n}\cdot)$.

We deduce from $w_n=[\pp^{(\nu_1)}(2^{-n}(\cdot+\tau_1)),\ldots,
\pp^{(\nu_r)}(2^{-n}(\cdot+\tau_r))]$ in item (3) that the function $\eta$ in \eqref{hsd:converg:face} must be $\pp$. Now item (i) follows directly from item (3). Using item (3), we have
\[
w_n(2k+2\tau_\ell-\tau_{\theta(\ell)})e_{\theta(\ell)}=
\pp^{(\nu_{\theta(\ell)})}(2^{-n}(2k+2\tau_\ell-\tau_{\theta(\ell)}+\tau_{\theta(\ell)}))
=\pp^{(\nu_\ell)}(2^{1-n}(k+\tau_\ell))=w_{n-1}(k)e_\ell,
\]
where we used $\nu_{\theta(\ell)}=\nu_\ell$. This proves item (ii).
\ep

We say that a matrix mask $a\in \dlrs{0}{r}{r}$ has \emph{order $M+1$ linear-phase moments of type $\ind$ and translation $T$} if the mask $a$ has order $M+1$ sum rules with a matching filter $\vgu_a\in \dlrs{0}{1}{r}$ satisfying \eqref{int:vgua}, that is, item (1) of \cref{thm:lpm} holds.
Because all $w_n$ in item (ii) of \cref{thm:lpm} are vector polynomials, $w_n$ can be defined on $\dR$ instead of just $\dZ$ and hence
the condition $2\tau_\ell-\tau_{\theta(\ell)}\in \dZ$ in \eqref{theta} can be dropped for the mapping $\theta$ in item (ii) of \cref{thm:lpm}.

\section{Existence and Examples of Generalized Hermite Subdivision Schemes}
\label{sec:ex}

To illustrate the theoretical results in previous sections,
in this section we shall first prove the existence of convergent generalized Hermite subdivision schemes of type $\ind$ for any given ordered multiset $\ind$ and then provide a few examples of symmetric generalized Hermite subdivision schemes.

\subsection{Existence of convergent generalized Hermite subdivision schemes}
To prove the existence, generalizing \cite[Proposition~6.2]{han09}, we have the following result.

\begin{lemma}\label{lem:convert:vec}
Let $\varphi=[\varphi_1,\ldots,\varphi_L]^\tp$ be a refinable vector function on $\dR$ such that $\wh{\varphi}(2\xi)=\wh{A}(\xi)\wh{\varphi}(\xi)$ for a matrix mask $A\in \dlrs{0}{L}{L}$ and $A$ has order $m+1$ sum rules with a matching filter $\vgu_A\in \dlrs{0}{1}{L}$ and $\wh{\vgu_A}(0)\ne 0$.
Let $r\in \N$ and take any $d\times d$ integer matrix $N$ satisfying $|\det (N)|=r$. Define $\Gamma_N:=\{\gamma_1,\ldots,\gamma_r\}:=(N[0,1)^d) \cap \dZ$ to be a set of all distinct representatives of cosets in $\dZ/(N\dZ)$. Then the new vector function $\phi:=[\varphi(N\cdot-\gamma_1)^\tp, \ldots, \varphi(N\cdot-\gamma_r)^\tp]^\tp$ must be a refinable vector function satisfying $\wh{\phi}(2\xi)=\wh{a}(\xi)\wh{\phi}(\xi)$,
$\sm_p(a)=\sm_p(A)$ for all $1\le p\le \infty$, and the mask $a\in \dlrs{0}{(rL)}{(rL)}$ has order $m+1$ sum rules with a matching filter $\vgu_a\in \dlrs{0}{1}{(rL)}$ given by
\be \label{vgua:A}
\wh{\vgu_a}(\xi)=\left[e^{i\gamma_1\cdot N^{-\tp} \xi} \wh{\vgu_A}(N^{-\tp}\xi),\ldots,
e^{i\gamma_r\cdot N^{-\tp} \xi} \wh{\vgu_A}(N^{-\tp}\xi)\right]+\bo(\|\xi\|^{m+1}),\quad \xi\to 0,
\ee
where the matrix mask $a\in \dlrs{0}{(rL)}{(rL)}$ is defined in the following way: the $L\times L$ block $(j,k)$-entry $[a(n)]_{j,k}$ of the matrix $a(n)$ is given by $A(Nn-2\gamma_j+\gamma_k)$
 for all $n\in \dZ$ and $j,k=1,\ldots, r$.
\end{lemma}

\bp For $j=1,\ldots,r$,
we define a vector function $\phi_j:=\varphi(N\cdot-\gamma_j)$ and a mapping $\theta_j: \Gamma_N \rightarrow \Gamma_N$ such that $\theta_j(\gamma_\ell)$ is the unique element in $\Gamma_N$ such that $2\gamma_j+\gamma_\ell-\theta_j(\gamma_\ell)\in N\dZ$ for $\ell=1,\ldots,r$.
Hence, $n_{j,\ell}:=N^{-1}(2\gamma_j+\gamma_\ell-\theta_j(\gamma_\ell))\in \dZ$ for $j,k=1,\ldots,r$.
Using the refinement equation $\varphi=2^d \sum_{k\in \dZ} A(k) \varphi(2\cdot-k)$ and noting that $2\gamma_j+\gamma_\ell=\theta_j(\gamma_\ell)+N n_{j,\ell}$, we deduce that
\begin{align*}
\phi_j&=\varphi(N\cdot-\gamma_j)
=2^d \sum_{\ell=1}^r \sum_{k\in \dZ}
A(Nk+\gamma_{\ell})\varphi(2(N\cdot-\gamma_j)-(Nk+\gamma_{\ell}))\\
&=2^d \sum_{\ell=1}^r
\sum_{k\in \dZ}
A(Nk+\gamma_{\ell}) \varphi(N(2\cdot-k-n_{j,\ell})-\theta_j(\gamma_{\ell}))\\
&=2^d \sum_{\ell=1}^r
\sum_{n\in \dZ}
A(Nn-Nn_{j,\ell}+\gamma_\ell) \varphi(N(2\cdot-n)-\theta_j(\gamma_\ell))\\
&=2^d \sum_{\ell=1}^r
\sum_{n\in \dZ}
A(Nn-2\gamma_j+\theta_j(\gamma_\ell)) \varphi(N(2\cdot-n)-\theta_j(\gamma_\ell)),
\end{align*}
where we used the identity $Nn_{j,\ell}=2\gamma_j+\gamma_\ell-\theta_j(\gamma_\ell)$.
Note that the mapping $\theta_j$ is bijective on $\Gamma_N$.
Using the substitution $\gamma_{k}=\theta_j(\gamma_\ell)$, we deduce from the above identity that
\[
\phi_j=2^d \sum_{k=1}^r \sum_{n\in \dZ}
A(Nn-2\gamma_j+\gamma_{k}) \phi_{k}(2\cdot-n),\qquad j,k=1,\ldots,r.
\]
Due to $\phi=[\phi_1^\tp,\ldots,\phi_r^\tp]^\tp$,
this proves that $\phi$ satisfies the refinement equation $\phi=2^d \sum_{n\in \dZ} a(n) \phi(2\cdot-n)$ such that
the $L\times L$ block $(j,k)$-entry $[a(n)]_{j,k}$ of the matrix $a(n)$ is given by $A(Nn-2\gamma_j+\gamma_k)$.

We now prove that $a$ must have order $m+1$ sum rules with the matching filter $\vgu_a$ in \eqref{vgua:A}.
Using the definition of $\vgu_a$ in \eqref{vgua:A}, for $\omega\in \Gamma:=[0,1]^d\cap \dZ$ as in \eqref{Gamma}, we have
\begin{align*}
\sum_{j=1}^r [\wh{\vgu_a}(2N^\tp \xi)]_j [\wh{a}(N^\tp \xi+\pi \omega)]_{j,k}
&=\sum_{j=1}^r \sum_{n\in \dZ}
e^{i2\gamma_j\cdot \xi}\wh{\vgu_A}(2\xi)
A(Nn-2\gamma_j+\gamma_k) e^{-i n\cdot (N^\tp \xi+\pi \omega)}+\bo(\|\xi\|^{m+1})\\
&=e^{i\gamma_k\cdot\xi} \wh{\vgu_A}(2\xi) B_{j,k,\omega}(\xi)+\bo(\|\xi\|^{m+1}),
\end{align*}
as $\xi\to 0$,  where
$B_{j,k,\omega}(\xi):=\sum_{j=1}^r \sum_{n\in \dZ}
A(Nn-2\gamma_j+\gamma_k) e^{-i(Nn-2\gamma_j+\gamma_k)\cdot \xi}
e^{-in\cdot \pi \omega}$.
Note that any $n\in \dZ$ can be uniquely expressed as $n=2n'+\gamma'$ with $n'\in \dZ$ and $\gamma'\in \Gamma$. Hence,
\begin{align*}
B_{j,k,\omega}(\xi)
&=\sum_{\gamma'\in \Gamma}
\sum_{n'\in \dZ} \sum_{j=1}^r
A(2Nn'-2\gamma_j+N\gamma'+\gamma_k)
e^{-i(2Nn'-2\gamma_j+N\gamma'+\gamma_k)}
e^{-i\gamma'\cdot \pi \omega}\\
&=\sum_{\gamma'\in \Gamma}
\Big(\sum_{q\in \dZ}
A(2q+N\gamma'+\gamma_k)
e^{-i(2q+N\gamma'+\gamma_k)\cdot\xi} \Big ) e^{-i\gamma'\cdot \pi \omega}\\
&=\sum_{\gamma'\in \Gamma}
\wh{A^{[N\gamma'+\gamma_k]}}(2\xi)
e^{-i(N\gamma'+\gamma_k)\cdot\xi}e^{-i \gamma' \cdot \pi \omega},
\end{align*}
where we used the substitution $q=Nn'-\gamma_j$ and the fact that any $q\in \dZ$ can be uniquely expressed as $q=Nn-\gamma$ for $n'\in \dZ$ and $\gamma\in \Gamma_N$.
Since $A$ has order $m+1$ sum rules with the matching filter $\vgu_A$, by \eqref{sr:coset} with $a$ being replaced by $A$, we deduce from the above identities that as $\xi\to 0$,
\begin{align*}
\sum_{j=1}^r [\wh{\vgu_a}(2N^\tp \xi)]_j [\wh{a}(N^\tp \xi+\pi \omega)]_{j,k}
&=e^{i\gamma_k\cdot\xi} \wh{\vgu_A}(2\xi) B_{j,k,\omega}(\xi)+\bo(\|\xi\|^{m+1})\\
&=e^{i\gamma_k\cdot\xi} \sum_{\gamma'\in \Gamma}
\wh{\vgu_A}(2\xi) \wh{A^{[N\gamma'+\gamma_k]}}(2\xi)
e^{-i(N\gamma'+\gamma_k)\cdot\xi}e^{-i\gamma'\cdot \pi \omega}+\bo(\|\xi\|^{m+1})\\
&=2^{-d}e^{i\gamma_k\cdot\xi}
 \wh{\vgu_A}(\xi)
\sum_{\gamma'\in \Gamma} e^{-i\gamma'\cdot \pi \omega}+\bo(\|\xi\|^{m+1}).
\end{align*}
Noting that $\sum_{\gamma'\in \Gamma} e^{-i\gamma'\cdot \pi \omega}=2^d \td(\omega)$ for all $\omega\in \Gamma$, we conclude from the above identity that
\[
[\wh{\vgu_a}(2N^\tp \xi) \wh{a}(N^\tp \xi+\pi \omega)]_k=
\sum_{j=1}^r [\wh{\vgu_a}(2N^\tp \xi)]_j [\wh{a}(N^\tp \xi+\pi\omega)]_{j,k}
=e^{i\gamma_k\cdot\xi} \wh{\vgu_A}(\xi) \td(\omega)+\bo(\|\xi\|^{m+1})
\]
as $\xi \to 0$.
Replacing $\xi$ by $N^{-\tp}\xi$ in the above identity and using the definition of $\vgu_a$ in \eqref{vgua:A}, we conclude that the mask $a$ has order $m+1$ sum rules satisfying \eqref{sr} with the matching filter $\vgu_a$ in \eqref{vgua:A}. Now the identity $\sm_p(a)=\sm_p(A)$ can be directly checked using the definition in \eqref{sma}.
\ep

To prove the existence of convergent generalized Hermite subdivision schemes of type $\ind$ for any given ordered multiset $\ind$, we recall the definition of the B-spline functions. For $n\in \N$,
the B-spline function $B_n$ of order $n$ is defined to be
\be \label{bspline}
B_1:=\chi_{(0,1]} \quad \mbox{and}\quad
B_n:=B_{n-1}*B_1=\int_0^1 B_{n-1}(\cdot-x) dx.
\ee
Then $B_n\in \CH{n-2}$ with support $[0,n]$, $\wh{B_n}(\xi)=(\frac{1-e^{-i\xi}}{i\xi})^n$,
and
$B_n|_{(k,k+1)}$ is a nonnegative polynomial of degree $n-1$ for every $k\in \Z$. Moreover, the B-spline function $B_n$ is refinable by satisfying
\be \label{bspline:mask}
\wh{B_n}(2\xi)=\wh{a_n^B}(\xi)\wh{B_n}(\xi)
\quad \mbox{with}\quad
\wh{a^B_n}(\xi):=2^{-n}(1+e^{-i\xi})^n.
\ee
Note that $\sm_p(B_n)=\sm_p(a^B_n)=n-1+1/p$ for $1\le p\le \infty$. Moreover,
the integer shifts of $B_n$ are linearly independent.

We are now ready to prove the existence of convergent generalized Hermite subdivision schemes.

\begin{theorem}\label{thm:existence}
Let $m\in \NN$ and $r\in \N$ with $r>1$.
Let $\ind=\{\nu_1,\ldots,\nu_r\} \subseteq \ind_m$
as in \eqref{indset} with $\nu_1=\mathbf{0}\in \dNN$ and $T=\{\tau_1,\ldots,\tau_r\}\subseteq \R^d$
be arbitrarily given ordered multisets.
For any mask $\mathring{a}\in \dlrs{0}{r}{r}$ with $\sm_\infty(\mathring{a})>m$, one can always constructively derive a generalized Hermite mask $a\in \dlrs{0}{r}{r}$ of type $\ind$ such that $\sm_\infty(a)=\sm_\infty(\mathring{a})>m$,
the mask $a$ has order $m+1$ linear-phase moments of type $\ind$ and translation $T$, and
the generalized Hermite subdivision scheme of type $\ind$ with mask $a$ is convergent with limit functions in $\dCH{m}$. Moreover, such a desired generalized Hermite mask $a\in \dlrs{0}{r}{r}$ of type $\ind$ always exists and can be constructed such that its basis vector function $\phi$  is a spline refinable vector function in $\dCH{m}$ and $\phi$ has linearly independent integer shifts.
\end{theorem}

\bp
By \cite[Theorem~4.3]{han03}, $\sm_\infty(\mathring{a})>m$ implies that the mask $\mathring{a}$ must have order $m+1$ sum rules with some matching filter $\mathring{\vgu}\in \dlrs{0}{1}{r}$ with $\wh{\mathring{\vgu}}(0)\ne 0$ and its refinable vector function $\mathring{\phi}$ belongs to $\dCH{m}$.
%If $r=1$, then the claim holds with $a=\mathring{a}$ by $\ind=\{0\}$ in \eqref{indset}. So, we assume $r>1$.
Let $\vgu_a\in \dlrs{0}{1}{r}$ satisfy
$\wh{\vgu_a}(\xi)=[(i\xi)^{\nu_1}e^{i\tau_1\cdot\xi},\ldots,
(i\xi)^{\nu_r}e^{i\tau_r\cdot\xi}]+\bo(\|\xi\|^{m+1})$ as $\xi \to 0$.
By \cref{lem:vgu} and $r>1$,
there must exist a strongly invertible sequence $U\in \dlrs{0}{r}{r}$ such that
\be \label{vgua:new}
\wh{\mathring{\vgu}}(\xi)
=\wh{\vgu_a}(\xi)\wh{U}(\xi)+\bo(\|\xi\|^{m+1}),\quad \xi \to 0.
\ee
Define $\wh{a}(\xi):=\wh{U}(2\xi)\wh{\mathring{a}}(\xi)\wh{U}(\xi)^{-1}$ and $\wh{\phi}(\xi):=\wh{U}(\xi)\wh{\mathring{\phi}}(\xi)$ as in \eqref{newa}.
Since $U$ is strongly invertible, we have $\sm_p(a)=\sm_p(\mathring{a})$ and
$\sm_p(\phi)=\sm_p(\mathring{\phi})$
for all $1\le p\le \infty$. Moreover, the mask $a$ has order $m+1$ sum rules with the matching filter $\vgu_a$.
By \cref{thm:converg} and $\sm_\infty(a)>m$, the generalized Hermite subdivision scheme of type $\ind$ with mask $a$ must be convergent with limit functions in $\dCH{m}$.
By \cref{thm:lpm}, the mask $a$ must have order $m+1$ linear-phase moments of type $\ind$ and translation $T$.

To show the existence of such desired masks $\mathring{a}$ and $a$,
we consider $\varphi:=\otimes^d B_{m+2}$ and $A:=\otimes^d a^B_{m+2}$, where $\otimes^d B_{m+2}$ is the tensor product spline defined by $[\otimes^d B_{m+2}](x_1,\ldots,x_d)=B_{m+2}(x_1)\cdots B_{m+2}(x_d)$. Note that $\sm_\infty(A)=\sm_\infty(a^B_{m+2})=m+1>m$.
Using \cref{lem:convert:vec}, we can construct a compactly supported vector refinable function $\mathring{\phi}$ and a mask $\mathring{a}\in \dlrs{0}{r}{r}$ such that $\wh{\mathring{\phi}}(2\xi)=\wh{\mathring{a}}(\xi)\wh{\mathring{\phi}}(\xi)$ and $\sm_\infty(\mathring{a})=\sm_a(A)=m+1>m$.
Moreover, since the integer shifts of $B_{m+2}$ are linearly independent, using the definition of linear independence,
we observe that the integer shifts of $\varphi$ and its derived $\mathring{\phi}$ in \cref{lem:convert:vec} must be linearly independent. Since $U$ in \eqref{vgua:new} is strongly invertible, the integer shifts of $\phi$ must be linearly independent as well and $\phi$ is obviously a spline vector function.
\ep

\subsection{Examples}
Though convergent generalized Hermite subdivision schemes can be constructed theoretically through \cref{thm:existence}, the support of the matrix mask $a$ in \cref{thm:existence} is often very large and the mask $a$ lacks symmetry. Therefore, it is not that useful to construct convergent generalized Hermite subdivision schemes by \cref{thm:existence} for practical purposes.
We  now discuss how to construct particular generalized Hermite subdivision schemes with short support, symmetry, and high smoothness.
Symmetry property of scalar multivariate subdivision schemes has been well understood (e.g., see \cite{han03smaa,hanbook} and references therein). However, symmetry property of vector/matrix subdivision schemes is a technical issue and has been only discussed for very special subdivision schemes, e.g., see \cite{hyp04,hyx05,hy06,hz09}. To present several examples of generalized Hermite subdivision schemes with the symmetry property, let us briefly discuss the symmetry property of generalized Hermite subdivision schemes.
We say that a finite set $G$ of $d\times d$ integer matrices is \emph{a symmetry group} in $\R^d$ if $|\det (E)|=1$ for all $E\in G$ and $G$ forms a group under the matrix multiplication. Typical examples of symmetry groups are $G=\{1,-1\}$ for $d=1$ and the following symmetry groups for $d=2$:
\begin{align*}
&D_4:=\left\{ \pm \begin{bmatrix} 1&0\\ 0&1\end{bmatrix},
\pm \begin{bmatrix} 1&0\\ 0&-1\end{bmatrix},
\pm \begin{bmatrix} 0&1\\ 1&0\end{bmatrix},
\pm \begin{bmatrix} 0&1\\ -1&0\end{bmatrix}\right\},\\
&D_6:=\left\{ \pm \begin{bmatrix} 1&0\\ 0&1\end{bmatrix},
\pm \begin{bmatrix} 0&-1\\ 1&-1\end{bmatrix},
\pm \begin{bmatrix} -1&1\\ -1&0\end{bmatrix}
\pm \begin{bmatrix} 0&1\\ 1&0\end{bmatrix},
\pm \begin{bmatrix} 1&-1\\ 0&-1\end{bmatrix},
\pm \begin{bmatrix} -1&0\\ -1&1\end{bmatrix}\right\},
\end{align*}
as well as their subgroups.
The symmetry groups $D_4$ and $D_6$ are often used for the quadrilateral mesh and the triangular mesh in CAGD, respectively, e.g., see \cite[Section~7.2]{hanbook}
%and \cite{fhs22}
for details.

Let $\phi=[\phi_1,\ldots,\phi_r]^\tp$ be a compactly supported refinable vector function satisfying $\wh{\phi}(2\xi)=\wh{a}(\xi)\wh{\phi}(\xi)$ for some matrix mask $a\in \dlrs{0}{r}{r}$.
Let an ordered multiset $\ind=\{\nu_1,\ldots,\nu_r\}$ as in \eqref{indset} and an ordered multiset $T=\{\tau_1,\ldots,\tau_r\}\subseteq \dR$ be the symmetry centers of the vector function $\phi$. Because the symmetry center of $\phi_\ell$ is $\tau_\ell$, the symmetry property of $\phi_\ell$ naturally requires
\[
\phi_\ell (E(\cdot-\tau_\ell)+\tau_\ell)=[S_E]_{\ell,1} \phi_1(\cdot+\tau_1-\tau_\ell)+\cdots+[S_E]_{\ell,r}\phi_r(\cdot+\tau_r-\tau_\ell),\qquad \forall\, E\in G, \ell=1,\ldots,r,
\]
where $S_E$ is an invertible $r\times r$ constant matrix. Using the Fourier transform, we observe that the above symmetry property is equivalent to
$\wh{\phi}(E^{-\tp}\xi)=D_T(-E^{-\tp}\xi) S_E D_T(\xi)\wh{\phi}(\xi)$ for all $E\in G$, i.e.,
\be \label{sym:phi}
\wh{\phi}(E^\tp \xi)=D_T(-E^\tp\xi)
S_{E^{-1}} D_T(\xi) \wh{\phi}(\xi),\quad \forall\, E\in G
\quad \mbox{with}\quad
D_T(\xi):=\mbox{diag}( e^{i\tau_1\cdot\xi},\ldots, e^{i\tau_r\cdot \xi}).
\ee
Note that $[D_T(\xi)]^{-1}=D_T(-\xi)$.
Using the refinement equation $\wh{\phi}(2\xi)=\wh{a}(\xi)\wh{\phi}(\xi)$, we observe that \eqref{sym:phi} naturally requires
the following symmetry property on the matrix mask $a$:
\be \label{sym:mask}
\wh{a}(E^\tp \xi)=
D_T(-2E^\tp\xi) S_{E^{-1}}
D_T(2\xi) \wh{a}(\xi) D_T(-\xi) S_{E^{-1}}^{-1} D_T(E^\tp\xi),\qquad \forall\, E\in G.
\ee
Now we assume that the mask $a$ has order $M+1$ linear-phase moments of type $\ind$ and translation $T$ as in \cref{thm:lpm}, i.e., the mask $a$ has order $M+1$ sum rules with a matching filter $\vgu_a$ satisfying \eqref{int:vgua}. As we explained in \cref{thm:hsd:mask}, under the natural condition \eqref{eig:one} with $m=M$, the matching filter $\vgu_a$ is essentially uniquely determined by the mask $a$ through \eqref{vgua:val}. Using \eqref{int:vgua} and \eqref{sym:mask}, we deduce from $\wh{\vgu_a}(2\xi)\wh{a}(\xi)=\wh{\vgu_a}(\xi)+\bo(\|\xi\|^{M+1})$ as $\xi \to 0$ that $\vgu_a$ has the following symmetry property:
\be \label{sym:vgua}
\wh{\vgu_a}(E^\tp\xi)=\wh{\vgu_a}(\xi)D_T(-\xi) S_{E^{-1}}^{-1} D_T(E^\tp \xi)+\bo(\|\xi\|^{M+1}),\quad \xi \to 0, \; \forall\; E\in G.
\ee
Noting  that the integer $M$ is arbitrary and
by \eqref{int:vgua} $\wh{\vgu_a}(\xi)=\left[ (i\xi)^{\nu_1},\ldots,(i\xi)^{\nu_r}\right] D_T(\xi)+\bo(|\xi|^{M+1})$ as $\xi\to 0$,
we see that
%the matrices $S_{E^{-1}}, E\in G$ in
\eqref{sym:vgua} becomes
\be \label{sym:SE}
S_{E}=S(E, \ind)
\quad\mbox{with}\quad
\left[(iE^\tp \xi)^{\nu_1},\ldots,
(iE^\tp \xi)^{\nu_r}\right]
=\left[(i\xi)^{\nu_1},\ldots,
(i\xi)^{\nu_r}\right] S(E,\ind), \quad E\in G.
\ee
Note that $S(E,\ind)_{j,k}=0$ if $|\nu_j|\ne |\nu_k|$ and hence $i\xi$ in \eqref{sym:SE} can be replaced by $\xi$.
Because the matching filter $\vgu_a$ for a generalized Hermite subdivision scheme must satisfy \eqref{vgua:hsd} by \cref{thm:hsd:mask}, the symmetry property in \eqref{sym:phi}, \eqref{sym:mask} and \eqref{sym:vgua} still holds for generalized Hermite subdivision schemes.
The matrix $S(E,\ind)$ defined in \eqref{sym:SE} is uniquely determined if all the elements in $\ind$ are not repeated. But when elements in the ordered multiset $\ind=\{\nu_1,\ldots,\nu_r\}$ can repeat, the definition of $S(E,\ind)$ is not unique and has free parameters which can be used for compatibility with the multiset $T$ and the mapping $\theta$ in \eqref{theta}. We shall not address this issue further. For $\ind=\ind_m$, by \cite[Lemma~5.1]{hz09}, the matrix $S(E,\ind_m)$ in \eqref{sym:SE} is the same as the matrix $S(E,\ind_m)$ defined in \cite[(2.1)]{han03} and \cite[(2.10)]{han03smaa}, which is used in
\cite[(2.5)]{hyp04} and \cite[(3.1)]{hz09} for studying interpolatory Hermite subdivision schemes.
The univariate spline generalized Hermite interpolants of type $\ind$ and translation $T$ in \eqref{phi:spline} and their associated matrix masks have the symmetry property. Hence, the above discussion on the symmetry property of generalized Hermite subdivision schemes agrees with and explains
\cite[Proposition~2.8]{hyp04},
\cite[Proposition~2.3]{hyx05}, \cite[Theorem~2.5]{hy06}, and \cite[Theorem~3.3]{hz09} for all these special cases.

A function $f$ on $\dR$ is called a spline
(or a piecewise polynomial with finite polygonal domains) in \cite{hm07}
if there exist mutually disjoint open sets $U_j, j=1,\ldots, J$ such that $f$ vanishes outside the closure of $\cup_{j=1}^J U_j$ and $f$ is a $d$-variate polynomial on each set $U_j$ for $j=1,\ldots,J$, where each $U_j$ is a finite intersection of some half-spaces $\{x\in \dR: x\cdot t>c\}$ with $t\in \dR$ and $c\in \R$.
%For a compactly supported vector function $\phi$ on $\dR$, we say that $\phi$ is a spline vector if every entry of $\phi$ is a piecewise polynomial with finite polynomial domains.
Refinable vector functions which are spline vectors are of particular theoretical interest due to desired properties of splines in applications.
Up to a multiplicative constant,
a refinable vector function $\phi$ with mask  $a\in \dlrs{0}{r}{r}$
is determined through the refinement equation \eqref{refeq:phi} but $\phi$ often has no analytic or explicit expressions at all.
It is a challenging task to find particular choices of the free parameters in a given matrix mask $a$ such that its associated refinable vector function $\phi$ is indeed a spline vector. Quite often, all its associated refinable vector functions $\phi$ are not spline vectors at all, regardless of how many free parameters in its matrix mask (see \cite{hm07}).
If a refinable vector function $\phi$ with a mask $a\in \dlrs{0}{r}{r}$ is indeed a spline vector and $\phi$ has linearly independent integer shifts, then \cite[Theorem~6]{hm07} or \cite[Theorem~6.1.9]{hanbook} tells us that all the nonzero eigenvalues of the transition operator $\mathcal{T}_a$ must take the form $2^{-j}, j\in \NN$, where the transition operator $\mathcal{T}_a: \dlrs{}{1}{r}\rightarrow \dlrs{0}{1}{r}$ is defined to be
\be \label{tz}
[\mathcal{T}_a v](n):=2^d\sum_{k\in \dZ} v(k) \ol{a(k-2n)}^\tp,\qquad n\in \dZ, v\in \dlrs{0}{1}{r}.
\ee
To find spline refinable vector functions $\phi$ from a given matrix mask $a\in \dlrs{0}{r}{r}$ with free parameters,
we first symbolically compute all the nonzero eigenvalues of $\mathcal{T}_a$.
Then we could search for particular parameters in a matrix mask $a$ such that the necessary condition on special eigenvalues of $\mathcal{T}_a$  holds. Finally, one can use the refinement equation \eqref{refeq:phi} to guess and check in a brute-force way about whether its refinable vector function $\phi$ is indeed a spline vector or not. It is often difficult to symbolically compute all the eigenvalues of $\mathcal{T}_a$ and then pick particular choices of parameters for the necessary condition if a matrix mask $a$ has a few free parameters and large supports.
A few examples of spline univariate vector refinable functions have been reported in \cite{hanbook,han20,hm07}.
For high dimensions, except the scalar refinable functions which are box splines,
so far the only known bivariate spline refinable vector functions in the literature are the $D_6$-symmetric Powell-Sabin spline (see \cite[Proposition~3.3]{hyp04} and \cite[Theorem~8]{hm07}) and
the $D_4$-symmetric Sibson spline (see \cite[Example~3.5]{hyp04} and \cite[Theorem~9]{hm07}), whose masks $a\in (l(\Z^2))^{3\times 3}$ have order $3$ sum rules
with support $[-1,1]^2$, and whose basis vector functions are Hermite interpolants of degree one. See \cite{hyp04,hm07,hanbook,han20} for details about some examples of spline refinable Hermite interpolants and how to find them.

Even in dimension one our results significantly generalize and extend known results on univariate Hermite subdivision schemes. Our introduction and analysis of multivariate generalized Hermite subdivision schemes also include tensor products of univariate generalized Hermite subdivision schemes as special cases.
Using the symmetry group $G=\{-1,1\}$, we first present two examples of univariate symmetric generalized Hermite subdivision schemes that are not known yet in the literature.

Recall that $\sr(a)$ is the largest possible integer $m\in \NN$ such that $a$ has order $m$ sum rules. Similarly, $\lpm(a)$ is the largest possible integer $n\in \NN$ such that $a$ has order $n$ linear-phase moments of type $\ind$ and translation $T$ in \cref{thm:lpm}. For convenience of discussion, we shall use $\fs(a)$ to denote the smallest interval in $d=1$ and the smallest rectangle in $d=2$ such that $a$ vanishes outside $\fs(a)$.

\begin{example}\label{ex1}
\normalfont
Let $d=1$, $\ind=\{0,2\}$ and $T=\{0,0\}$.
We consider three families of generalized Hermite subdivision schemes of type $\ind=\{0,2\}$ (more precisely, Birkhoff subdivision schemes). Using \cref{thm:lpm} with $\wh{\vgu_a}(\xi)=[1,(i\xi)^2]+\bo(|\xi|^6)$ as $\xi\to 0$, we find that all the symmetric generalized Hermite masks $a \in \lrs{0}{2}{2}$ of type $\ind$ with
$\fs(a)=[-3,3]$ and $\lpm(a)\ge 6$ are given by
\be \label{mask:birkhoff}
\begin{split}
a=&\left\{
\begin{bmatrix}
\frac{5}{128}-t_1 &-\frac{3}{32}-t_2\\[0.1em]
\frac{1}{12}t_1 &\frac{1}{12}t_2 \end{bmatrix},
\begin{bmatrix}
-t_3 &-t_4\\[0.1em]
\frac{1}{12}t_3 &\frac{1}{12}t_4 \end{bmatrix},
\begin{bmatrix}
\frac{27}{128}+t_1 &\frac{3}{32}+ t_2\\[0.1em]
-\frac{9}{128}+\frac{11}{12}t_1 &\frac{5}{32}+\frac{11}{12}t_2 \end{bmatrix},
\begin{bmatrix}
\frac{1}{2}+2t_3 &2t_4\\[0.1em]
\frac{5}{6}t_3 &\frac{1}{8}+\frac{5}{6}t_4 \end{bmatrix},\right.\\
&\left.\;\;
\begin{bmatrix}
\frac{27}{128}+t_1 &\frac{3}{32}+ t_2\\[0.1em]
-\frac{9}{128}+\frac{11}{12}t_1 &\frac{5}{32}+\frac{11}{12}t_2 \end{bmatrix},
\begin{bmatrix}
-t_3&-t_4\\[0.1em]
\frac{1}{12}t_3 &\frac{1}{12}t_4 \end{bmatrix},
\begin{bmatrix}
\frac{5}{128}-t_1 &-\frac{3}{32}- t_2\\[0.1em]
\frac{1}{12}t_1 &\frac{1}{12}t_2 \end{bmatrix}
\right\}_{[-3,3]},
\end{split}
\ee
where $t_1,\ldots,t_4\in \R$.
Hence the Birkhoff subdivision schemes with masks $a$ have the polynomial-interpolation property of order $6$ by interpolating all polynomials of degree less than $6$ in \cref{thm:lpm}.
For $t_1=\frac{5}{128}, t_2=-\frac{3}{16}$ and $t_3=t_4=-\frac{3}{32}$, $\sm_2(a)\approx 4.3522$ is almost the highest and hence, $\sm_\infty(a)\ge \sm_2(a)-0.5>3$. By \cref{thm:converg,thm:lpm}, the Birkhoff subdivision scheme with mask $a$ is convergent with limit functions in $\CH{3}$ and has the polynomial-interpolation property of order $6$. Moreover, $\sr(a)=10$ if and only if $t_1 = \frac{91}{1024}$, $t_2 = -\frac{15}{64}$, $t_3 = -\frac{17}{512}$, and $t_4 = -\frac{9}{64}$, for which $\sm_2(a)\approx 2.53079$ and $\wh{\vgu_a}(\xi)=[1-\frac{17}{12096}(i\xi)^6+\frac{1}{4320}(i\xi)^8, (i\xi)^2-\frac{1}{40}(i\xi)^6+\frac{1}{252}(i\xi)^8]+\bo(|\xi|^{10})$ as $\xi \to 0$.

The identity mapping $\theta: \{1,2\}\rightarrow \{1,2\}$ obviously satisfies \eqref{theta} with $T=\{0,0\}$.
All symmetric interpolatory generalized Hermite masks $a\in \lrs{0}{2}{2}$ of type $\ind$ and translation $T$ (i.e., interpolatory Birkhoff masks) with $\fs(a)=[-3,3]$ and $\sr(a)\ge 6$ must be given in \eqref{mask:birkhoff} with $t_3=t_4=0$.
For $t_1=\frac{25}{256}$ and $t_2=-\frac{1}{4}$, $\sm_2(a)\approx 2.6943$ is nearly the highest and
hence $\sm_\infty(a)\ge \sm_2(a)-0.5>2$. By \cref{thm:int,thm:converg},
the interpolatory Birkhoff subdivision scheme with mask $a$ is convergent with limit functions in $\CH{2}$ and its refinable vector function $\phi=[\phi_1,\phi_2]^\tp\in (\CH{2})^2$ is a refinable Birkhoff interpolant satisfying $\phi_1(k)=\phi_2''(k)=\td(k)$ and $\phi''_1(k)=\phi_2(k)=0$ for all $k\in \Z$.
Moreover, for $t_3=t_4=0$, $\sr(a)=8$ if and only if $t_1=-\frac{27}{256}$ and $t_2=\frac{33}{128}$, for which $\sm_2(a)\approx 0.02797$.

The above two families of generalized Hermite masks show that the interpolation conditions in \eqref{interpolation} and \eqref{int:mask} are much stronger than the linear-phase moment conditions for the polynomial-interpolation condition in \cref{thm:lpm}. If we give up the polynomial-interpolation property, then we can achieve better smoothness even with shorter support. We find that all symmetric generalized Hermite mask $a\in \lrs{0}{2}{2}$ of type $\ind$ with $\fs(a)=[-2,2]$ and $\sr(a)\ge 8$ are given by
%
%{\footnotesize{
%\be \label{mask:birkhoff2}
\begin{align*}
a=&\left\{
\begin{bmatrix}
\frac{151}{3136}+\frac{5}{56}t &\frac{15}{64}\\[0.1em]
-\frac{1887}{307328}-\frac{321}{10976}t
-\frac{5}{147}t^2
&-\frac{359}{12544}-\frac{5}{56}t \end{bmatrix},
\begin{bmatrix}
\frac{1}{4}&0\\[0.1em]
-\frac{31}{1568}-\frac{1}{14}t &\frac{1}{16} \end{bmatrix},
\begin{bmatrix}
\frac{633}{1568}-\frac{5}{28}t
&-\frac{15}{32}\\[0.1em]
-\frac{31}{2401}-\frac{953}{5488}t
+\frac{10}{147}t^2 &-\frac{327}{6272}+\frac{5}{28}t \end{bmatrix},\right.\\
&\left.\;\;
\begin{bmatrix}
\frac{1}{4}&0\\[0.1em]
-\frac{31}{1568}-\frac{1}{14}t &\frac{1}{16} \end{bmatrix},
\begin{bmatrix}
\frac{151}{3136}+\frac{5}{56}t &\frac{15}{64}\\[0.1em]
-\frac{1887}{307328}-\frac{321}{10976}t-\frac{5}{147}t^2
&-\frac{359}{12544}-\frac{5}{56}t \end{bmatrix}
\right\}_{[-2,2]},
\end{align*}
where $t\in \R$ and $a$ has order $8$ sum rules with a matching filter $\vgu_a\in \lrs{0}{1}{r}$ satisfying
\[
\wh{\vgu_a}(\xi)=[1+(-\tfrac{3}{49}+\tfrac{8}{21}t)(i\xi)^2
+(\tfrac{11}{5880}-\tfrac{4}{63}t)(i\xi)^4
+\tfrac{t}{135}(i\xi)^6, (i\xi)^2-\tfrac{1}{6}(i\xi)^4+\tfrac{7}{360}(i\xi)^6]+\bo(|\xi|^8),\quad \xi\to 0.
\]
By calculation, $\sm_2(a)=5.5$ and its basis vector function $\phi=[\phi_1,\phi_2]^\tp$ is a spline vector function with symmetry and support $[-2,2]$ such that $\phi_1=\phi_1(-\cdot)$ and $\phi_2=\phi_2(-\cdot)$ and
{\small{\[
\phi|_{[0,1]}(x)=\begin{bmatrix}
-\frac{1}{7}x^7 + \frac{2}{3}x^6 - \frac{3}{2}x^5 + \frac{5}{3}x^4 - \frac{4}{3}x^2 + \frac{16}{21}\\
(\tfrac{-85}{4116}+\tfrac{8}{147}t)x^7
+(\tfrac{61}{294}-\tfrac{16}{63}t)x^6+
(-\tfrac{386}{735}+\tfrac{4}{7}t)x^5
+(\tfrac{64}{147}-\tfrac{40}{63}t) x^4+(-\tfrac{4}{49}+\tfrac{32}{63}t)x^2-(\tfrac{152}{5145}+\tfrac{128}{441}t)\end{bmatrix}
\]}}
and
\[
\phi|_{(1,2]}(x)=
\begin{bmatrix}
\frac{(x-2)^5(2x^2-8x+1)}{42}\\
\frac{(x-2)^5((555+1120t)x^2-(2220+4480t)x+792+560t)}{61740}
\end{bmatrix}.
\]
Hence, we have $\sm_\infty(\phi)=\sm_\infty(a)=5$.
By \cref{thm:converg},
the Birkhoff subdivision scheme of type $\ind$ with mask $a$ is convergent with limit functions in $\CH{4}$.
%See \cref{fig:ex1} for the graphs of the refinable vector functions $\phi$ with or without interpolation properties.
\end{example}

\begin{example}\label{ex2}
\normalfont
Let $d=1$, $\ind=\{0,1\}$ and $T=\{\frac{1}{2},\frac{1}{2}\}$.
We consider generalized Hermite subdivision schemes of type $\ind=\{0,1\}$ (more precisely, a dual (or face-based) Hermite subdivision scheme of degree $1$).
Using \cref{thm:lpm} with $\wh{\vgu_a}(\xi)=[e^{i\xi/2}, i\xi e^{i\xi/2}]+\bo(|\xi|^4)$ as $\xi \to 0$,
we find that all such symmetric Hermite masks $a\in \lrs{0}{2}{2}$ with $\fs(a)=[-1,2]$ and $\lpm(a)\ge 4$ are given by
\be \label{ex2:mask1}
a=\left\{
\begin{bmatrix}
\frac{5}{64} &\frac{9}{32}\\[0.1em]
-\frac{3}{128} &-\frac{5}{64} \end{bmatrix},
\begin{bmatrix}
\frac{27}{64}&\frac{9}{32}\\[0.1em]
-\frac{9}{128} &\frac{3}{64} \end{bmatrix},
\begin{bmatrix}
\frac{27}{64} &-\frac{9}{32}\\[0.1em]
\frac{9}{128} &\frac{3}{64} \end{bmatrix},
\begin{bmatrix}
\frac{5}{64} &-\frac{9}{32}\\[0.1em]
\frac{3}{128} &-\frac{5}{64} \end{bmatrix}
\right\}_{[-1,2]}.
\ee
By calculation, we have $\sm_2(a)\approx 3.33904$ and hence $\sm_\infty(a)\ge \sm_2(a)-0.5>2$. By \cref{thm:lpm,thm:converg},
the (dual) Hermite subdivision scheme of type $\ind$ with mask $a$ is convergent with limit functions in $\CH{2}$ and has the polynomial-interpolation property of order $4$ by interpolating all polynomials of degree less than $4$.

If we give up the polynomial-interpolation property,
then we find that a symmetric Hermite mask $a\in \lrs{0}{2}{2}$ with $\fs(a)=[-1,2]$ and $\sr(a)\ge 6$ is given by
\be \label{ex2:mask2}
a=\left\{
\begin{bmatrix}
\frac{13}{128} &\frac{15}{64}\\[0.1em]
-\frac{33}{1280} &-\frac{7}{128} \end{bmatrix},
\begin{bmatrix}
\frac{51}{128}&\frac{15}{64}\\[0.1em]
-\frac{63}{1280} &\frac{9}{128} \end{bmatrix},
\begin{bmatrix}
\frac{51}{128} &-\frac{15}{64}\\[0.1em]
\frac{63}{1280} &\frac{9}{128} \end{bmatrix},
\begin{bmatrix}
\frac{13}{128} &-\frac{15}{64}\\[0.1em]
\frac{33}{1280} &-\frac{7}{128} \end{bmatrix}
\right\}_{[-1,2]}.
\ee
The above mask $a$ has order $6$ sum rules with a matching filter $\vgu_a\in \lrs{0}{1}{r}$ given by
\[
\wh{\vgu_a}(\xi)=[1+\tfrac{1}{2}(i\xi)+\tfrac{1}{10} (i\xi)^2+\tfrac{1}{120}(i\xi)^3,
i\xi+\tfrac{1}{2}(i\xi)^2+\tfrac{1}{12}(i\xi)^3]+
\bo(|\xi|^6),\quad \xi \to 0.
\]
By calculation, $\sm_2(a)=4.5$ and its basis vector function $\phi=[\phi_1,\phi_2]^\tp$ is a spline vector function with symmetry and support $[-1,2]$ such that $\phi_1=\phi_1(1-\cdot)$ and $\phi_2=-\phi_2(1-\cdot)$ and
{\footnotesize{\[
\phi|_{[-1,0]}(x)=\begin{bmatrix}
\frac{(x+1)^4(2-3x)}{4}\\
\frac{(x+1)^4(11x-4)}{40}
\end{bmatrix},\quad
\phi|_{(0,1]}(x)=
\begin{bmatrix}
\frac{5}{4}x^4-\frac{5}{2}x^3+\frac{5}{4}x+\frac{1}{2}\\
\frac{1}{40}(2x-1)(19x^4-38x^3+6x^2+13x+4)
\end{bmatrix},\quad
\phi|_{(1,2]}(x)=
\begin{bmatrix}
\frac{(x-2)^4(3x-1)}{4}\\
\frac{(x-2)^4(11x-7)}{40}
\end{bmatrix}.
\]}}
Hence, we have $\sm_\infty(\phi)=\sm_\infty(a)=4$.
By \cref{thm:converg},
the (dual) Hermite subdivision scheme of type $\ind$ with mask $a$ is convergent with limit functions in $\CH{3}$.
%See \cref{fig:ex2} for the graphs of the refinable vector functions $\phi$.
\end{example}

For $j=1,2$, assume that a generalized Hermite subdivision scheme of type $\ind_j$ and translation $T_j$ with mask $a_j\in (l_0(\Z^{d_j}))^{r_j\times r_j}$
is convergent with limit functions in $\mathscr{C}^{m_j}(\R^{d_j})$. Then the generalized Hermite subdivision scheme of type $\ind_1\otimes \ind_2$ and translation $T_1\otimes T_2$ with mask $a_1\otimes a_2\in (l_0(\Z^{d_1+d_2}))^{(r_1+r_2)\times (r_1+r_2)}$ must converge with limit functions in $\mathscr{C}^m(\R^{d_1+d_2})$, where $m:=\min(m_1,m_2)$ and $\otimes$ stands for the tensor product. Moreover, if both generalized Hermite subdivision schemes have the interpolation property or linear-phase moments in \cref{sec:hsdint}, then so does their tensor product one. Consequently, examples of multivariate generalized Hermite subdivision schemes can be easily obtained from univariate ones through tensor product. In the following, we are only interested in non-tensor product ones. A few examples of bivariate interpolatory and non-interpolatory Hermite subdivision schemes have been reported in \cite{dhmm05,hm07,hyp04,hyx05,hy06}, while a few bivariate interpolating generalized Hermite subdivision schemes are given in \cite{hz09}. Because the presentation of high-dimensional examples is often quite messy, here we only provide four relatively simple bivariate examples.

\begin{example}\label{ex3}
\normalfont
Let $d=2$, $\ind=\{(0,0), (1,0), (0,1)\}$ and $T=\{(0,0), (0,0), (0,0)\}$.
We consider three families of bivariate generalized Hermite subdivision schemes of type $\ind$ (more precisely, standard Hermite subdivision schemes of degree $1$). Using \cref{thm:lpm} with $\wh{\vgu_a}(\xi_1, \xi_2)=[1,i\xi_1, i\xi_2]+\bo(\|\xi\|^4)$ as $\xi=(\xi_1,\xi_2)^\tp \to 0$, we find that all bivariate $D_6$-symmetric Hermite masks $a\in (l_0(\Z^2))^{3\times 3}$ of type $\ind$ with $\fs(a)=[-2,2]^2$ and $\lpm(a)\ge 4$ are given by $a(0,0)=\mbox{diag}(\frac{1}{4}-12 t_2, \frac{1}{8}+6t_3+6t_4, \frac{1}{8}+6t_3+6t_4)$,
{\footnotesize{\[
a(1,0)=\begin{bmatrix}
-4t_1 &-\frac{3}{16} &\frac{3}{32}\\[0.1em]
-t_1 &-\frac{1}{32} &\frac{1}{64}\\[0.1em]
0 &0 &0
\end{bmatrix},\;
a(2,0)=\begin{bmatrix}
2t_2 &2t_3+4t_4 &-t_3-2t_4\\[0.1em]
t_2 &t_3+t_4 &-t_3\\[0.1em]
0 &0 &-t_3+t_4\end{bmatrix},\;
a(2,1)=\begin{bmatrix}
\frac{1}{8}+4t_1 &-\frac{3}{32} &0\\[0.1em]
\frac{1}{16}+2t_1 &-\frac{1}{32} &0\\[0.1em]
\frac{1}{32}+t_1 &-\frac{1}{64} &0
\end{bmatrix},
\]}}
where $t_1, t_2, t_3, t_4\in \R$
and all other nonzero $a(k), k\in \Z^2$ are determined through symmetry in \eqref{sym:mask}, i.e., $a(Ek)=S(E, \ind) a(k)S(E,\ind)^{-1}$ for all $k\in \Z^2$ and $E\in D_6$. Because $-\frac{1}{8}$ is an eigenvalue of the transition operator $\mathcal{T}_a$ in \eqref{tz},
by \cite[Theorem~6]{hm07}, we conclude that the refinable vector function $\phi$ with mask $a$ cannot be a spline vector for any choice of $t_1,\ldots,t_4\in \R$. For $t_1=-\frac{15}{512}, t_2=\frac{1}{512}, t_3=-\frac{1}{128}$ and $t_4=\frac{1}{256}$, we have $\sm_2(a)\approx 3.13452$ and
hence, $\sm_\infty(a)\ge \sm_2(a)-1>2$. By \cref{thm:lpm,thm:converg}, the Hermite subdivision scheme of type $\ind$ with mask $a$ is convergent with limit functions in $\mathscr{C}^2(\R^2)$ and interpolates all the bivariate polynomials of (total) degree less than $4$.

Note that the identity mapping $\theta: \{1,2,3\}\rightarrow \{1,2,3\}$ satisfies \eqref{theta} with $T=\{(0,0), (0,0),(0,0)\}$. We find that all $D_6$-symmetric interpolatory Hermite masks $a\in (l_0(\Z^2))^{3\times 3}$ of type $\ind$ and translation $T$ with $\fs(a)\subseteq
[-2,2]^2$ and $\sr(a)\ge 4$ must be given by the above mask $a$ with $t_2=t_3=t_4=0$. In particular, for $t_1=-\frac{1}{32}$, we have $\sm_2(a)\approx 2.71094$ and hence $\sm_\infty(a)\ge \sm_2(a)-1>1$.
By \cref{thm:int,thm:converg}, the interpolatory Hermite subdivision scheme of type $\ind$ and translation $T$ with mask $a$ is convergent with limit functions in $\mathscr{C}^1(\R^2)$ and its basis vector function $\phi$ is a refinable Hermite interpolant of type $\ind$ and translation $T$ satisfying \eqref{int:phi}.

If we give up the linear-phase moments for the polynomial-interpolation property,
then we find three families of $D_6$-symmetric Hermite masks $a\in (l_0(\Z^2))^{3\times 3}$ of type $\ind$ with $\fs(a)\subseteq [-2,2]^2$ and $\sr(a)\ge 5$, one of these three families is given by $a(0,0)=\mbox{diag}(\frac{47}{128}-\frac{27}{8} t_3, \frac{1}{8}-6t_1-6t_2, \frac{1}{8}-6t_1-6t_2)$,
{\small{\begin{align*}
&a(1,0)=\begin{bmatrix}
\frac{21}{128}-\frac{9}{8}t_3
&-\frac{3}{8}+\frac{9}{2} t_3
&\frac{3}{16}-\frac{9}{4} t_3\\[0.1em]
\frac{3}{64}-\frac{9}{16} t_3
&-\frac{7}{64}+\frac{15}{8} t_3
&\frac{1}{32}-\frac{3}{8}t_3\\[0.1em]
0 &0 &-\frac{3}{64}+\frac{9}{8} t_3
\end{bmatrix},\quad
a(2,1)=\begin{bmatrix}
-\frac{5}{128}+\frac{9}{8}t_3
&-\frac{3}{16}+\frac{9}{4} t_3
&0\\[0.1em]
-\frac{1}{64}+\frac{3}{8}t_3
&-\frac{5}{64}+\frac{9}{8} t_3
&0\\[0.1em]
-\frac{1}{128}+\frac{3}{16} t_3
&-\frac{1}{32}+\frac{3}{8} t_3 &-\frac{1}{64}+\frac{3}{8} t_3,
\end{bmatrix},\\
&a(2,0)=\begin{bmatrix}
-\frac{5}{256}+\frac{9}{16}t_3
&-2t_1 &t_1\\[0.1em]
-\frac{1}{128}+\frac{3}{16}t_3
&-\frac{5}{4}t_1-t_2+\frac{3}{16}t_3
&\frac{5}{4} t_1+2t_2-\frac{3}{16} t_3\\[0.1em]
0 &0 &\frac{5}{4}t_1+3t_2-\frac{3}{16}t_3\end{bmatrix},
\end{align*}
}}%
where $t_1, t_2, t_3\in \R$
and all other nonzero $a(k), k\in \Z^2$ are determined through symmetry in \eqref{sym:mask}, i.e., $a(Ek)=S(E, \ind) a(k)S(E,\ind)^{-1}$ for all $k\in \Z^2$ and $E\in D_6$. Moreover, $a$ has order $5$ sum rules with the matching filter $\vgu_a\in (l_0(\Z^2))^{1\times 3}$ satisfying
\begin{align*}
\wh{\vgu_a}(\xi_1,\xi_2)=
&\left[1+(2t_3-\tfrac{1}{24})
(\xi_1^2+\xi_1 \xi_2+\xi_2^2)
+(\tfrac{16}{5} t_3^2-\tfrac{5}{15}t_3-\tfrac{1}{720})
(\xi_1^4+2\xi_1^3 \xi_2+3\xi_1^2\xi_2^2+2\xi_1 \xi_2^3+\xi_2^4),\right.\\
&\left.
i\xi_1
(1+2t_3(\xi_1^2+\xi_1\xi_2+\xi_2^2)),
i\xi_2
(1+2t_3(\xi_1^2+\xi_1\xi_2+\xi_2^2))\right]+\bo(\|\xi\|^5),
\end{align*}
as $\xi=(\xi_1,\xi_2)^\tp \to 0$.
We notice that both $9t_3-\frac{1}{2}$ and
$\frac{3}{2}t_3-\frac{1}{16}$ are eigenvalues of $\mathcal{T}_a$ in \eqref{tz}.
If $9t_3-\frac{1}{2}=2^{-n}$ for some $n\in \NN\cup\{\infty\}$, then $t_3=\frac{1}{18}+\frac{2^{-n}}{9}$ and $\frac{3}{2}t_3-\frac{1}{16}=\frac{1}{48}+\frac{2^{-n}}{6}
\not\in \{2^{-j}\}_{j\in \NN\cup\{\infty\}}$.
By \cite[Theorem~6]{hm07}, we conclude that the refinable vector function $\phi$ with mask $a$ cannot be a spline vector for any choice of $t_1,t_2,t_3\in \R$.
For $t_1=\frac{5}{256}, t_2=-\frac{1}{256}$ and $t_3=\frac{29}{512}$, we have $\sm_2(a)\approx 4.81514$ and
hence, $\sm_\infty(a)\ge \sm_2(a)-1>3$. By \cref{thm:converg}, the (standard) Hermite subdivision scheme of type $\ind$ and translation $T$ with mask $a$ is convergent with limit functions in $\mathscr{C}^3(\R^2)$.
\end{example}

\begin{example}\label{ex4}
\normalfont
Let $d=2$, $\ind=\{(0,0), (1,0), (0,1)\}$ and $T=\{(\frac{1}{2},\frac{1}{2}), (\frac{1}{2},\frac{1}{2}),(\frac{1}{2},\frac{1}{2})\}$.
We consider two families of bivariate generalized Hermite subdivision schemes of type $\ind$ (i.e., dual Hermite subdivision schemes of degree $1$). Using \cref{thm:lpm} with $\wh{\vgu_a}(\xi_1, \xi_2)=[e^{i(\xi_1+\xi_2)/2},i\xi_1 e^{i(\xi_1+\xi_2)/2}, i\xi_2e^{i(\xi_1+\xi_2)/2}]+\bo(\|\xi\|^4)$ as $\xi=(\xi_1,\xi_2)^\tp \to 0$, we find that all bivariate $D_4$-symmetric Hermite masks $a\in (l_0(\Z^2))^{3\times 3}$ of type $\ind$ with $\fs(a)=[-1,2]^2$ and $\lpm(a)\ge 4$ are given by
{\tiny{
\begin{align*}
&a(1,1)=\begin{bmatrix}
\frac{21}{128}-4t_3 &-\frac{9}{64}+t_2 &-\frac{9}{64}+t_2\\[0.1em]
\frac{3}{128}-t_3 &-\frac{3}{256}+t_1+\frac{1}{2} t_2 &-t_1\\[0.1em]
\frac{3}{128}-t_3 &-t_1 &-\frac{3}{256}+t_1+\frac{1}{2} t_2
\end{bmatrix},\quad
a(2,1)=\begin{bmatrix}
\frac{3}{64}+4t_3 &-\frac{9}{64}+t_2 &-t_2\\[0.1em]
\frac{3}{256}+t_3 &-\frac{15}{256}+t_1+\frac{1}{2}t_2 &-\frac{3}{128}+t_1\\[0.1em]
\frac{3}{256}+t_3 &-t_1 &\frac{9}{256}-t_1-\frac{1}{2}t_2\end{bmatrix},\quad\\
&a(2,2)=\begin{bmatrix}
-\frac{1}{128}-4t_3 &-t_2 &-t_2\\[0.1em]
-t_3 &\frac{5}{256}-t_1-\frac{1}{2}t_2 &-\frac{3}{128}+t_1\\[0.1em]
-t_3 &-\frac{3}{128}+t_1 &\frac{5}{256}-t_1-\frac{1}{2} t_2
\end{bmatrix},
\end{align*}
}}
and all other nonzero $a(k), k\in \Z^2$ are determined through symmetry in \eqref{sym:mask}, i.e.,
\be \label{D4}
a(E(k-(\tfrac{1}{2},\tfrac{1}{2}))+(\tfrac{1}{2},\tfrac{1}{2}))
=S(E, \ind) a(k)S(E,\ind)^{-1},\qquad \forall\, k\in \Z^2, E\in D_4.
\ee
Because the imaginary numbers
$\pm \frac{\sqrt{2}}{32}i$ are eigenvalues of $\mathcal{T}_a$ in \eqref{tz},
by \cite[Theorem~6]{hm07}, we conclude that the refinable vector function $\phi$ with mask $a$ cannot be a spline vector for any choice of $t_1,t_2,t_3\in \R$.
For $t_1=\frac{1}{64}, t_2=\frac{5}{128}$ and $t_3=0$, $\sm_2(a)\approx 3.33904$ and so
$\sm_\infty(a)\ge \sm_2(a)-1>2$. By \cref{thm:lpm,thm:converg}, the dual Hermite subdivision scheme of degree $1$ with mask $a$ is convergent with limit functions in $\mathscr{C}^2(\R^2)$ and interpolates all the bivariate polynomials of degree less than $4$.

If we give up the polynomial-interpolation property, then we find two families of $D_4$-symmetric Hermite masks $a\in (l_0(\Z^2))^{3\times 3}$ of type $\ind$ with $\fs(a)\subseteq [-1,2]^2$ and $\sr(a)\ge 5$, one of which is
{\tiny{\begin{align*}
&a(1,1)=\begin{bmatrix}
\frac{5}{32}+4t_1
&-\frac{5}{128}-4t_2
&-\frac{5}{128}-4 t_2\\[0.1em]
\frac{1}{64}+t_1
&\frac{11}{256}-t_2
&-t_2\\[0.1em]
\frac{1}{64}+t_1 &-t_2 &\frac{11}{256}- t_2
\end{bmatrix},\quad
a(2,1)=\begin{bmatrix}
\frac{1}{32}-4t_1
&-\frac{5}{128}-4t_2
&-\frac{7}{128}+4t_2\\[0.1em]
\frac{1}{128}-t_1
&-\frac{1}{256}- t_2
&-\frac{1}{64}+t_2\\[0.1em]
-t_1
&-t_2 &\frac{1}{256}+t_2,
\end{bmatrix},\\
&a(2,2)=\begin{bmatrix}
\frac{1}{32}+4t_1
&-\frac{7}{128}+4t_2 &-\frac{7}{128}+4t_2\\[0.1em]
\frac{1}{128}+t_1
&-\frac{3}{256}+t_2
&-\frac{1}{64}+ t_2\\[0.1em]
 \frac{1}{128}+t_1& -\frac{1}{64}+t_2 &-\frac{3}{256}+t_2\end{bmatrix},
\end{align*}
}}%
and all other nonzero $a(k), k\in \Z^2$ are determined through symmetry in \eqref{D4}.
Moreover, the mask $a$ has order $5$ sum rules with the matching filter $\vgu_a\in (l_0(\Z^2))^{1\times 3}$ satisfying
{\small{\begin{align*}
\wh{\vgu_a}(\xi_1,\xi_2)=
&\left[
1+\tfrac{i}{2}(\xi_1+\xi_2)
-\tfrac{1}{12}(\xi_1^2+3\xi_1\xi_2+\xi_2^2)
-\tfrac{i}{24}\xi_1\xi_2(\xi_1+\xi_2),
i\xi_1 (1+\tfrac{i}{2}(\xi_1+\xi_2)
-\tfrac{1}{12}\xi_1^2
-\tfrac{1}{4}\xi_1\xi_2\right.\\
&\qquad \left.+\tfrac{i}{24} \xi_2^3
-\tfrac{i}{24} \xi_1^2\xi_2),
i\xi_2
(1+\tfrac{i}{2}(\xi_1+\xi_2)
-\tfrac{1}{12}\xi_2^2
-\tfrac{1}{4}\xi_1\xi_2+\tfrac{i}{24} \xi_1^3
-\tfrac{i}{24} \xi_1\xi_2^2
\right]+\bo(\|\xi\|^5),
\end{align*}
}}%
as $\xi=(\xi_1,\xi_2)^\tp \to 0$.
We notice that $-32t_1-\frac{3}{32}, 16t_1+8t_2+\frac{1}{64}$ and $64t_1-32t_2+\frac{9}{16}$ are eigenvalues of $\mathcal{T}_a$ in \eqref{tz}.
If $-32t_1-\frac{3}{32}=2^{-j}$ and $16t_1+8t_2+\frac{1}{64}=2^{-k}$ for some $j,k\in \NN\cup{\infty}$, then $t_1=-\frac{3}{1024}-\frac{2^{-j}}{32}, t_2=\frac{1}{256}+\frac{2^{-j}}{16}+\frac{2^{-k}}{8}$ and consequently, $64t_1-32t_2+\frac{9}{16}=\frac{1}{4}-2^{2-j}-2^{2-k}$, which belongs to $\{2^{-n}\}_{n\in \NN\cup{\infty}}$ if and only if $j,k\in \{5,6\}$, or $(j,k)=(4,\infty), (\infty,4), (5,\infty), (5,\infty), (\infty,\infty)$.
For all these cases, we have $\sm_2(a)=3$.
By \cite[Theorem~6]{hm07}, we conclude that these refinable vector functions $\phi$ with mask $a$ cannot be a spline vector for any choice of $t_1,t_2\in \R$. For $t_1=0$ and $t_2=\frac{1}{64}$, we have $\sm_2(a)\approx 3.0$ and hence, $\sm_\infty(a)\ge \sm_2(a)-1>1$. By \cref{thm:converg}, the dual (or face-based) Hermite subdivision scheme of degree $1$ with mask $a$ is convergent with limit functions in $\mathscr{C}^1(\R^2)$. This agrees with the example in \cite[Section~3.2]{hy06} and is not surprising because under $D_4$ symmetry and $\ind=\{(0,0), (1,0), (0,1)\}=\ind_1$, the sufficient condition in \cite[Theorem~2.3]{hy06} agrees with the necessary and sufficient condition in \cref{thm:hsd:mask}.
\end{example}

\begin{example}\label{ex5}
\normalfont
Let $d=2$, $\ind=\{(0,0), (1,1)\}$ and $T=\{(\frac{1}{2},\frac{1}{2}), (\frac{1}{2},\frac{1}{2})\}$.
We consider bivariate generalized Hermite subdivision schemes of type $\ind$ (more precisely, Birkhoff subdivision schemes).
We find that all bivariate $D_4$-symmetric generalized Hermite masks $a\in (l_0(\Z^2))^{2\times 2}$ of type $\ind$ with $\fs(a)=[-2,3]^2$ and $\sr(a)\ge 5$ are given by three families.
To reduce too many free parameters in these families, we artificially set $a(3,3)=a(3,2)=0$. Two of the three families are given by
{\tiny{
\begin{align*}
&a_1(1,1)=\begin{bmatrix}
\frac{17}{256}-6t_3 &\frac{3}{32}-6t_2\\[0.1em]
t_1 &-\frac{3}{256}-\frac{5}{2} t_2
\end{bmatrix},\quad
a_1(2,1)=\begin{bmatrix}
\frac{19}{256}+6t_3 &\frac{3}{32}-6t_2\\[0.1em]
-\frac{1}{512}+t_1 &-\frac{3}{256}+\frac{5}{2}t_2 \end{bmatrix},\quad\\
&a_1(2,2)=\begin{bmatrix}
\frac{1}{256}-6t_3 &\frac{3}{32}-6t_2 \\[0.1em]
\frac{1}{128}+t_1+4t_3&-\frac{5}{256}+
\frac{3}{2}t_2
\end{bmatrix},\quad
a_1(3,1)=\begin{bmatrix}
\frac{1}{64} &0 \\[0.1em]
t_3&t_2\end{bmatrix},
\end{align*}
}}
and
{\tiny{\begin{align*}
&a_2(1,1)=\begin{bmatrix}
\frac{21}{256}-2t_2-2t_3 &\frac{1}{64}-4t_1\\[0.1em]
t_3 & t_1
\end{bmatrix},\quad
a_2(2,1)=\begin{bmatrix}
\frac{15}{256}+2t_2+2t_3 &\frac{1}{64}-4t_1\\[0.1em]
-t_2 &\frac{1}{64}+t_1 \end{bmatrix},\quad\\
&a_2(2,2)=\begin{bmatrix}
\frac{5}{256}-2t_2-2t_3 &\frac{1}{64}-4t_1 \\[0.1em]
t_3&t_1
\end{bmatrix},\quad
a_2(3,1)=\begin{bmatrix}
\frac{1}{64} &0 \\[0.1em]
-\frac{1}{512}&\frac{1}{128}\end{bmatrix},
\end{align*}
}}%
where $t_1, t_2, t_3\in \R$,
and all other nonzero $a_1(k)$ and $a_2(k), k\in \Z^2$ are determined through symmetry in \eqref{D4}. The mask $a_1$ has order $5$ sum rules with a matching filter $\vgu_{a_1}$ satisfying
\begin{align*}
\wh{\vgu_{a_1}}(\xi_1,\xi_2)=&\left[
1+\tfrac{i}{2}(\xi_1+\xi_2)+
\tfrac{1}{12}(\xi_1^2-3\xi_1 \xi_2+\xi_2^2)+\tfrac{i}{24}(2\xi_1^3+\xi_2^1 \xi_2+\xi_1\xi_2^2+2\xi_2^2)
-\tfrac{1}{24}(\xi_1^3 \xi_2+\xi_1\xi_2^3)
\right.\\
&\left.
(\tfrac{1}{144}-\tfrac{64}{15}t_1 -\tfrac{32}{3}t_3)\xi_1^2\xi_2^2,
(-\xi_1\xi_2)(
1+\tfrac{i}{2}(\xi_1+\xi_2)+
\tfrac{1}{3}\xi_1^2-\tfrac{1}{4}
\xi_1\xi_2+\tfrac{1}{3}\xi_2^2)\right]+\bo(\|\xi\|^5)
\end{align*}
as $\xi=(\xi_1,\xi_2)^\tp\to 0$.
The transition operator $\mathcal{T}_a$ in \eqref{tz} has eigenvalues $4t_2, \frac{1}{32}-16t_2$ and a few complicated ones. By calculation we check that the eigenvalue condition in \cite[Theorem~6]{hm07} fails and the refinable vector function $\phi$ with mask $a_1$ cannot be a spline vector for any choice of $t_1,t_2,t_3\in \R$.
For $t_1=\frac{1}{512}, t_2=\frac{1}{128}$ and $t_3=-\frac{1}{256}$, $\sm_2(a_1)\approx 3.41080$ and hence
$\sm_\infty(a_1)\ge \sm_2(a_1)-1>2$.
The mask $a_2$ has order $5$ sum rules with a matching filter $\vgu_{a_2}$ satisfying
\begin{align*}
\wh{\vgu_{a_2}}(\xi_1,\xi_2)=&\left[
1+\tfrac{i}{2}(\xi_1+\xi_2)+
\tfrac{1}{12}(\xi_1^2-3\xi_1 \xi_2+\xi_2^2)+\tfrac{i}{24}(2\xi_1^3+
\xi_1^2 \xi_2+\xi_1\xi_2^2+2\xi_2^2)
-\tfrac{1}{24}(\xi_1^3 \xi_2+\xi_1\xi_2^3)
\right.\\
&\left.
(\tfrac{17}{720}+\tfrac{32}{15}t_2 -\tfrac{32}{15}t_3)\xi_1^2\xi_2^2,
(-\xi_1\xi_2)(
1+\tfrac{i}{2}(\xi_1+\xi_2)+
\tfrac{1}{6}\xi_1^2-\tfrac{1}{4}
\xi_1\xi_2+\tfrac{1}{6}\xi_2^2)\right]+\bo(\|\xi\|^5)
\end{align*}
as $\xi=(\xi_1,\xi_2)^\tp\to 0$.
Because $-\frac{3}{32}$ is an eigenvalue of $\mathcal{T}_a$ in \eqref{tz},
by \cite[Theorem~6]{hm07}, we conclude that the refinable vector function $\phi$ with mask $a_2$ cannot be a spline vector for any choice of $t_1,t_2,t_3\in \R$. For $t_1=-\frac{1}{128}, t_2=\frac{1}{256}$ and $t_3=-\frac{1}{256}$, $\sm_2(a_2)\approx 3.59632$ and hence
$\sm_\infty(a_2)\ge \sm_2(a_2)-1>2$.
By \cref{thm:converg}, the bivariate Birkhoff subdivision schemes of type $\ind$ with both masks $a_1$ and $a_2$ are convergent with limit functions in $\mathscr{C}^2(\R^2)$.
\end{example}

\begin{example}\label{ex6:D4}
\normalfont
Let $d=2$, $\ind=\{(0,0), (0,0)\}$ and $T=\{(0,0), (\frac{1}{2},\frac{1}{2})\}$.
We consider bivariate generalized Hermite subdivision schemes of type $\ind$ (more precisely, Lagrange subdivision schemes).
We consider all bivariate $D_4$-symmetric generalized Hermite masks $a\in (l_0(\Z^2))^{2\times 2}$ with $\fs(a)=[-2,2]^2$ under the extra condition $a(-2,2)=0$ for reducing freedoms. For $\sr(a)\ge 4$, we obtain
four families such that one of them with the smallest number of free parameters is given by
{\tiny{
\begin{align*}
&a(2,-2)=\begin{bmatrix}
0 &0\\[0.1em]
0 &0\end{bmatrix},\quad
a(2,-1)=\begin{bmatrix}
\frac{1}{128} &0\\[0.1em]
0 &t_2-\frac{1}{64}\end{bmatrix},\quad
a(2,0)=\begin{bmatrix}
\frac{1}{64} &0\\[0.1em]
t_1 &t_2 \end{bmatrix},\quad
a(2,1)=\begin{bmatrix}
\frac{1}{128} &0 \\[0.1em]
\frac{1}{32} &t_2\end{bmatrix},\quad
a(2,2)=\begin{bmatrix}
0 &0 \\[0.1em]
t_1 &t_2-\frac{1}{64}\end{bmatrix},\\
&a(1,-2)=\begin{bmatrix} \frac{1}{128} &\frac{1}{128} \\
0 &0
\end{bmatrix},\quad
a(1,-1)=\begin{bmatrix} \frac{3}{64} &\frac{3}{128}\\[0.1em]
0 &t_2\end{bmatrix},\quad
a(1,0)=\begin{bmatrix} \frac{5}{64} &\frac{3}{128}\\[0.1em]
\frac{1}{32} &\frac{3}{64}+t_2\end{bmatrix},\quad
a(1,1)=\begin{bmatrix} \frac{3}{64} &\frac{1}{128}\\[0.1em]
\frac{1}{16} &\frac{3}{64}+t_2\end{bmatrix},\quad
a(1,2)=\begin{bmatrix} \frac{1}{128} &0\\[0.1em]
\frac{1}{32} &t_2 \end{bmatrix},\\
&a(0,-2)=\begin{bmatrix} \frac{1}{64} &\frac{3}{128}\\[0.1em]
0 &0\end{bmatrix},\quad
a(0,-1)=\begin{bmatrix} \frac{5}{64} &\frac{21}{128}-4t_2\\[0.1em]
0 &t_2\end{bmatrix},\quad
a(0,0)=\begin{bmatrix} \frac{3}{16}-4t_1 &\frac{21}{128}-4t_2\\[0.1em]
t_1 &\frac{3}{64}+t_2\end{bmatrix},\quad
a(0,1)=\begin{bmatrix} \frac{5}{64} &\frac{3}{128}\\[0.1em]
\frac{1}{32} &\frac{3}{64}+t_2\end{bmatrix},\quad
a(0,2)=\begin{bmatrix} \frac{1}{64} &0\\[0.1em]
t_1 &t_2\end{bmatrix},\\
&a(-1,-2)=\begin{bmatrix} \frac{1}{128} &\frac{3}{128}\\[0.1em]
0 &0\end{bmatrix},\quad
a(-1,-1)=\begin{bmatrix} \frac{3}{64} &\frac{21}{128}-4t_2\\[0.1em]
0 &t_2-\frac{1}{64}\end{bmatrix},\quad
a(-1,0)=\begin{bmatrix} \frac{5}{64} &\frac{21}{128}-4t_2\\[0.1em]
0 &t_2\end{bmatrix},\quad
a(-1,1)=\begin{bmatrix} \frac{3}{64} &\frac{3}{128}\\[0.1em]
0 &t_2\end{bmatrix},\quad
a(-1,2)=\begin{bmatrix} \frac{1}{128} &0\\[0.1em]
0 &t_2-\frac{1}{64} \end{bmatrix},\\
&a(-2,-2)=\begin{bmatrix} 0 &\frac{1}{128}\\[0.1em]
0 &0\end{bmatrix},\quad
a(-2,-1)=\begin{bmatrix} \frac{1}{128} &\frac{3}{128}\\[0.1em]
0 &0\end{bmatrix},\quad
a(-2,0)=\begin{bmatrix}\frac{1}{64} &\frac{3}{128}\\[0.1em]
0 &0 \end{bmatrix},\quad
a(-2,1)=\begin{bmatrix} \frac{1}{128} &\frac{1}{128}\\[0.1em]
0 &0\end{bmatrix},\quad
a(-2,2)=\begin{bmatrix} 0 &0\\ 0 &0 \end{bmatrix}.
\end{align*}
}}
for all $t_1,t_2\in \R$. The mask $a$ has order $4$ sum rules with a matching filter $\vgu_a\in (l_0(\Z^2))^{1\times 2}$ given by
\[
\wh{\vgu_a}(\xi_1,\xi_2)=[1+\tfrac{1}{12}(\xi_1^2+\xi_2^2),
1+\tfrac{i}{2}(\xi_1+\xi_2)+\tfrac{1}{12}(\xi_1^2-3\xi_1\xi_2+\xi_2^2)
+\tfrac{i}{24}(2\xi_1^3+\xi_1^2 \xi_2+\xi_1\xi_2^2+2\xi_2^3)]+\bo(\|\xi\|^4)
\]
as $\xi\to 0$. The necessary condition on special eigenvalues of $\mathcal{T}_a$ in \cite[Theorem~6]{hm07} is satisfied if and only if either $t_1=\frac{1}{64}, t_2=\frac{3}{128}$ with $\sm_2(a)=2$ or $t_1=t_2=\frac{1}{64}$ with $\sm_2(a)=3$. But further investigation indicates that their basis vector functions are not vector splines. For $t_1=\frac{1}{64}$ and $t_2=\frac{1}{128}$, $\sm_2(a)=4$ and hence $\sm_\infty(a)\ge \sm_2(a)-1\ge 3$. By \cref{thm:converg}, the generalized Hermite subdivision scheme of type $\ind$ with $t_1=\frac{1}{64}$ and $t_2=\frac{1}{128}$ is convergent with limit functions in $\mathscr{C}^2(\R^2)$.

Using \cref{thm:lpm} with $\wh{\vgu_{a_1}}(\xi_1,\xi_2)=[1, e^{i(\xi_1+\xi_2)/2}]+\bo(\|\xi\|^4)$ as $\xi\to 0$ for the polynomial-interpolation property, we find that all masks $a_1\in (l_0(\Z^2))^{2\times 2}$ with $\fs(a_1)=[-2,2]^2$, $a_1(-2,2)=0$, $\sr(a_1)\ge 4$ and $\lpm(a_1)\ge 4$ are given by
{\tiny{
\begin{align*}
&a_1(2,-2)=\begin{bmatrix}
0 &0\\[0.1em]
0 &0\end{bmatrix},\quad
a_1(2,-1)=\begin{bmatrix}
-\frac{1}{64} &0\\[0.1em]
0 &-\frac{1}{64}\end{bmatrix},\quad
a_1(2,0)=\begin{bmatrix}
t &0\\[0.1em]
-2t &0 \end{bmatrix},\quad
a_1(2,1)=\begin{bmatrix}
-\frac{1}{64} &0 \\[0.1em]
\frac{1}{8} &0 \end{bmatrix},\quad
a_1(2,2)=\begin{bmatrix}
0 &0 \\[0.1em]
-2t &-\frac{1}{64}\end{bmatrix},\\
&a_1(1,-2)=\begin{bmatrix} -\frac{1}{64} &-\frac{1}{64} \\
0 &0
\end{bmatrix},\quad
a_1(1,-1)=\begin{bmatrix} 0 &0\\[0.1em]
0 &0\end{bmatrix},\quad
a_1(1,0)=\begin{bmatrix} \frac{1}{32} &0\\[0.1em]
\frac{1}{8} &\frac{9}{64}\end{bmatrix},\quad
a_1(1,1)=\begin{bmatrix} 0 &-\frac{1}{64}\\[0.1em]
\frac{1}{4} &\frac{9}{64}\end{bmatrix},\quad
a_1(1,2)=\begin{bmatrix} -\frac{1}{64} &0\\[0.1em]
\frac{1}{8} &0 \end{bmatrix},\\
&a_1(0,-2)=\begin{bmatrix} t &0\\[0.1em]
0 &0\end{bmatrix},\quad
a_1(0,-1)=\begin{bmatrix} \frac{1}{32} &\frac{9}{64}\\[0.1em]
0 &0\end{bmatrix},\quad
a_1(0,0)=\begin{bmatrix} \frac{1}{4}-4t &\frac{9}{64}\\[0.1em]
-2t &\frac{9}{64}\end{bmatrix},\quad
a_1(0,1)=\begin{bmatrix} \frac{1}{32} &0\\[0.1em]
\frac{1}{8} &\frac{9}{64}\end{bmatrix},\quad
a_1(0,2)=\begin{bmatrix} t &0\\[0.1em]
-2t &0\end{bmatrix},\\
&a_1(-1,-2)=\begin{bmatrix} -\frac{1}{64} &0\\[0.1em]
0 &0\end{bmatrix},\quad
a_1(-1,-1)=\begin{bmatrix}0 &\frac{9}{64}\\[0.1em]
0 &-\frac{1}{64}\end{bmatrix},\quad
a_1(-1,0)=\begin{bmatrix} \frac{1}{32} &\frac{9}{64}\\[0.1em]
0 &0\end{bmatrix},\quad
a_1(-1,1)=\begin{bmatrix} 0 &0\\[0.1em]
0 &0\end{bmatrix},\quad
a_1(-1,2)=\begin{bmatrix} -\frac{1}{64} &0\\[0.1em]
0 &-\frac{1}{64} \end{bmatrix},\\
&a_1(-2,-2)=\begin{bmatrix} 0 &-\frac{1}{64}\\[0.1em]
0 &0\end{bmatrix},\quad
a_1(-2,-1)=\begin{bmatrix} -\frac{1}{64} &0\\[0.1em]
0 &0\end{bmatrix},\quad
a_1(-2,0)=\begin{bmatrix}t &0 \\[0.1em]
0 &0 \end{bmatrix},\quad
a_1(-2,1)=\begin{bmatrix} -\frac{1}{64} &-\frac{1}{64}\\[0.1em]
0 &0\end{bmatrix},\quad
a_1(-2,2)=\begin{bmatrix} 0 &0\\ 0 &0 \end{bmatrix},
\end{align*}
}}
where $t\in \R$.
Because $-\frac{1}{16}$ is an eigenvalue of the transition operator $\mathcal{T}_{a_1}$ in \eqref{tz},
by \cite[Theorem~6]{hm07}, we conclude that the refinable vector function $\phi$ with mask $a_1$ cannot be a spline vector for any choice of $t\in \R$. For $t=\frac{1}{64}$, we have $\sm_2(a_1)\approx 2.26564$ and hence $\sm_\infty(a_1)\ge \sm_2(a_1)-1>1$. Consequently,
the generalized Hermite subdivision scheme of type $\ind$ with mask $a_1$ and $t=\frac{1}{64}$ is convergent with limit functions in $\mathscr{C}^1(\R^2)$ and interpolates all the bivariate polynomials of degree less than $4$.

Using the mapping $\theta:\{1,2\}\rightarrow \{1,2\}$ with $\theta(1)=\theta(2)=1$ which satisfies \eqref{theta}, we find that such masks $a_1$ are interpolatory of type $\ind$ and translation $T$ if and only if $t=0$.
For $t=0$, $\sm_2(a_1)\approx 1.91361$ and by \cref{thm:int}, its basis vector function $\phi\in (\mathscr{C}^0(\R^2))^2$ is a generalized Hermite interpolant of type $\ind$ and translation $T$.
If we consider $\fs(a_2)=[-1,1]^2$, then
there is a unique interpolatory generalized Hermite mask $a_2\in (l_0(\Z^2))^{2\times 2}$ of type $\ind$ and translation $T$ with $\sr(a_2)=2$ and $\sm_2(a_2)=1.5$, whose nonzero elements are given by
{\tiny{\[
a_2(1,0)=a_2(0,1)=\begin{bmatrix}
\frac{1}{8} &0\\[0.1em]
0 &\frac{1}{8}\end{bmatrix},\quad
a_2(0,-1)=a_2(-1,0)=\begin{bmatrix}
\frac{1}{8} &\frac{1}{8}\\[0.1em]
0 &0\end{bmatrix},\quad
a_2(1,1)=\begin{bmatrix}
0 &0\\[0.1em]
\frac{1}{4} &\frac{1}{8}\end{bmatrix},\quad
a_2(0,0)=\begin{bmatrix}
\frac{1}{4} &\frac{1}{8}\\[0.1em]
0 &\frac{1}{8}\end{bmatrix},\quad
%&a_2(0,1)=\begin{bmatrix}
%\frac{1}{8} &0\\[0.1em]
%0 &\frac{1}{8}\end{bmatrix},\quad
a_2(-1,-1)=\begin{bmatrix}
0 &\frac{1}{8}\\[0.1em]
0 &0\end{bmatrix}.
%a_2(-1,0)=\begin{bmatrix}
%\frac{1}{8} &\frac{1}{8}\\[0.1em]
%0 &0\end{bmatrix},\quad
\]
}}
By \cref{thm:int}, its basis vector function $\phi=[\phi_1,\phi_2]^\tp \in (\mathscr{C}^0(\R^2))^2$ is a generalized Hermite interpolant of type $\ind$ and translation $T$, where $\phi_1,\phi_2$ are piecewise linear functions given by $\phi_1(0,0)=\phi_2(1/2,1/2)=1$, $\mbox{supp}(\phi_1)=\{(x,y)^\tp\in \R^2 \setsp |x+y|\le 1, |x-y|\le 1\}$ and $\mbox{supp}(\phi_2)=[0,1]^2$.
\end{example}


\begin{thebibliography}{99}

%\bibitem{bhr93}
%C.~de Boor, K.~H\"ollig, and S.~Riemenschneider, Box splines. \emph{Applied Mathematical Sciences}, \textbf{98}. Springer-Verlag, New York, 1993. xviii+200 pp.


\bibitem{cdm91}
A.~Cavaretta, W.~Dahmen and C.~ Micchelli, Stationary subdivision, \emph{Memoirs AMS}.
Amer. Math. Soc. \textbf{93} (1991).

\bibitem{ccs05}
M.~Charina, C.~Conti, and T.~Sauer, Regularity of multivariate vector subdivision schemes. \emph{Numer. Algorithms} \textbf{39} (2005),  97--113.

\bibitem{cj03}
C.~Chui and Q.~Jiang, Surface subdivision schemes generated by refinable bivariate spline function vectors. \emph{Appl. Comput. Harmon. Anal.} \textbf{15} (2003), 147--162.

\bibitem{ch19}
C.~Conti and S.~H\"uning, An algebraic approach to polynomial reproduction of Hermite subdivision schemes. \emph{J. Comput. Appl. Math.} \textbf{349} (2019), 302--315.

\bibitem{cmr14}
C.~Conti, J.-L.~Merrien and L.~Romani, Dual Hermite subdivision schemes of de Rham-type. \emph{BIT} \textbf{54} (2014), 955--977.

\bibitem{ccs16}
C.~Conti, M.~Cotronei, and T.~Sauer, Factorization of Hermite subdivision operators preserving exponentials and polynomials. \emph{Adv. Comput. Math.} \textbf{42} (2016), 1055--1079.

\bibitem{cmss19}
M.~Cotronei, C.~Moosm\"uller, T.~Sauer and N.~Sissouno, Level-dependent interpolatory Hermite subdivision schemes and wavelets. \emph{Constr. Approx.} \textbf{50} (2019), 341--366.

\bibitem{dah97}
W.~Dahmen, Wavelet and multiscale methods for operator equations. Acta Numer., \textbf{6} (1997), 55--228.

\bibitem{dhjk00}
W.~Dahmen, B.~Han, R.-Q.~Jia, and A.~Kunoth, Biorthogonal multiwavelets on the interval: cubic Hermite splines. \emph{Constr. Approx.} \textbf{16} (2000), 221--259.

\bibitem{dhmm05}
S.~Dubuc, B.~Han, J.-L.~Merrien and Q.~Mo,
Dyadic $C^2$ Hermite interpolation on a square mesh. \emph{Comput. Aided Geom. Design} \textbf{22} (2005), 727--752.


\bibitem{dm06}
S.~Dubuc and J.-L.~Merrien, Convergent vector and Hermite subdivision schemes. \emph{Constr. Approx.} \textbf{23} (2006), 1--22.

\bibitem{dm09}
S.~Dubuc and J.-L.~Merrien, Hermite subdivision schemes and Taylor polynomials. \emph{Constr. Approx.} \textbf{29} (2009), 219--245.


\bibitem{dl95}
N.~Dyn and D.~Levin, Analysis of Hermite-type subdivision schemes. Approximation theory VIII, Vol. 2, 117--124, Ser. Approx. Decompos., 6, World Sci. Publ., River Edge, NJ. 1995.


\bibitem{dl02}
N.~Dyn and D.~Levin, Subdivision schemes in geometric modelling. \emph{Acta Numer.} \textbf{11} (2002), 73--144.

%\bibitem{fhs22}
%L.~Fang, B.~Han, and Y.~Shen, Quasi-interpolating bivariate dual $\sqrt{2}$-subdivision using 1D stencils. \emph{Comput. Aided Geom. Design} \textbf{98} (2022), Paper No. 102139, 18 pp.

\bibitem{han01}
B.~Han, Approximation properties and construction of Hermite interpolants and biorthogonal multiwavelets. \emph{J. Approx. Theory} \textbf{110} (2001), 18--53.

\bibitem{han03}
B.~Han, Vector cascade algorithms and refinable function vectors in Sobolev spaces. \emph{J. Approx. Theory} \textbf{124} (2003), 44--88.

\bibitem{han03smaa}
B.~Han, Computing the smoothness exponent of a symmetric multivariate refinable function. \emph{SIAM J. Matrix Anal. Appl.} \textbf{24} (2003), 693--714.

\bibitem{han06}
B.~Han, Solutions in Sobolev spaces of vector refinement equations with a general dilation matrix. \emph{Adv. Comput. Math.} \textbf{24} (2006), 375--403.

\bibitem{han09}
B.~Han, Dual multiwavelet frames with high balancing order and compact fast frame transform. \emph{Appl. Comput. Harmon. Anal.} \textbf{26} (2009), 14--42.

\bibitem{han10}
B.~Han, Symmetric orthonormal complex wavelets with masks of arbitrarily high linear-phase moments and sum rules. \emph{Adv. Comput. Math.} \textbf{32} (2010), 209--237.

\bibitem{han10mc}
B.~Han, The structure of balanced multivariate biorthogonal multiwavelets and dual multiframelets, \emph{Math. Comp.}, \textbf{79} (2010), 917--951.

\bibitem{han13}
B.~Han, Properties of discrete framelet transforms. \emph{Math. Model. Nat. Phenom.} \textbf{8} (2013), 18--47.

\bibitem{hanbook}
B.~Han, Framelets and wavelets: Algorithms, analysis, and applications. \emph{Applied and Numerical Harmonic Analysis}. Birkh\"auser/Springer, Cham. xxxiii + 724 pp. 2017.

\bibitem{han20}
B.~Han, Analysis and convergence of Hermite subdivision schemes, \emph{Found. Comput. Math.}, (2021), published online.

\bibitem{hkz09}
B.~Han, S.~Kwon and X.~Zhuang, Generalized interpolating refinable function vectors. \emph{J. Comput. Appl. Math.} \textbf{227} (2009), 254--270.

\bibitem{hj06}
B.~Han and R.-Q.~Jia, Optimal $C^2$ two-dimensional interpolatory ternary subdivision schemes with two-ring stencils, \emph{Math. Comp.}, \textbf{75} (2006), 1287--1308.

\bibitem{hj98}
B.~Han and R.-Q.~Jia, Multivariate refinement equations and convergence of subdivision schemes. \emph{SIAM J. Math. Anal.} \textbf{29} (1998), 1177--1199.

\bibitem{hm21}
B.~Han and M.~Michelle, Wavelets on intervals derived from arbitrary compactly supported biorthogonal multiwavelets. \emph{Appl. Comput. Harmon. Anal.} \textbf{53} (2021), 270--331.

\bibitem{hm07}
B.~Han and Q.~Mo, Analysis of optimal bivariate symmetric refinable Hermite interpolants. \emph{Commun. Pure Appl. Anal.} \textbf{6} (2007), 689--718.

\bibitem{hy06}
B.~Han and T.~Yu, Face-based Hermite subdivision schemes. \emph{J. Concr. Appl. Math.} \textbf{4} (2006), 435--450.

\bibitem{hyp04}
B.~Han, T.~Yu and B.~Piper, Multivariate refinable Hermite interpolant. \emph{Math. Comp.} \textbf{73} (2004), 1913--1935.

\bibitem{hyx05}
B.~Han, T.~Yu and Y.~ Xue, Noninterpolatory Hermite subdivision schemes. \emph{Math. Comp.} \textbf{74} (2005), 1345--1367.

\bibitem{hz09}
B.~ Han and X.~Zhuang, Analysis and construction of multivariate interpolating refinable function vectors. \emph{Acta Appl. Math.} \textbf{107} (2009), 143--171.

\bibitem{jy19}
B.~Jeong and J.~Yoon, Analysis of non-stationary Hermite subdivision schemes reproducing exponential polynomials. \emph{J. Comput. Appl. Math.} \textbf{349} (2019), 452--469.

\bibitem{jj03}
R.-Q.~Jia and Q.~Jiang, Spectral analysis of the transition operator and its applications to smoothness analysis of wavelets. \emph{SIAM J. Matrix Anal. Appl.} \textbf{24} (2003), 1071--1109.

%\bibitem{jjl02}
%R.-Q.~Jia, Q.~Jiang and S.~Lee, Convergence of cascade algorithms in Sobolev spaces and integrals of wavelets. \emph{Numer. Math.} \textbf{91} (2002),  453--473.

\bibitem{jm91}
R.-Q.~Jia and C. A.~Micchelli, Using the refinement equations for the construction of pre-wavelets. II. Powers of two. \emph{Curves and surfaces} (Chamonix-Mont-Blanc, 1990), 209--246, Academic Press, Boston, MA, 1991.

\bibitem{jrz98}
R.-Q.~Jia, S.~Riemenschneider and D.-X.~Zhou, Vector subdivision schemes and multiple wavelets. \emph{Math. Comp.} \textbf{67} (1998), 1533--1563.

\bibitem{ju02}
B.~J\"uttler and U.~Schwanecke,
Analysis and design of Hermite subdivision schemes, \emph{Visual Computer} \textbf{18} (2002), 326--342.

\bibitem{mer92}
J.-L.~Merrien, A family of Hermite interpolants by bisection algorithms. \emph{Numer. Algorithms} \textbf{2} (1992), 187--200.

\bibitem{ms11}
J.-L.~Merrien and T.~Sauer, A generalized Taylor factorization for Hermite subdivision schemes. \emph{J. Comput. Appl. Math.} \textbf{236} (2011), 565--574.

\bibitem{ms17}
J.-L.~Merrien and T.~Sauer, Extended Hermite subdivision schemes. \emph{J. Comput. Appl. Math.} \textbf{317} (2017), 343--361.

\bibitem{ms19}
Merrien J.-L. \& Sauer T. (2019) Generalized Taylor operators and polynomial chains for Hermite subdivision schemes. \emph{Numer. Math.} \textbf{142}, 167--203.


\bibitem{ms98}
C.~Micchelli and T.~Sauer, On vector subdivision. \emph{Math. Z.} \textbf{229} (1998), 621--674.

\bibitem{md19}
C.~Moosm\"uller and N.~Dyn, Increasing the smoothness of vector and Hermite subdivision schemes. \emph{IMA J. Numer. Anal.} \textbf{39} (2019), 579--606.

\bibitem{mhc20}
C.~Moosm\"uller, S.~H\"uning and C.~Conti, Stirling numbers and Gregory coefficients for the factorization of Hermite subdivision operators, \emph{IMA J. Numer. Anal.}, \textbf{41} (2021), 2936--2961.


\bibitem{zhou00}
D.-X.~Zhou, Multiple refinable Hermite interpolants. \emph{J. Approx. Theory} \textbf{102} (2000), 46--71.
\end{thebibliography}
\end{document}